\documentclass{article}
\usepackage[utf8]{inputenc}
\usepackage{amsmath}
\usepackage{amsfonts}
\usepackage{amssymb}
\usepackage[linesnumbered,ruled,vlined]{algorithm2e}
\DontPrintSemicolon 
\usepackage{multirow}
\usepackage{float}
\usepackage{array}
\usepackage[left=0.75in,right=0.75in,top=0.75in,bottom=0.75in]{geometry}
\usepackage{comment}
\usepackage{graphicx}
\usepackage{url}
\usepackage{natbib}
\usepackage{color}
\usepackage{pgfplots}
\usepackage{caption}
\usepackage{subcaption}
\usepackage{pdflscape}
\usepackage{booktabs} 
\usepackage{enumerate}
\usepackage{authblk}
\usepackage{blkarray}
\usepackage{bm}
\usepackage{marvosym} 
\usepackage{soul}
\usepackage{textcomp, xspace}
\usetikzlibrary{matrix}
\usepackage{ifthen}
\usepackage{wrapfig}
\usepackage{booktabs}
\usepackage{eurosym}
\usepackage{lineno}
\usepackage{multirow}
\usepackage{lipsum}  

\SetCommentSty{mycommfont}
\definecolor{ForestGreen}{RGB}{34, 139, 34}
\usepackage{tikz}

\usepgfplotslibrary{dateplot}
\usepackage[colorlinks]{hyperref}
\hypersetup{bookmarks=true, citecolor=blue, urlcolor=blue, linkcolor=blue}

\definecolor{greenxllite}{RGB}{196,215,155}
\definecolor{bluexllite}{RGB}{149,179,215}
\definecolor{redxllite}{RGB}{218,150,148}
\definecolor{greenxl}{RGB}{155, 187, 89}
\definecolor{bluexl}{RGB}{79,129,189}
\definecolor{redxl}{RGB}{192, 80, 77}
\definecolor{OliveGreen}{RGB}{70,136,52}

\usetikzlibrary{patterns}
\allowdisplaybreaks

\newcommand{\Keywords}[1]{\par\noindent
{\small{\em \textbf{Keywords}\/}: #1}}

\usepackage{varwidth}

\title{On the Impact of Co-Optimizing Station Locations, Trip Assignment, and Charging Schedules for Electric Buses}

\author[1] {Rito Brata Nath*}
\author[1,2] {Tarun Rambha*}
\author[3] {Maximilian Schiffer}

\affil[1]{\small Department of Civil Engineering, Indian Institute of Science (IISc), Bengaluru, India}
\affil[2]{\small Center for infrastructure, Sustainable Transportation and Urban Planning (C\textit{i}STUP), Indian Institute of Science (IISc), Bengaluru, India}  
\affil[3]{\small School of Management \& Munich Data Science Institute, Technical University of Munich (TUM), Munich, Germany}

\date{ }

\newcommand*{\electricityVar}{w}
\newcommand*{\busIndex}{b}
\newcommand*{\locationIndex}{s}
\newcommand*{\locationAnotherIndex}{s'}
\newcommand*{\chargingIndex}{k}
\newcommand*{\levelVar}{l}
\newcommand*{\consumption}{e}
\newcommand*{\deadheadconsumption}{e}
\newcommand*{\electricityPrice}{c^{elec}}

\newcommand*{\capacityPrice}{c^{cap}}

\newcommand*{\stopSet}{S}
\newcommand*{\tripSet}{I}
\newcommand*{\routeSet}{R}
\newcommand*{\candidateLocationSet}{S^{cand}}
\newcommand*{\tripIndexI}{i}
\newcommand*{\tripIndexJ}{j}
\newcommand*{\tripIndexK}{k}
\newcommand*{\EndTime}{\beta}
\newcommand*{\StartTime}{\alpha}
\newcommand*{\deadheadDur}{\gamma}
\newcommand*{\idleTime}{\theta}
\newcommand*{\graph}{G}
\newcommand*{\nodeSet}{N}
\newcommand*{\arcSet}{A}
\newcommand*{\compatibleSet}{A^{comp}}
\newcommand*{\depotSet}{D}
\newcommand*{\numTrips}{n}
\newcommand*{\numDepots}{m}
\newcommand*{\OriginStopsSet}{O}
\newcommand*{\DestinationStopsSet}{E}
\newcommand*{\TripStartStop}{i^{start}}
\newcommand*{\TripEndStop}{i^{end}}
\newcommand*{\BusSet}{B}
\newcommand*{\arcVar}[3]{x_{#1,#2,#3}}
\newcommand*{\LocVar}{z}
\newcommand*{\contractedVar}{q}
\newcommand*{\LocationSet}{Z}
\newcommand*{\EndChargeVar}[1]{y_{#1}^{end}}
\newcommand*{\StartChargeVar}[1]{y_{#1}^{start}}
\newcommand*{\EndChargeLocVar}[2]{r_{#1, #2}^{end}}
\newcommand*{\StartChargeLocVar}[2]{r_{#1, #2}^{start}}
\newcommand*{\EndChargeDurVar}[2]{t_{#1,#2}^{end}}
\newcommand*{\StartChargeDurVar}[2]{t_{#1,#2}^{start}}
\newcommand*{\EndLevelVar}[2]{l_{#1,#2}^{end}}
\newcommand*{\StartLevelVar}[2]{l_{#1,#2}^{start}}
\newcommand*{\energyTrip}[1]{d_{#1}}
\newcommand*{\energyDeadhead}[2]{d_{#1,#2}}

\newcommand*{\BigMConst}{M}
\newcommand*{\BusCost}{c^{bus}}
\newcommand*{\LocationCost}{c^{loc}}
\newcommand*{\DeadheadCost}[2]{c_{#1,#2}^{km}}
\newcommand*{\chargingRate}{\lambda}
\newcommand*{\BatteryCap}{l_{max}}
\newcommand*{\MinEnergy}{l_{min}}
\newcommand*{\timeIndex}{t}
\newcommand*{\timeAnotherIndex}{t'}
\newcommand*{\OppSet}{K_b}
\newcommand*{\SingleOppSet}{K_b^{1}}
\newcommand*{\DualOppSet}{K_b^{2}}
\newcommand*{\timeSet}{T}
\newcommand*{\StartStopTimeSet}{T^{start}}
\newcommand*{\EndStopTimeSet}{T^{end}}
\newcommand*{\timeperiodSet}{P}
\newcommand*{\periodIndex}{p}
\newcommand*{\startLoc}{s^{start}}
\newcommand*{\EndLoc}{s^{end}}
\newcommand*{\maxTransfer}{\psi_b}
\newcommand*{\electricityStartVar}[3]{w^{start}_{#1,#2,#3}}
\newcommand*{\electricityEndVar}[3]{w^{end}_{#1,#2,#3}}
\newcommand*{\EndChargeTimeVar}[2]{y^{end}_{#1,#2}}
\newcommand*{\StartChargeTimeVar}[2]{y^{start}_{#1,#2}}
\newcommand*{\busrotationList}{V}
\newcommand*{\busrotationAnotherList}{V^{temp}}
\newcommand*{\bestbusrotationList}{V^{best}}

\newcommand*{\ExchangetriprotationList}{V^{exst}}
\newcommand*{\ExchangedepotrotationList}{V^{exd}}

\newcommand*{\ShifttriprotationList}{V^{sst}}

\newcommand*{\busrotationb}{V_b}

\newcommand*{\numTripsBusb}{n_b}
\newcommand*{\numTripsBus}[1]{n_{#1}}

\newcommand*{\bustrip}[2]{i_{#1}^{#2}}
\newcommand*{\nearestdepotTrip}[1]{\Delta_{#1}}
\newcommand*{\objectivefunction}[2]{f(#1,#2)}

\newcommand*{\LocationAnotherSet}{Z^{temp}}

\newcommand*{\currentUtilization}[1]{\kappa_{#1}}
\newcommand*{\potentialUtilization}[1]{\pi_{#1}}
\newcommand*{\currentChargeCapacity}[1]{q^*_{#1}}

\newcommand{\fewThreshold}{\zeta}

\newcommand*{\PreviousTimeSet}[1]{T^{prev}_{#1}}
\newcommand*{\NextTimeSet}[1]{T^{next}_{#1}}

\begin{document}
\maketitle
\let\thefootnote\relax\footnotetext{* These authors contributed equally to this manuscript.}
\vspace{-5mm}
\begin{abstract}
As many public transportation systems around the world transition to electric buses, the planning and operation of fleets can be improved via tailored decision-support tools. In this work, we study the impact of jointly locating charging facilities, assigning electric buses to trips, and determining when and where to charge the buses. We propose a mixed integer linear program that co-optimizes planning and operational decisions jointly and an iterated local search heuristic to solve large-scale instances. Herein, we use a concurrent scheduler algorithm to generate an initial feasible solution, which serves as a starting point for our iterated local search algorithm. In the sequential case, we first optimize trip assignments and charging locations. Charging schedules are then determined after fixing the optimal decisions from the first level. The joint model, on the other hand, integrates charge scheduling within the local search procedure. The solution quality of the joint and sequential iterated local search models are compared for multiple real-world bus transit networks. Our results demonstrate that joint models can help further improve operating costs by 14.1\% and lower total costs by about 4.1\% on average compared with sequential models. In addition, energy consumption costs and contracted power capacity costs have been reduced significantly due to our integrated planning approach. 

\vspace{3mm}
\end{abstract}

\section{Introduction}
\label{sec:intro}
Electric Vehicles (EVs) have the potential to significantly reduce the environmental impacts of tail-pipe emissions. While several government subsidies promote electric mobility \citep{sierzchula2014influence}, most are personal-vehicle centric and hence do not address congestion issues \citep{wan2015china}. For cities to grow sustainably, it is necessary to electrify public transportation systems \citep{pelletier2019electric} that offer eco-friendly travel alternatives and at the same time limit the use of personal travel modes. However, large-scale deployment of electric buses requires proper charging infrastructure and strategies for day-to-day operations. To this end, one must factor in the buses' state of charge (SoC) at different time points and the energy needs of transit trips to optimize operational and planning decisions. 

In this context, the range of EV buses can vary widely depending on the battery capacity. For example, the BYD's K7M buses have an advertised range of 258 km, and their K8M models have a higher range of 315 km \citep{BYD}. The true range is usually lower and depends on several other factors, such as driving styles, passenger loads, and weather conditions. Unlike diesel/gasoline buses, electric buses may require charging between operations to replenish their battery levels due to limited range. Furthermore, it is recommended not to discharge the battery below a certain threshold for handling unforeseen events and improving the battery life cycle \citep{ellingsen2022life}.

\begin{figure}[H]
	\centering
	\includegraphics[scale=0.48]{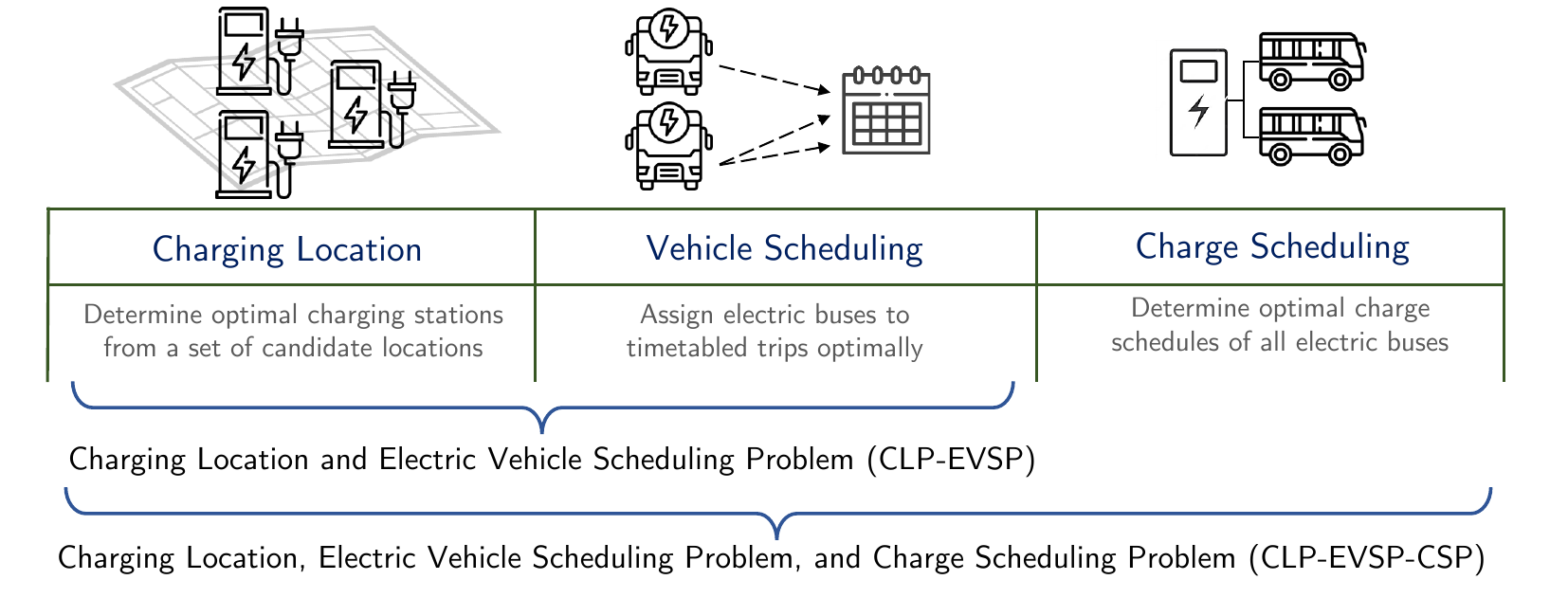}
	\caption{Different tasks involved in strategic and operational planning of electric bus fleets}
    \label{fig:infographic_problems}
\end{figure}

If the layover times between trips are short, buses may not have sufficient time to recharge their batteries. This issue can be addressed by understanding where to locate charging stations -- the \textit{Charging Location Problem} (CLP), how to assign buses to trips -- the \textit{Vehicle Scheduling Problem} (VSP), and when and how much to charge buses at a charging station -- the \textit{Charge Scheduling Problem} (CSP). Figure \ref{fig:infographic_problems} illustrates these problem components. While the decisions to these problems can be made sequentially, strategic choices can have a significant role on the operational ones \citep{dirks2022integration}. For instance, the bus-to-trip assignment is more constrained if the charging locations are pre-determined. Likewise, the bus-to-trip assignment determines the charging demand and the time slots for recharging. Consequently, a holistic framework that co-optimizes these planning tasks in an integrated fashion promises cost savings over their sequential counterparts. However, finding high-quality solutions for joint models presents significant challenges due to increased search space, raising the question: \textit{Can we design algorithms that co-optimize the aforementioned planning tasks and are computationally tractable? If so, what is the added benefit of such an integrated approach?} This paper explores solution techniques for such integrated frameworks and demonstrates their potential to achieve significant cost savings while remaining tractable.
Specifically, our paper addresses the Charging Location and Electric Vehicle Scheduling Problem (CLP-EVSP) and the Charging Location, Electric Vehicle Scheduling Problem, and Charge Scheduling Problem (CLP-EVSP-CSP) jointly. Numerous research studies have focused on these planning tasks in isolation, but there has been very limited effort in developing a holistic framework, due to the aforementioned challenges. We seek solutions that minimize the fixed costs -- total investment cost of electric buses and charging facilities -- and operational costs due to deadheading and recharging. In the charge scheduling problem, we permit recharging multiple buses simultaneously. 

A list of features modeled in this research and comparisons with relevant papers are shown in Table \ref{tab:intro_model}. Our research makes multiple contributions to the literature on the planning and operations of electric bus fleets, as summarized below: 
\begin{itemize}
\item We propose a Mixed Integer Linear Programming (MILP) formulation for the CLP-EVSP-CSP for a multiple depot setting. Our formulation captures several other practical characteristics of electric transit systems, such as grid capacity, time-of-day pricing, and split charging.
\item We design an iterated local search (ILS) procedure for the CLP-EVSP to optimize locations and bus rotations. We also propose a novel MILP and two surrogate Linear Programming (LP) formulations to solve the charge scheduling problem. These LPs allow us to extend the ILS framework to the joint CLP-EVSP-CSP while maintaining tractability.   
\item We demonstrate the practical value of our integrated planning approach using multiple real-world \href{https://gtfs.org/}{GTFS} datasets. With reasonable assumptions about input parameters, our joint CLP-EVSP-CSP model exhibited average savings of 41.4\% in scheduling costs, 39.6\% in operational costs, and 14.4\% in total costs over all the test networks compared to the concurrent scheduler. The corresponding average savings compared to the sequential CLP-EVSP model were 17.5\%, 14.1\%, and 4.1\%, respectively.
\end{itemize} 

The rest of the paper is structured as follows: Section \ref{sec:litrev} explores relevant literature on the CLP, EVSP, and CSP, including a few studies on joint modeling. Section \ref{sec:formulation} formally defines the CLP-EVSP-CSP problem. In Section \ref{sec:csp_joint}, we present CSP models and a joint formulation for the CLP-EVSP-CSP. Section \ref{sec:localsearch} introduces algorithms for the CLP-EVSP and CLP-EVSP-CSP ILS heuristics. We study the benefits from applying these heuristics to different networks in Section \ref{sec:results}. Finally, Section \ref{sec:conc} summarizes the paper's findings and discusses future research directions.

 \begin{table}[t]
 \small
  \centering
  \caption{Modeling features of relevant works ($\bullet$ \, Fully modeled, \, $\circ$ \, Partially modeled). In ambiguous situations where it was not clear how a certain feature was represented, we considered it to be fully modeled.}
    \begin{tabular}{lp{11.5mm}p{11.5mm}p{11.5mm}p{11.5mm}p{11.5mm}p{11.5mm}p{11.5mm}p{11.5mm}}
    \hline
     \textbf{References} & \textbf{Bus acquisition} & \textbf{Charging station costs} & \textbf{Vehicle scheduling} & \textbf{Dynamic pricing} & \textbf{Multiple depots} & \textbf{Power load costs} & \textbf{Partial charging} & \textbf{Split charging}\\
    \hline
    \cite{rogge2018electric} & $\bullet$ & $\circ$ & $\bullet$ & & & $\bullet$ & &\\
    \cite{janovec2019exact} & $\bullet$ & & $\bullet$ & & & & $\bullet$ &\\
    \cite{li2019mixed} & $\bullet$ & $\bullet$ & $\bullet$ & & $\bullet$ & & &\\
    \cite{liu2019regional} & $\bullet$ & $\bullet$ & $\bullet$ & & $\bullet$ & & $\bullet$ &\\
    \cite{abdelwahed2020evaluating} & & & & $\bullet$ & & & $\bullet$ & $\bullet$\\
    \cite{li2020joint} & $\bullet$ & $\circ$ & $\bullet$ & $\bullet$ & $\bullet$ & & $\bullet$ &\\
    \cite{liu2020battery} &  & $\circ$ & $\bullet$ & & $\bullet$ & & $\bullet$ &\\
    \cite{olsen2020scheduling} & $\bullet$ & & $\bullet$ & & & & $\bullet$ &\\
    \cite{teng2020integrated} & $\bullet$ & & $\bullet$ & $\bullet$ & & & &\\
    \cite{yao2020optimization} & $\bullet$ &  & $\bullet$ & & $\bullet$ & & &\\
    \cite{zhou2020bi} & & & $\bullet$ & $\bullet$ & & & $\bullet$ &\\
    \cite{alvo2021exact} & $\circ$ & & $\bullet$ & & & $\bullet$ & $\bullet$ &\\
    \cite{alwesabi2021novel} & $\bullet$ & $\bullet$ & & & & $\circ$ & $\bullet$ &\\
    \cite{jiang2021multi} & $\bullet$ & & $\bullet$ & $\bullet$ & $\bullet$ & & $\bullet$ &\\
    \cite{lee2021optimal} & $\bullet$ & $\bullet$ & $\bullet$ & & & $\bullet$ & $\bullet$ &\\
    \cite{stumpe2021study} & $\bullet$ & $\bullet$ & $\bullet$ & & & & $\bullet$ & \\
    \cite{olsen2022location} & $\bullet$ & $\bullet$ & $\bullet$ & & & & $\bullet$ &\\
    \cite{zhang2022long} & $\bullet$ & $\bullet$ & $\bullet$ & & & & $\bullet$ &\\
    \cite{gairola2023optimization} & $\bullet$ & $\bullet$ & $\circ$ & $\bullet$ & $\circ$ & $\circ$ & $\bullet$ &\\
      \hline
     \textbf{Our work} & $\bullet$ & $\bullet$ & $\bullet$ & $\bullet$ & $\bullet$ & $\bullet$ & $\bullet$ & $\bullet$ \\
    \hline 
    \end{tabular}%
  \label{tab:intro_model}%

\end{table}%

\section{Literature Review}
\label{sec:litrev}
Strategic and operational problems in managing electric buses and charging stations have been widely researched in recent years. Strategic problems include locating charging and battery swapping stations and fleet transition decisions. On the other hand, operational problems deal with electric bus scheduling, charging schedule optimization, and battery swapping. A detailed survey on these problems can be found in \cite{perumal2022electric}. This section discusses the state-of-the-art models for CLP, EVSP, and CSP of the joint framework proposed in this paper. 

\subsection{Optimizing Charging Locations}
Locating EV charging stations resembles the classic facility location problem \citep{schettini2023locating}. The objective of the CLP is to establish charging sites from a set of candidate locations to fulfill the demand of electric vehicles. For transit networks, charging stations are typically located at the start and end of bus routes, areas where buses tend to have extended parking times. Station capacities are also usually optimized along with location decisions. E.g., a new charging facility is redundant if a nearby station can serve the peak demand \citep{metais2022too}. Capacities, often determined by the number of chargers, are constrained by available parking spaces. Some studies ignore this effect and assume infinite capacities \citep{wang2010locating, he2015deploying}, while others associate station capacities to power grid limits \citep{zhang2016pev}. CLP formulations can optimize objectives such as distance traveled by EVs \citep{shahraki2015optimal}, number of EVs charged \citep{cavadas2015mip}, and deadheading distances \citep{xu2013optimal}. For electric buses, specific objectives include energy consumption of the system \citep{xylia2017locating}, number and locations of chargers \citep{kunith2017electrification}, and infrastructure cost \citep{he2019fast}. Most of these problems are solved using MILP models \citep{alwesabi2020electric} and meta-heuristics \citep{iliopoulou2019integrated}. 

\subsection{Electric Vehicle Scheduling Problem}
The VSP assigns a set of timetabled trips to vehicles originating from one or more depots based on time compatibility \citep{bunte2009overview}. The multi-depot VSP (MDVSP) can be formulated either as a \textit{multi-commodity flow model} \citep{forbes1994exact, kliewer2006time} or as a \textit{set partitioning model} \citep{ribeiro1994column} and is known to be an NP-hard problem \citep{bertossi1987some}. 
The EVSP is an extension of the VSP with additional range constraints \citep{reuer2015electric, liu2020battery, olsen2020scheduling}. EVSP objectives are similar to those of VSP, but they can also include transit-specific aims, such as minimizing the total cost of ownership \citep{rogge2018electric}. Most studies in the literature differ in solution techniques and problem scales. Exact methods using MILP models are common \citep{li2014transit, janovec2019exact, li2019mixed}. The Column Generation (CG) method can also be applied to the EVSP for improving tractability \citep{adler2017vehicle, tang2019robust}. However, large problem instances require heuristics or meta-heuristics such as GAs \citep{rogge2018electric, liu2019regional, li2020joint, yao2020optimization}, particle swarm optimization \citep{teng2020integrated}, simulated annealing \citep{zhou2020collaborative}, and ant colony optimization \citep{wang2007heuristic}. Heuristics based on constraint generation \citep{haghani2002heuristic}, concurrent scheduling \citep{adler2017vehicle}, adaptive large neighborhood search \citep{wen2016adaptive, perumal2021solution}, and iterative neighborhood search \citep{zhou2020collaborative} have also proven to scale successfully to larger instances.  

\subsection{Charge Scheduling Problem}
The stand-alone version of the CSP assumes that the allocation of buses to their respective trips and the charging stations' locations are predetermined. The goal of the CSP is usually to minimize charging costs under time-varying electricity prices. While charging schedules are typically designed not to exceed grid capacity \citep{zhang2016pev}, some studies allow grid reinforcement and exclude such constraints \citep{sadeghi2014optimal}. Among literature that considered dynamic electricity pricing, \cite{leou2017optimal} presented two MILP models for a single charging station to minimize capacity and energy charges. \cite{abdelwahed2020evaluating} proposed discrete time- and discrete event-based MILP models to minimize charging costs. Although they considered split charging events, they did not extend their model to multiple charging locations. A similar approach to minimize total charging costs was adopted by \cite{he2020optimal} using a linearized model that could be solved using commercial solvers. \cite{ke2020battery} employed a GA to minimize electricity costs. They also considered selling electricity back to the power company. 

Several other studies have focused on power load management at charging stations. \cite{jahic2019charging} proposed a greedy algorithm and a heuristic to minimize the peak load at a central depot. \cite{houbbadi2019optimal} analyzed overnight depot charging, considering battery aging, using nonlinear programming. The operational cost of power distribution systems was minimized by \cite{bagherinezhad2020spatio} using a relaxed cone programming model. Lastly, \cite{yang2017charging} focused on minimizing electricity consumption costs in a wireless charging system.

\subsection{Joint Models}
The decision variables in the problems discussed above are interconnected. Several papers have jointly studied CLP-EVSP and EVSP-CSP to optimize electric bus operations fully. These studies vary in their objective functions, methodologies, assumptions about charging profiles, and the sizes of test networks. This section offers a discussion on the joint modeling of CLP-EVSP and EVSP-CSP. Since the solution to the CSP inherently includes the solution for the CLP, we do not discuss integrated CLP-CSP models separately. Table \ref{tab:lit_review_table} provides an overview of the papers addressing the CLP-EVSP.

\begin{table}[h]
 \small
  \centering
  \caption{Overview of different networks and methods used to solve the CLP-EVSP. (Data on the largest instance used is shown in the table. Missing information has been marked as `-')}
    \begin{tabular}{lp{9mm}p{9mm}p{9mm}p{16mm}p{58mm}}
    \hline
    \textbf{Reference} & \textbf{Trips} & \textbf{Routes} & \textbf{Stops} & \textbf{Candidate Locations} & \textbf{Solution Techniques} \\
    \hline
    \cite{rogge2018electric} & 200 & 3 & - & - & Grouping GA \\
    \cite{li2019mixed} & 288 & 6 & - & 2 & MILP model \\
    \cite{liu2019regional} & 544 & 4 & - & 2 & Bin-packing and GA \\
    \cite{li2020joint} & 867 & 8 & 164 & 5 & Adaptive GA \\ 
     \cite{yao2020optimization} & 931 & 4 & - & 2 & GA with a fitness computation algorithm  \\
    \cite{alwesabi2021novel} & 102 & 6 & - & 57 & MILP model\\
    \cite{lee2021optimal} & 74 & 1 & 38 & 38 & A two-stage nonlinear integer model\\
    \cite{stumpe2021study} & 1296 & - & - & 88  & MILP and VNS\\
    \cite{hu2022joint} & 213 & 3 & 111 & 111 & MILP model\\
    \cite{olsen2022location} & 3067 & - & 209 & 209 & VNS algorithm\\
    \cite{zhang2022long} & - & 16 & 238 & 7 & MILP model\\
    \hline
    \end{tabular}%
  \label{tab:lit_review_table}%
\end{table}%

\cite{stumpe2021study} and \cite{olsen2022location} introduced MILP models for the CLP-EVSP. For larger instances, they used a variable neighborhood search (VNS) algorithm, originally proposed by \cite{mladenovic1997variable} and  \cite{hansen2010variable}. \cite{alwesabi2021novel} integrated vehicle scheduling with charging planning using dynamic wireless charging. \cite{li2020joint} used a GA variant for the Stationary Charger Deployment and MDVSP, with partial charging and time-of-day electricity pricing. 

Both in-depot and en-route charging can replenish electric buses' battery levels \citep{zhou2020bi}. A few studies optimize charge schedules and bus-to-trip assignments for a given set of charging stations. \cite{jiang2021multi} optimized the MDVSP along with charging events while allowing partial charging. A linearized version of a non-convex programming model was used by \cite{zhou2022electric} to address the bus and charge scheduling for a single route. \cite{lee2021optimal} studied joint bus and charge scheduling for small networks using a two-stage exact approach that considered electricity consumption costs. Similar problems were investigated with fixed electricity pricing \citep{janovec2019exact, olsen2020scheduling, gkiotsalitis2023exact} and time-varying costs \citep{teng2020integrated, zhou2020collaborative, jiang2021multi, klein2023electric}. Integrated approaches that jointly optimize charger configurations and vehicle routing decisions have also shown significant benefits in logistics \citep{schiffer2019vehicle, schiffer2021integrated}. Yet, the application of such integrated strategies in electrification of bus networks remains less explored.

 In conclusion, many CLP-EVSP models overlook charge scheduling, which can yield additional cost savings. EVSP-CSP, on the other hand, assumes known charging stations whose locations could have been sub-optimal. The problem of integrating CLP-EVSP with CSP has not been widely researched, primarily because of its computational challenges. Our joint CLP-EVSP-CSP addresses this gap and is designed for tractability for network sizes commonly found in the literature. We also improve upon simplifying assumptions frequently found in CSP literature by considering multiple charging locations and power load capacities. 

\section{Problem Description}
\label{sec:formulation}
EV bus fleet operations vary across transit agencies. Buses typically start their journey from a depot and return to the same or a different depot at the end of the day. Two types of charging models are common: overnight-slow and opportunity-fast charging. Overnight charging is usually done at depots, and buses start with a fully charged battery the next day. Opportunity charging can be performed at any charging location during the day when the battery levels are low or in response to dynamic prices and availability of charging opportunities. Our model focuses on opportunity charging and, in this context, makes the following assumptions.
\begin{itemize}
\itemsep 0pt
    \item The candidate charging stations are chosen from terminal bus stops of routes where most buses stay for a significant duration. We make this assumption because charging buses at intermediate bus stops adds to passenger delays \citep{iliopoulou2021robust}.
    \item All buses are electric and homogeneous with the same range and can perform any scheduled trip. Buses can charge multiple times during their layovers with linear charging rates. We disallow deadheading to a charging station not in a bus itinerary, but allow partial charging. 
    \item Due to the challenge of obtaining accurate information on the locations and capacities of depots, we assume a predetermined set of depots. We chose these depots strategically based on the major terminal stops, which serve as starting or ending points for multiple routes.
    \item Interlining is permitted, i.e., a bus can serve trips along multiple routes. A bus does not need to return to its starting depot at the end of the day. However, since schedules are periodic over different days, we require the initial and final distribution of buses to be the same across depots.
    \item Charging costs vary by the time of the day. Hence, it may benefit to fully charge a bus during an off-peak period or wait when the electricity prices are high. We assume all stations are equipped with smart charging technology that dynamically regulates the amount of energy transferred to the buses \citep{sadeghian2022comprehensive}. 
    \item Overnight charging costs at the depots are ignored since they do not vary across bus-to-trip assignments; this is due to deadheading being already accounted for and prices typically remaining constant during this period.
    \item Transit schedules are deterministic. While it is also possible to re-design trip schedules \citep{tang2023optimization}, we only optimize supply components and assume that the timetabled trips (that are typically developed using transit demand estimates) are not altered.  Effects of congestion and delays due to traffic or vehicle breakdowns are out of scope of the current work. Still, they can be integrated into our methodology to some degree by adding appropriate slacks. 
\end{itemize}

\begin{figure}[H]
\centering
\begin{subfigure}{.5\textwidth}
  \centering
  \includegraphics[width=0.9\textwidth]{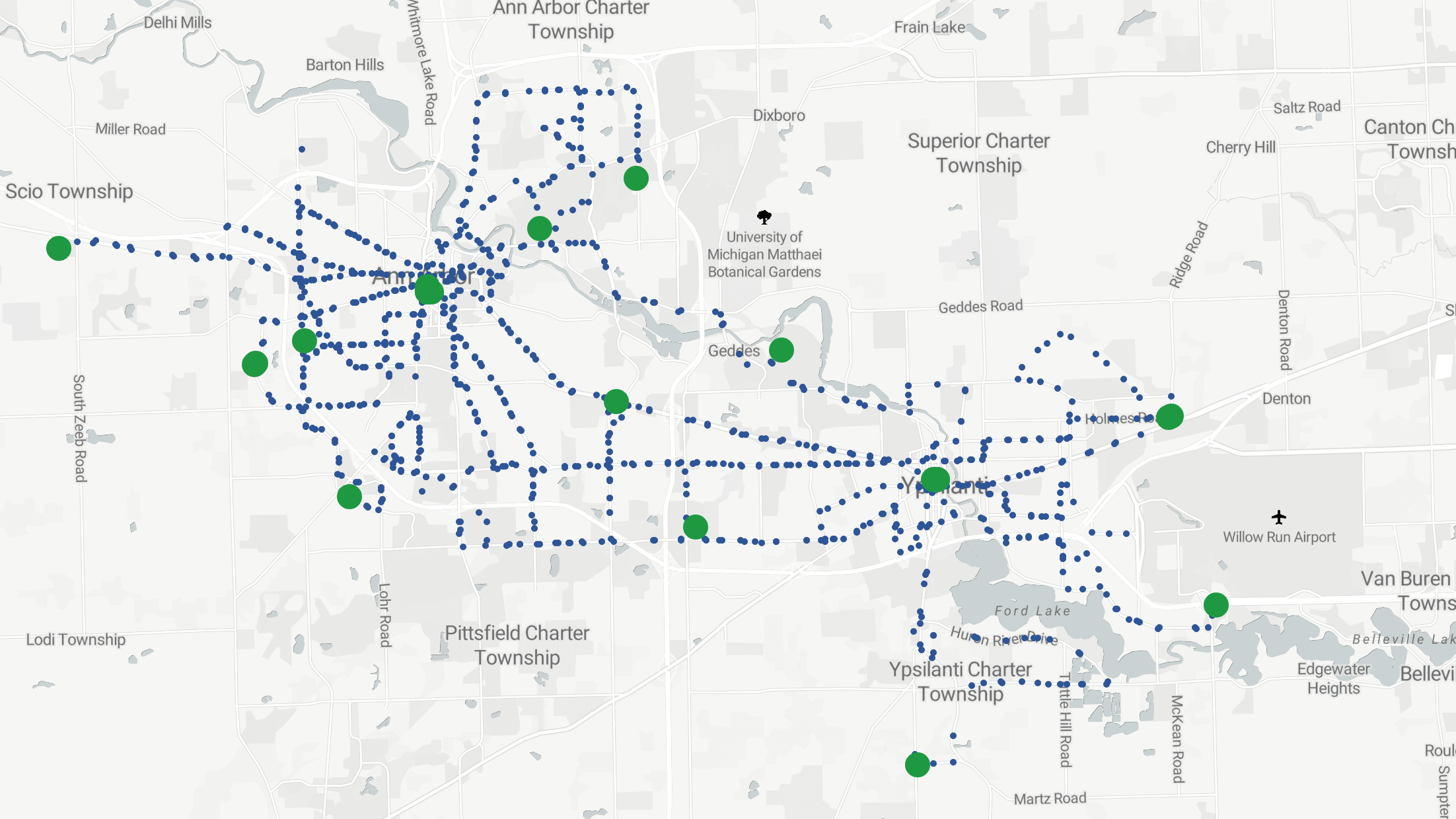}
  \caption{Bus stops and candidate charging locations}
  \label{fig:potential_locations}
\end{subfigure}%
\begin{subfigure}{.5\textwidth}
  \centering
  \includegraphics[width=0.9\textwidth]{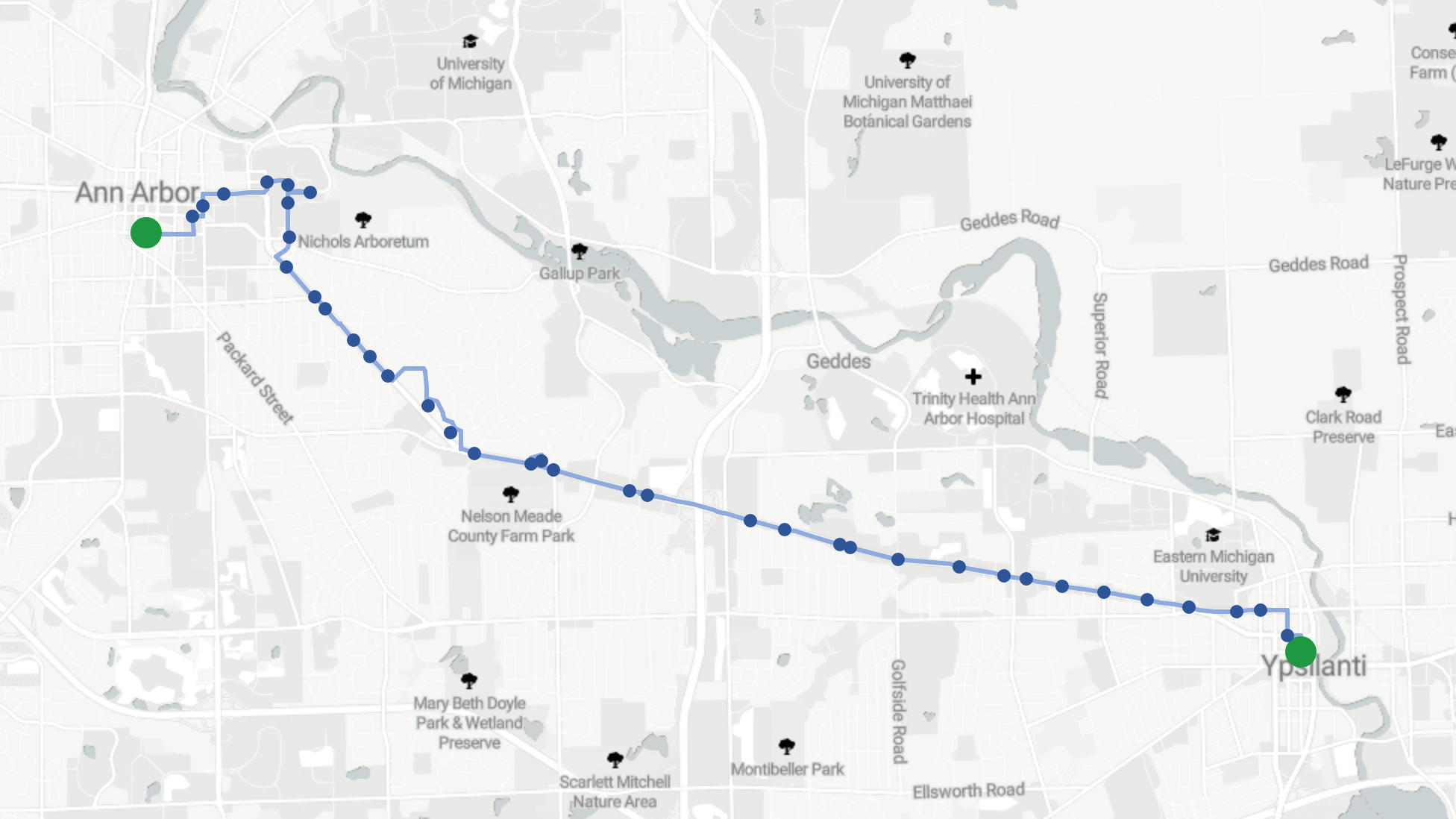}
  \caption{Washtenaw bus route}
  \label{fig:ann_route}
\end{subfigure}
\caption{Bus stops and candidate charging locations in the Ann Arbor Area Transportation Authority network, US.}
\label{fig:route_potential_locations}
\end{figure}

The transit data includes a set of stops $\stopSet$, daily trips $\tripSet$, and bus routes $\routeSet$. We assume that the candidate charging locations $\candidateLocationSet = \{\locationIndex_1, \locationIndex_2, \ldots, \locationIndex_p\}$ are the start and end locations of all routes. Figure \ref{fig:potential_locations} shows an example network. Stops in $\stopSet$ are indicated in blue, and those in $\candidateLocationSet$ are shown as larger green nodes. The bus stops along a sample route are shown in Figure \ref{fig:ann_route}.  
A pair of trips $\tripIndexI$ and $\tripIndexJ$ are \textit{compatible} if trip $\tripIndexJ$ can be carried out after trip $\tripIndexI$ by the same bus, i.e., $\EndTime_i + \deadheadDur_{\tripIndexI,\tripIndexJ} + \idleTime_{\tripIndexI,\tripIndexJ} \leq \StartTime_{\tripIndexJ}$, where $\EndTime_{\tripIndexI}$ is the end time of $\tripIndexI$, $\deadheadDur_{\tripIndexI,\tripIndexJ}$ is the deadhead trip duration from the end stop of $\tripIndexI$ to the start stop of $\tripIndexJ$, $\idleTime_{\tripIndexI,\tripIndexJ}$ is the idle time before or after the deadhead trip 
(layover time $-$ deadhead time) in minutes and rounded down to the nearest integer, and $\StartTime_{\tripIndexJ}$ is the start time of $\tripIndexJ$. The idle time is set to zero if one of the trip ends is a depot.
\begin{figure}[H]
	\centering
	\includegraphics[scale=0.55]{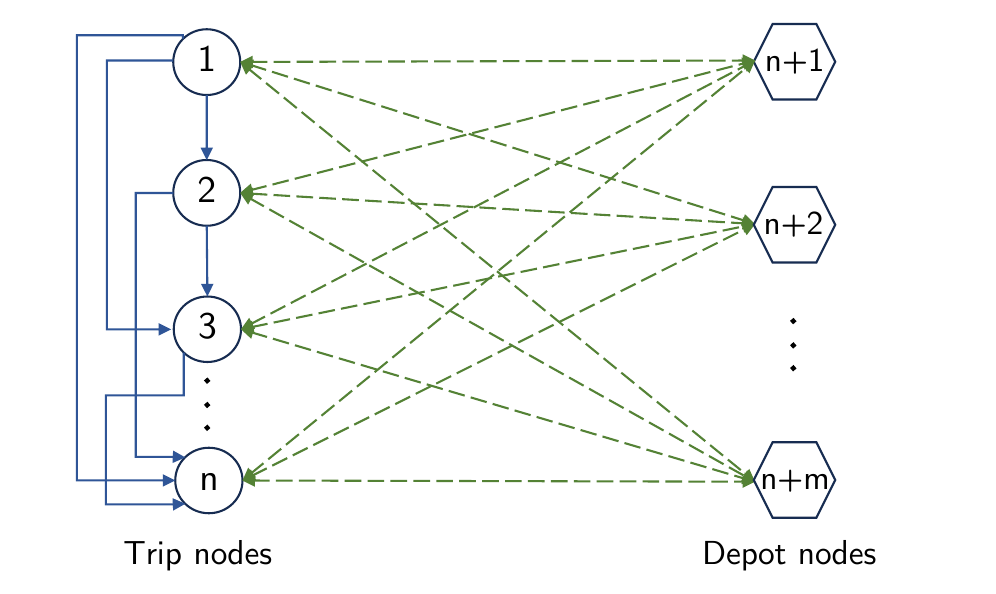}
	\caption{Network diagram for the EVSP}
\label{fig:network_diagram}
\end{figure}
Trip compatibility can be visualized using a graph $\graph = (\nodeSet, \arcSet)$, where $\nodeSet$ and $\arcSet$ are the set of nodes and arcs as shown in Figure \ref{fig:network_diagram}. We define $\nodeSet = \tripSet \cup \depotSet$, where $\tripSet = \{1, 2, \ldots, \numTrips\}$ is the set of $\numTrips$ trips and $\depotSet = \{\numTrips+1, \numTrips+2, \ldots , \numTrips+\numDepots\}$ is the set of $\numDepots$ depots. The set of compatible trip pairs, $\compatibleSet$, are denoted using solid blue arcs, and the connections between the depots and the trip nodes are shown using dashed green arcs. We do not require additional nodes for candidate charging locations since the trip itineraries capture the sequence of visits to charging station locations.

The set of all arcs $\arcSet$ can be written as $\arcSet = \compatibleSet \cup_{\tripIndexI \in \tripSet} \{(\numTrips+1, \tripIndexI), (\numTrips+2, \tripIndexI), \ldots , (\numTrips+\numDepots, \tripIndexI)\} \cup_{\tripIndexI \in \tripSet} \{(\tripIndexI, \numTrips+1), (\tripIndexI, \numTrips+2), \ldots , (\tripIndexI, \numTrips+\numDepots)\}$. Bus charging levels are tracked only at trip terminals. We denote $\OriginStopsSet_{\locationIndex}$ and $\DestinationStopsSet_{\locationIndex}$ as the set of trips with $\locationIndex$ as their start/origin and end/destination stops, respectively. Mathematically, they can be written as: $\OriginStopsSet_{\locationIndex} = \{\tripIndexI \in \tripSet: \TripStartStop = \locationIndex\}$ and $\DestinationStopsSet_{\locationIndex} = \{\tripIndexI \in \tripSet: \TripEndStop = \locationIndex\}$, where $\TripStartStop$ and $\TripEndStop$ are the starting and ending stops of trip $\tripIndexI$, respectively. We assume an initial set of buses $\BusSet = \{\busIndex_1, \busIndex_2, \ldots, \busIndex_n\}$ (not all of which may be used in the optimal solution). The maximum number of buses available can be set to the number of trips $\numTrips$ as each bus can be assigned to exactly a single trip, but tighter bounds or existing fleet data can also be used if available.

The joint CLP-EVSP-CSP model is formulated as a MILP in Section \ref{sec:cl_evsp_csp}. The objective is to minimize the fixed cost of buses and charging stations, along with operational costs encompassing deadheading, electricity consumption, and grid power capacity costs. Specifically, through bus scheduling, we aim to ensure that all the trip nodes depicted in Figure \ref{fig:network_diagram} are served exactly once. Charging stations are opened based on the charge feasibility of buses, and charge schedules are devised considering charging station availability, bus idle time, and time-of-day electricity prices. Constraints include trip compatibility, charging level feasibility, and battery level thresholds. The decision variables in the model capture the arc flows in $G$, the candidate bus stops that are converted into charging stations, and the energy transferred to buses every minute while charging. We present a formal exposition of the CSP and the joint model in Section \ref{sec:csp_joint}.

\section{Charge Scheduling and Joint Models}
\label{sec:csp_joint}

The CLP-EVSP model commonly found in the literature \citep{li2019mixed, stumpe2021study} does not consider electricity costs and grid capacity in deciding when and where buses should charge. This could increase the need for simultaneous charging and lead to higher operating costs. To tackle this issue, we first formulate a general CSP that minimizes both electricity consumption and contracted power capacity costs in Section \ref{sec:csp}. In this CSP, buses have the flexibility to charge between trips, either at the last stop of a trip, the starting stop of its next trip, or at both ends if charging stations are available. This approach necessitates a MILP framework, which is subsequently integrated into a CLP-EVSP formulation, resulting in our joint optimization model, CLP-EVSP-CSP, in Section \ref{sec:cl_evsp_csp}. However, the MILP version of the CSP is not computationally tractable for realistic instance sizes. Therefore, we introduce surrogate CSP models in Section \ref{sec:surrogate_lp}. These models prioritize where to charge between trip ends, enabling us to solve the CSP using linear programs.

\subsection{CSP Model}
\label{sec:csp}
\begin{figure}[H]
    \centering
    \includegraphics[width=0.95\textwidth]{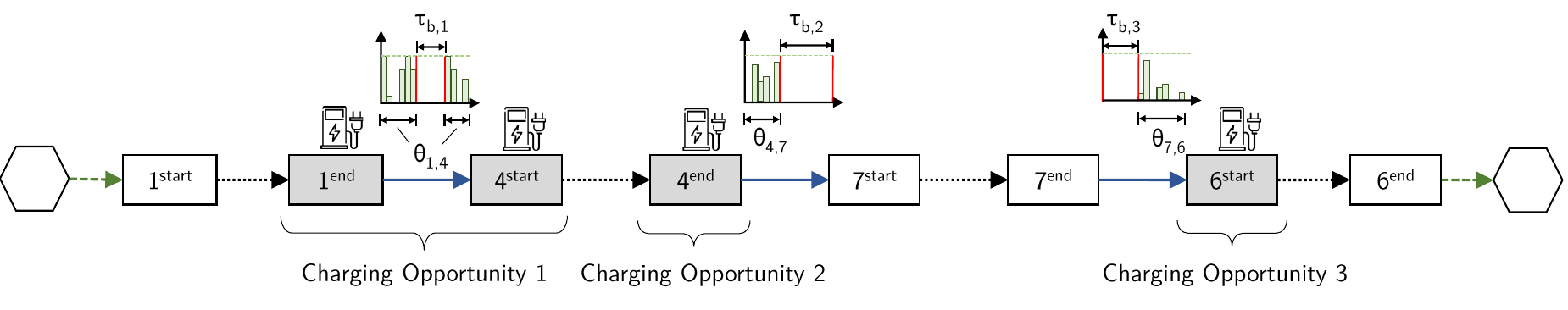}
	\caption{An illustration of the CEE charging strategy (The green dashed arrows indicate depot trips, black dotted arrows indicate service trips, and blue solid arrows indicate deadhead trips. Charging is allowed at the gray terminals.)}
\label{fig:diagram_cee}
\end{figure}

Considering time-of-day electricity prices and the contracted capacity of charging locations, a bus may choose to charge at the end stop of a trip, then deadhead to the starting stop of the next trip and charge again if a charging opportunity exists. We address this aspect through a MILP model, where the time-step at which a bus begins deadheading if charging opportunities are available at both ends is also a decision variable. We refer to this model as the \textit{Charge at Either End} (CEE) version of the CSP. We split time into one-minute intervals and consider the energy provided to buses at charging stations every minute as a decision variable. Figure \ref{fig:diagram_cee} illustrates the CEE charging. Trip ends where charging stations are located are shown using an icon, and the nodes are shaded when charging is permitted. In cases where the bus can charge at both the ending stop of a trip and the starting stop of the next trip (e.g., trips 1 and 4), the idle time $\idleTime_{1,4}$ is divided into two charging sessions, and the model determines the optimal time to begin deadheading. The solid red lines on the plots between trips indicate the start and end times of deadheading. We denote the deadheading duration for bus $\busIndex$ during its $\chargingIndex^\text{th}$ charging opportunity by $\tau_{b,k}$. The dotted green lines represent the maximum energy transferable ($\maxTransfer$)  in a single time period. Note that not charging in a particular time-step may be optimal, depending on the energy needs of other buses at the charging station. For the charging opportunities between trips (4,7) and (7,6), the process is straightforward as the charging station is only at one end of the trips. Therefore, the bus will either charge and deadhead to the next trip or deadhead first and charge before starting the next trip.

Table \ref{tab:notation_cf} summarizes the notation used in the MILP formulation of CSP. Charging opportunities are indexed consecutively for each bus in a set $\OppSet$, and subsets $\SingleOppSet$ and $\DualOppSet$ are created depending on the type of charging opportunity. For instance, in Figure \ref{fig:diagram_cee}, $\OppSet = \{0, 1, 2, 3\}$, $\SingleOppSet = \{2, 3\}$, and $\DualOppSet=\{1\}$. Opportunities in $\SingleOppSet$ allow charging at only one trip end, whereas ones in $\DualOppSet$ allow charging at both ends. For opportunities in $\DualOppSet$ such as $k=1$, we refer to 1 as the \textit{current trip} and 4 as the \textit{next trip}.

\begin{table}[H]
 \small
  \centering
  \caption{Notation used in the CSP formulation}
    \begin{tabular}{p{14mm}p{145mm}}
    \hline
    \multirow{1}{*}{\textbf{Notation Description}} \\
    \hline
    \multirow{1}{*}{\textbf{Decision variables:}} \\
    \hline
    $\electricityVar_{\busIndex,\locationIndex,\timeIndex}$ & Energy (kWh) provided to bus $\busIndex$ at charging location $\locationIndex$ and time-step $\timeIndex$  \\
    $\EndChargeLocVar{\busIndex}{\timeIndex}$ & Binary variable which is $1$ if bus $\busIndex$ can charge at time-step $\timeIndex$ at the end stop location of a trip, when charging is allowed at both the end points \\
    $\StartChargeLocVar{\busIndex}{\timeIndex}$ & Binary variable which is $1$ if bus $\busIndex$ can charge at time-step $\timeIndex$ at the start stop location of a trip, when charging is allowed at both the end points \\
    $\contractedVar_{\locationIndex}$ & Contracted power capacity (kW) at charging location $\locationIndex$ \\
    \hline
    \multirow{1}{*}{\textbf{Auxiliary variables:}} \\
    \hline
    $\levelVar_{\busIndex,\chargingIndex}$ & Charging level (kWh) of bus $\busIndex$ at the end of its $\chargingIndex^\text{th}$ charging opportunity\\
    \hline
    \textbf{Sets:} \\
    \hline
    $\BusSet$ & Set of electric buses\\
    $\LocationSet$ & Set of charging locations \\
    $\OppSet$ & Set of sequentially enumerated charging opportunities, i.e., trip changes/layovers by $\busIndex$, each with at least a charging station at trip ends. We include 0 to keep track of the charge levels from the depot. \\
    $\SingleOppSet$ & Indices of charging opportunities that have a charging station at only one of the trip ends \\
    $\DualOppSet$ & Indices of charging opportunities that have charging stations at both trip ends \\
    $\timeSet_{\busIndex,\chargingIndex}$ & Set of time-steps of bus $\busIndex$ at charging opportunity $\chargingIndex \in \SingleOppSet$ \\
    $\EndStopTimeSet_{\busIndex,\chargingIndex}$ & Set of time-steps when bus $\busIndex$ can charge at the ending stop of a trip at charging opportunity $\chargingIndex \in \DualOppSet$ \\
    $\StartStopTimeSet_{\busIndex,\chargingIndex}$ & Set of time-steps of bus $\busIndex$ at the starting stop of a trip at charging opportunity $\chargingIndex \in \DualOppSet$ \\
    $\timeperiodSet$ & Set of time periods with different energy prices. Each period is a collection of successive time-steps \\
    $\timeSet_{\periodIndex}$ & Set of time-steps in time period $\periodIndex \in \timeperiodSet$ \\
    $\timeSet$ & Set of all time-steps that cover the period of operations\\
    \hline
    \textbf{Data/Parameters:} \\
    \hline
    $\consumption_{\busIndex,\chargingIndex, \chargingIndex + 1}$ & Energy consumed (kWh) by bus $\busIndex$ to cover the distance from the end of charging opportunity $\chargingIndex \in \OppSet$ to the start of its $(\chargingIndex + 1)^\text{th}$ charging opportunity \\
    $\deadheadconsumption_{\busIndex,\chargingIndex}$ & Energy consumed (kWh) by bus $\busIndex$ for deadheading during charging opportunity $\chargingIndex \in \DualOppSet$  \\
    $\locationIndex_{\busIndex,\chargingIndex}$ & Charging location of bus $\busIndex$ at charging opportunity $\chargingIndex \in \SingleOppSet$\\
    $\EndLoc_{\busIndex,\chargingIndex}$ & End stop of the current trip on bus $\busIndex$ at charging opportunity $\chargingIndex \in \DualOppSet$  \\
    $\startLoc_{\busIndex,\chargingIndex}$ & Start stop of next trip on bus $\busIndex$ at charging opportunity $\chargingIndex \in \DualOppSet$ \\
    $\tau_{b,k}$ & Deadheading time for bus $b$ during charging opportunity $k \in \OppSet$ \\
    $\maxTransfer$ & Maximum amount of charge that can be provided to a bus $\busIndex$ in a single time-step (kWh/min) \\
    $\electricityPrice_{\periodIndex}$ & Electricity price (\$/kWh) at time period $\periodIndex$ \\
    $\capacityPrice$ & Unit price of contracted capacities (\$/kW) of charging locations \\
    $\BatteryCap$ & Battery capacity (kWh) of an electric bus \\
    $\MinEnergy$ & Lower limit of energy level (kWh) of buses \\
    \hline
\end{tabular}%
  \label{tab:notation_cf}%
\end{table}%

\begin{figure}[h]
    \centering
    \includegraphics[scale=0.7]{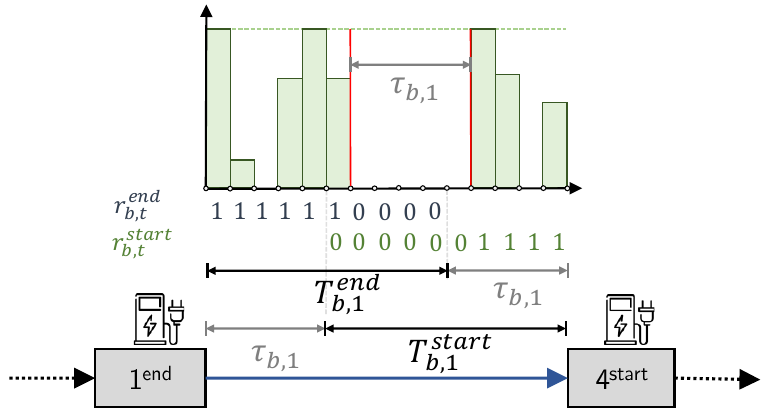}
    \caption{Decision variables at charging opportunities with stations at both trip ends}
    \label{fig:layover}
\end{figure}

We divide the time horizon of operations into one-minute steps and label time intervals where energy prices vary (e.g., 0900--1400, 1400--1600, etc.) as periods. We denote the sets of time steps where a bus can potentially charge at a station as  $\timeSet_{\busIndex,\chargingIndex}$, $\EndStopTimeSet_{\busIndex,\chargingIndex}$, and $\StartStopTimeSet_{\busIndex,\chargingIndex}$. We populate these time steps by excluding the deadheading time associated with a charging opportunity, $\tau_{b,k}$, from the layover time. Figure \ref{fig:layover} shows an example of such sets for a charging opportunity in $\DualOppSet$. 

The objective \eqref{eqn:lp_obj_split_cagac} of the CSP consists of four terms. The first term considers the electricity consumption costs for charging at opportunities in $\SingleOppSet$. The second and third terms represent the electricity costs for charging at the ending stop of the current trip and the starting stop of the next trip, respectively, for charging opportunities in $\DualOppSet$. The final term of the objective function accounts for the contracted power capacity costs. This represents the fixed cost associated with establishing charging infrastructure at different stations. The $\contractedVar_{\locationIndex}$ variable represents the maximum power required to meet the charging needs of buses recharging at station $\locationIndex$. 
\begin{equation}
\min \sum\limits_{\periodIndex \in \timeperiodSet} \electricityPrice_{\periodIndex} \sum\limits_{\busIndex \in \BusSet} \Biggl( \sum\limits_{\substack{\locationIndex \in \LocationSet, \timeIndex \in \timeSet_{\periodIndex} \cap \timeSet_{\busIndex,\chargingIndex}: \\ \chargingIndex \in \SingleOppSet, \locationIndex_{\busIndex,\chargingIndex} = \locationIndex}}  \electricityVar_{\busIndex,\locationIndex,\timeIndex} + 
\sum\limits_{\substack{\locationIndex \in \LocationSet, \timeIndex \in \timeSet_{\periodIndex} \cap \EndStopTimeSet_{\busIndex,\chargingIndex}: \\ \chargingIndex \in \DualOppSet, \EndLoc_{\busIndex,\chargingIndex} = \locationIndex}} \electricityVar_{\busIndex,\locationIndex,\timeIndex} +  \sum\limits_{\substack{\locationIndex \in \LocationSet, \timeIndex \in \timeSet_{\periodIndex} \cap \StartStopTimeSet_{\busIndex,\chargingIndex}: \\ \chargingIndex \in \DualOppSet, \startLoc_{\busIndex,\chargingIndex} = \locationIndex}} \electricityVar_{\busIndex,\locationIndex,\timeIndex} \Biggl) +  \capacityPrice \sum\limits_{\locationIndex \in \LocationSet} \contractedVar_{\locationIndex} \label{eqn:lp_obj_split_cagac} 
\end{equation}
\begin{flalign}
\text{s.t.} \, \, &\levelVar_{\busIndex, \chargingIndex+1} = \levelVar_{\busIndex, \chargingIndex} - \consumption_{\busIndex, \chargingIndex, \chargingIndex + 1} + \sum\limits_{\timeIndex \in \timeSet_{\busIndex, \chargingIndex + 1}} \electricityVar_{\busIndex, \locationIndex,\timeIndex} && \forall \, \busIndex \in \BusSet, \, (\chargingIndex + 1) \in \SingleOppSet, \, \locationIndex = \locationIndex_{\busIndex,\chargingIndex + 1} \label{eqn:cons_split_milp_1}\\
&\levelVar_{\busIndex, \chargingIndex+1} = \levelVar_{\busIndex, \chargingIndex} - \consumption_{\busIndex, \chargingIndex, \chargingIndex + 1} + \sum\limits_{\timeIndex \in \EndStopTimeSet_{\busIndex, \chargingIndex + 1}} \electricityVar_{\busIndex,\locationIndex,\timeIndex} - \deadheadconsumption_{\busIndex,\chargingIndex + 1} + \sum\limits_{\timeIndex \in \StartStopTimeSet_{\busIndex, \chargingIndex + 1}} \electricityVar_{\busIndex,\locationAnotherIndex,\timeIndex} && {\forall \, \busIndex \in \BusSet, \, (\chargingIndex + 1) \in \DualOppSet, \, \locationIndex = \EndLoc_{\busIndex,\chargingIndex + 1}, \, \locationAnotherIndex = \startLoc_{\busIndex,\chargingIndex + 1}} \label{eqn:cons_split_milp_2} \\
\, \, &{\levelVar_{\busIndex,\chargingIndex} - \consumption_{\busIndex, \chargingIndex, \chargingIndex + 1} \geq \MinEnergy} && {\forall \, \busIndex \in \BusSet, \, \chargingIndex \in \OppSet} \label{eqn:cons_split_milp_3a}\\
\, \, &{\levelVar_{\busIndex,\chargingIndex} - \consumption_{\busIndex, \chargingIndex, \chargingIndex + 1} + \sum\limits_{\timeIndex \in \EndStopTimeSet_{\busIndex, \chargingIndex + 1}} \electricityVar_{\busIndex,\locationIndex,\timeIndex} - \deadheadconsumption_{\busIndex,\chargingIndex + 1} \geq \MinEnergy} && {\forall \, \busIndex \in \BusSet, \, (\chargingIndex + 1) \in \DualOppSet, \locationIndex = \EndLoc_{\busIndex,\chargingIndex + 1}} \label{eqn:cons_split_milp_3}\\
\, \, &{\levelVar_{\busIndex,\chargingIndex} \leq \BatteryCap} && {\forall \, \busIndex \in \BusSet, \, \chargingIndex \in \OppSet} \label{eqn:cons_split_milp_4}\\
\, \, &{\levelVar_{\busIndex,\chargingIndex} - \consumption_{\busIndex, \chargingIndex, \chargingIndex + 1} + \sum\limits_{\timeIndex \in \EndStopTimeSet_{\busIndex, \chargingIndex + 1}} \electricityVar_{\busIndex,\locationIndex,\timeIndex} \leq \BatteryCap} && {\forall \, \busIndex \in \BusSet, \, (\chargingIndex + 1) \in \DualOppSet, \locationIndex = \EndLoc_{\busIndex,\chargingIndex + 1}} \label{eqn:cons_split_milp_5}\\
&\levelVar_{\busIndex,0} = \BatteryCap && \forall \, \busIndex \in \BusSet \label{eqn:cons_split_milp_6}\\
& 0 \leq \electricityVar_{\busIndex,\locationIndex,\timeIndex} \leq \maxTransfer && \forall \, \busIndex \in \BusSet, \, \chargingIndex \in \SingleOppSet, \, \locationIndex = \locationIndex_{\busIndex,\chargingIndex}, \, \timeIndex \in \timeSet_{\busIndex,\chargingIndex} \label{eqn:cons_split_milp_7} \\
& 0 \leq \electricityVar_{\busIndex,\locationIndex,\timeIndex} \leq \maxTransfer \, \EndChargeLocVar{\busIndex}{\timeIndex} \label{eqn:cons_split_milp_8} && {\forall \, \busIndex \in \BusSet, \, \chargingIndex \in \DualOppSet, \, \locationIndex = \EndLoc_{\busIndex,\chargingIndex}, \, \timeIndex \in \EndStopTimeSet_{\busIndex,\chargingIndex}} \\
& 0 \leq \electricityVar_{\busIndex,\locationIndex,\timeIndex} \leq \maxTransfer \, \StartChargeLocVar{\busIndex}{\timeIndex} \label{eqn:cons_split_milp_9} && {\forall \, \busIndex \in \BusSet, \, \chargingIndex \in \DualOppSet, \, \locationIndex = \startLoc_{\busIndex,\chargingIndex}}, \, \timeIndex \in \StartStopTimeSet_{\busIndex,\chargingIndex} \\
&\EndChargeLocVar{\busIndex}{\timeIndex+1} \leq \EndChargeLocVar{\busIndex}{\timeIndex} \label{eqn:cons_split_milp_10} && {\forall \, \busIndex \in \BusSet, \, \chargingIndex \in \DualOppSet, \, \timeIndex, \, (\timeIndex + 1) \in \EndStopTimeSet_{\busIndex,\chargingIndex}} \\
&\StartChargeLocVar{\busIndex}{\timeIndex+1} \geq \StartChargeLocVar{\busIndex}{\timeIndex} \label{eqn:cons_split_milp_11} && {\forall \, \busIndex \in \BusSet, \, \chargingIndex \in \DualOppSet}, \, \timeIndex, \, (\timeIndex + 1) \in \StartStopTimeSet_{\busIndex,\chargingIndex}\\
&\StartChargeLocVar{\busIndex}{\timeIndex +\tau_{b, k}} \leq 1 - \EndChargeLocVar{\busIndex}{\timeIndex} \label{eqn:cons_split_milp_13} && \forall \, \busIndex \in \BusSet, \, \chargingIndex \in \DualOppSet, \,  \timeIndex \in \EndStopTimeSet_{\busIndex,\chargingIndex} \\
&\sum\limits_{\substack{\busIndex \in \BusSet, \chargingIndex \in \SingleOppSet: \\ \locationIndex_{\busIndex,\chargingIndex} = \locationIndex, \timeIndex \in \timeSet_{\busIndex,\chargingIndex}}} \electricityVar_{\busIndex,\locationIndex,\timeIndex} +  \sum\limits_{\substack{\busIndex \in \BusSet, \chargingIndex \in \DualOppSet: \\ \EndLoc_{\busIndex,\chargingIndex} = \locationIndex, \timeIndex \in \EndStopTimeSet_{\busIndex,\chargingIndex}}} \electricityVar_{\busIndex,\locationIndex,\timeIndex} + \sum\limits_{\substack{\busIndex \in \BusSet, \chargingIndex \in \DualOppSet: \\ \startLoc_{\busIndex,\chargingIndex} = \locationIndex, \timeIndex \in \StartStopTimeSet_{\busIndex,\chargingIndex}}} \electricityVar_{\busIndex,\locationIndex,\timeIndex} \leq  \contractedVar_{\locationIndex} / 60 && \forall \, \locationIndex \in \LocationSet, \, \timeIndex \in \timeSet \label{eqn:cons_split_milp_14} \\
&\EndChargeLocVar{\busIndex}{\timeIndex} \in \{0,1\}, \, \StartChargeLocVar{\busIndex}{\timeAnotherIndex} \in \{0,1\} \label{eqn:cons_split_milp_15} && {\forall \, \busIndex \in \BusSet, \,  \chargingIndex \in \DualOppSet, \, \timeIndex \in \EndStopTimeSet_{\busIndex,\chargingIndex}}, \timeAnotherIndex \in \StartStopTimeSet_{\busIndex,\chargingIndex}
\end{flalign}
Constraints \eqref{eqn:cons_split_milp_1} and \eqref{eqn:cons_split_milp_2} update the charging levels between two consecutive charging opportunities for each bus. Specifically, \eqref{eqn:cons_split_milp_1} applies when charging at the ($\chargingIndex+1)^\text{th}$ opportunity is permitted only at one of the end stops, while \eqref{eqn:cons_split_milp_2} applies when charging is allowed at both end stops. Constraints \eqref{eqn:cons_split_milp_3a} and \eqref{eqn:cons_split_milp_3} ensure that the buses' charging levels do not drop below a minimum threshold at any time. Constraints \eqref{eqn:cons_split_milp_4} and \eqref{eqn:cons_split_milp_5} prevent buses from charging beyond their battery capacity. Buses start their daily operations with a fully charged battery due to \eqref{eqn:cons_split_milp_6}. The maximum energy transferred to a bus within a single time-step is governed by the availability of the bus at charging locations, which is captured in constraints \eqref{eqn:cons_split_milp_7}--\eqref{eqn:cons_split_milp_9}. When a charging opportunity $k$ is in $\SingleOppSet$, we simply set the upper bound based on the maximum rate $\maxTransfer$ at which the batteries on the bus can be charged. For a charging opportunity $k \in \DualOppSet$, we use binary variables $\EndChargeLocVar{\busIndex}{\timeIndex}$ and $\StartChargeLocVar{\busIndex}{\timeIndex}$ to decide on the time-steps where charging is permitted at the trip ends (see Figure \ref{fig:layover}). Constraint \eqref{eqn:cons_split_milp_10} ensures that the $\EndChargeLocVar{\busIndex}{\timeIndex}$ variables exhibit a sequence of consecutive ones followed by zeros. Similarly, \eqref{eqn:cons_split_milp_11} forces $\StartChargeLocVar{\busIndex}{\timeIndex}$ variables to assume a series of zeros followed by ones. 
Furthermore, charging at the start stop of the next trip can begin only after finishing charging at the current trip's end stop and the deadheading operation. Hence, the $\StartChargeLocVar{\busIndex}{\timeIndex}$ variables in $\StartStopTimeSet_{\busIndex,\chargingIndex}$ can be active after considering the deadheading time from the last time-step where $\EndChargeLocVar{\busIndex}{\timeIndex}$ is one, as specified in \eqref{eqn:cons_split_milp_13}. Finally, \eqref{eqn:cons_split_milp_14} and the minimization objective help set $\contractedVar_{\locationIndex}$ to the maximum amount of energy used by the charging station at $\locationIndex$ in a single time step. Binary restrictions on the deadheading decision variables are imposed using \eqref{eqn:cons_split_milp_15}.

\subsection{CLP-EVSP-CSP Model}
\label{sec:cl_evsp_csp}
We enhance CLP-EVSP formulations found in the literature by incorporating CSP-specific variables. Instead of charging opportunities, we utilize trip indices since the locations are not predetermined and are determined jointly by the model. Table \ref{tab:joint_mip} introduces additional notation to describe the MILP. The notation from the CSP formulation remains applicable.

\begin{table}[h]
 \small
  \centering
  \caption{Additional notation used in the joint optimization formulation of CLP-EVSP-CSP}
    \begin{tabular}{p{14mm}p{153mm}}
    \hline
   \multirow{1}{*}{\textbf{Notation Description}} \\
    \hline
    \multirow{1}{*}{\textbf{Decision variables:}} \\
    \hline
    $\arcVar{\tripIndexI}{\tripIndexJ}{\busIndex}$ & Binary variable which is $1$ if node $\tripIndexJ$ is visited after node $\tripIndexI$ by bus $\busIndex$, and $0$ otherwise \\
    $\LocVar_{\locationIndex}$ & Binary variable which is $1$ if a charging facility is set up at stop $\locationIndex$, and $0$ otherwise \\
    $\electricityStartVar{\busIndex}{\tripIndexI}{\timeIndex}$ & Energy (kWh) provided to bus $\busIndex$ at the starting stop of trip $\tripIndexI$ at time-step $\timeIndex$\\
    $\electricityEndVar{\busIndex}{\tripIndexI}{\timeIndex}$ & Energy (kWh) provided to bus $\busIndex$ at the ending stop of trip $\tripIndexI$ at time-step $\timeIndex$\\
     \hline
     \multirow{1}{*}{\textbf{Auxiliary variables:}} \\
    \hline
    $\StartChargeTimeVar{\tripIndexI}{\timeIndex}$ & Binary variable which is $1$ if charging can be done at the starting stop of trip $\tripIndexI$ at time-step $\timeIndex$, and $0$ otherwise \\
    $\EndChargeTimeVar{\tripIndexI}{\timeIndex}$ & Binary variable which is $1$ if charging can be done at the ending stop of trip $\tripIndexI$ at time-step $\timeIndex$, and $0$ otherwise\\
    $\StartLevelVar{\tripIndexI}{\busIndex}$ & Energy level (kWh) of bus $\busIndex$ at the starting stop of trip $\tripIndexI$ or at the end of depot $\tripIndexI$\\
    
    $\EndLevelVar{\tripIndexI}{\busIndex}$ & Energy level (kWh) of bus $\busIndex$ at the ending stop of trip $\tripIndexI$ or at the start of depot $\tripIndexI$ \\
    \hline
    \textbf{Sets:} \\
    \hline
    $\tripSet$ & Set of trips \\
    $\depotSet$ & Set of depots \\
    $\nodeSet$ & Set of nodes for the EVSP \\
    $\arcSet$ & Set of arcs for the EVSP \\
    $\compatibleSet$ & Set of compatible trip pairs \\
    $\candidateLocationSet$ & Set of candidate charging locations\\
    $\OriginStopsSet_{\locationIndex}$ & Set of trips with stop $\locationIndex$ as starting/origin stop \\
    $\DestinationStopsSet_{\locationIndex}$ & Set of trips with stop $\locationIndex$ as ending/destination stop \\
    $\StartStopTimeSet_{\tripIndexI,\tripIndexJ}$ & Set of time-steps during which a bus can charge at the starting stop of node $\tripIndexJ$ if the preceding node was $\tripIndexI$\\
    $\EndStopTimeSet_{\tripIndexI, \tripIndexJ}$ & Set of time-steps during which a bus can charge at the ending stop of node $\tripIndexI$ before moving to node $\tripIndexJ$\\
    $\PreviousTimeSet{\tripIndexI}$ & Set of time-steps that precedes the starting time of trip $\tripIndexI$ \\
    $\NextTimeSet{\tripIndexI}$ & Set of time-steps that succeeds the ending time of trip $\tripIndexI$ \\
    \hline
    \textbf{Data/Parameters:} \\
    \hline
    $\energyTrip{\tripIndexI}$ & Energy consumed (kWh) to travel the distance of trip $\tripIndexI$\\
    $\energyDeadhead{\tripIndexI}{\tripIndexJ}$ & Energy consumed (kWh) in deadheading from the end of node $\tripIndexI$ to the start of node $\tripIndexJ$ \\
    $\BigMConst$ & A large positive number \\
    $\TripStartStop$ & Starting stop of trip $\tripIndexI$\\
    $\TripEndStop$ & Ending stop of trip $\tripIndexI$\\
    $\deadheadDur_{\tripIndexI,\tripIndexJ}$ & Deadhead trip duration from the end stop of $\tripIndexI$ to the start stop of $\tripIndexJ$ \\
    $\BusCost$ & Acquisition cost (\$) of an electric bus  \\
    $\LocationCost$ & Fixed investment cost (\$) per charging location \\
    $\DeadheadCost{\tripIndexI}{\tripIndexJ}$ & Total energy consumption costs (\$) from the end location of node $\tripIndexI$ to the start location of node $\tripIndexJ$ (which equals the per km cost multiplied by the corresponding deadhead distance traveled) \\
    \hline 
    \end{tabular}%
  \label{tab:joint_mip}%
\end{table}%

The objective \eqref{eqn:obj_joint} comprises five terms: the acquisition cost of electric buses required to serve all the trips, the fixed cost of a charging station, the deadheading costs, the total electricity consumption cost, and the total cost for installing sufficient power capacity at all charging stations. The deadheading and electricity consumption costs are scaled appropriately to account for the life cycle of electric buses. 

\begin{equation}
\label{eqn:obj_joint}
\min \sum\limits_{\busIndex \in \BusSet} \sum_{\substack{(\tripIndexI,\tripIndexJ) \in \arcSet: \\ \tripIndexI \in \depotSet}}   \BusCost \, \arcVar{\tripIndexI}{\tripIndexJ}{\busIndex} + \sum\limits_{\locationIndex \in \candidateLocationSet} \LocationCost \LocVar_{\locationIndex} + \sum\limits_{\busIndex \in \BusSet} \sum\limits_{(\tripIndexI,\tripIndexJ) \in \arcSet} \DeadheadCost{\tripIndexI}{\tripIndexJ} \arcVar{\tripIndexI}{\tripIndexJ}{\busIndex} + \sum\limits_{\periodIndex \in \timeperiodSet} \sum\limits_{\timeIndex \in \timeSet_{\periodIndex}} \sum\limits_{\tripIndexI \in \tripSet} \sum\limits_{\busIndex \in \BusSet} \electricityPrice_{\periodIndex}  (\electricityStartVar{\busIndex}{\tripIndexI}{\timeIndex} + \electricityEndVar{\busIndex}{\tripIndexI}{\timeIndex}) + \capacityPrice \sum\limits_{\locationIndex \in \candidateLocationSet} \contractedVar_{\locationIndex}
\end{equation}
\begin{flalign}
\label{eqn:cons1}
\text{s.t.} & \sum\limits_{\tripIndexJ: (\tripIndexI,\tripIndexJ) \in \arcSet} \arcVar{\tripIndexI}{\tripIndexJ}{\busIndex} - \sum\limits_{\tripIndexK: (\tripIndexK,\tripIndexI) \in \arcSet} \arcVar{\tripIndexK}{\tripIndexI}{\busIndex} = 0 && \, \forall \, \tripIndexI \in \tripSet, \, \busIndex \in \BusSet\\
\label{eqn:cons2a}
& \sum\limits_{\busIndex \in \BusSet} \sum\limits_{\tripIndexJ: (\tripIndexI,\tripIndexJ) \in \arcSet} \arcVar{\tripIndexI}{\tripIndexJ}{\busIndex} - \sum\limits_{\busIndex \in \BusSet} \sum\limits_{\tripIndexK: (\tripIndexK,\tripIndexI) \in \arcSet} \arcVar{\tripIndexK}{\tripIndexI}{\busIndex} = 0 && \, \forall \, \tripIndexI \in \depotSet\\
\label{eqn:cons2}
& \sum\limits_{\busIndex \in \BusSet} \sum\limits_{\tripIndexJ: (\tripIndexI,\tripIndexJ) \in \arcSet} \arcVar{\tripIndexI}{\tripIndexJ}{\busIndex} = 1 && \, \forall \, \tripIndexI \in \tripSet\\
\label{eqn:cons3}
& \sum\limits_{\tripIndexI \in \depotSet} \sum\limits_{\tripIndexJ: (\tripIndexI,\tripIndexJ) \in \arcSet} \arcVar{\tripIndexI}{\tripIndexJ}{\busIndex} \leq 1 && \, \forall \, \busIndex \in \BusSet\\
\label{eqn:cons7}
& \StartLevelVar{\tripIndexI}{\busIndex} \leq \BatteryCap \sum\limits_{\tripIndexK: (\tripIndexK,\tripIndexI) \in \arcSet} \arcVar{\tripIndexK}{\tripIndexI}{\busIndex} && \, \forall \, \tripIndexI \in \tripSet,  \, \busIndex \in \BusSet \\
\label{eqn:cons8}
& \EndLevelVar{\tripIndexI}{\busIndex} \leq \BatteryCap \sum\limits_{\tripIndexJ: (\tripIndexI,\tripIndexJ) \in \arcSet} \arcVar{\tripIndexI}{\tripIndexJ}{\busIndex} && \, \forall \, \tripIndexI \in \tripSet, \, \busIndex \in \BusSet \\
\label{eqn:cons9}
& \EndLevelVar{\tripIndexI}{\busIndex} = \BatteryCap \sum\limits_{\tripIndexJ: (\tripIndexI,\tripIndexJ) \in \arcSet} \arcVar{\tripIndexI}{\tripIndexJ}{\busIndex} && \, \forall \, \busIndex \in \BusSet, \, \tripIndexI \in \depotSet\\
\label{eqn:cons10}
& \StartLevelVar{\tripIndexI}{\busIndex} - \sum\limits_{\tripIndexK: (\tripIndexK,\tripIndexI) \in \arcSet} \energyTrip{\tripIndexI} \, \arcVar{\tripIndexK}{\tripIndexI}{\busIndex} \geq \MinEnergy \sum\limits_{\tripIndexK:(\tripIndexK,\tripIndexI) \in \arcSet} \arcVar{\tripIndexK}{\tripIndexI}{\busIndex} && \, \forall \, \tripIndexI \in \tripSet,  \, \busIndex \in \BusSet \\
\label{eqn:cons11}
& \EndLevelVar{\tripIndexI}{\busIndex} - \sum\limits_{\tripIndexJ: (\tripIndexI,\tripIndexJ) \in \arcSet} \energyDeadhead{\tripIndexI}{\tripIndexJ} \, \arcVar{\tripIndexI}{\tripIndexJ}{\busIndex} \geq \MinEnergy \sum\limits_{\tripIndexJ:(\tripIndexI,\tripIndexJ) \in \arcSet} \arcVar{\tripIndexI}{\tripIndexJ}{\busIndex} && \, \forall  \, \tripIndexI \in \tripSet,  \,  \busIndex \in \BusSet \\
\label{eqn:cons12}
& \StartLevelVar{\tripIndexJ}{\busIndex} \geq \MinEnergy \, \sum\limits_{\tripIndexI \in \tripSet} \arcVar{\tripIndexI}{\tripIndexJ}{\busIndex} && \,  \forall  \, \tripIndexJ \in \depotSet,  \, \busIndex \in \BusSet \\
\label{eqn:cons_joint_4}
& \EndLevelVar{\tripIndexI}{\busIndex} \leq \StartLevelVar{\tripIndexI}{\busIndex} - \energyTrip{\tripIndexI} + \sum\limits_{\timeIndex \in \EndStopTimeSet_{\tripIndexI, \tripIndexJ}} \electricityEndVar{\busIndex}{\tripIndexI}{\timeIndex} + \BigMConst (1- \arcVar{\tripIndexI}{\tripIndexJ}{\busIndex}) && \, \forall \, \busIndex \in \BusSet, \, \, \tripIndexI \in \tripSet, \, (\tripIndexI, \tripIndexJ) \in \arcSet \\
\label{eqn:cons_joint_5}
& \EndLevelVar{\tripIndexI}{\busIndex} \geq \StartLevelVar{\tripIndexI}{\busIndex} - \energyTrip{\tripIndexI} +\sum\limits_{\timeIndex \in \EndStopTimeSet_{\tripIndexI, \tripIndexJ}} \electricityEndVar{\busIndex}{\tripIndexI}{\timeIndex} - \BigMConst (1- \arcVar{\tripIndexI}{\tripIndexJ}{\busIndex}) && \, \forall \, \busIndex \in \BusSet, \tripIndexI \in \tripSet, \, (\tripIndexI, \tripIndexJ) \in \arcSet \\
\label{eqn:cons_joint_6}
& \StartLevelVar{\tripIndexJ}{\busIndex} \geq \EndLevelVar{\tripIndexI}{\busIndex} - \energyDeadhead{\tripIndexI}{\tripIndexJ} \, \arcVar{\tripIndexI}{\tripIndexJ}{\busIndex} + \sum\limits_{\timeIndex \in \StartStopTimeSet_{\tripIndexI, \tripIndexJ}} \electricityStartVar{\busIndex}{\tripIndexJ}{\timeIndex} - \BigMConst (1 - \arcVar{\tripIndexI}{\tripIndexJ}{\busIndex}) && \, \forall \, \busIndex \in \BusSet, \, (\tripIndexI,\tripIndexJ) \in \arcSet\\
\label{eqn:cons_joint_7}
& \StartLevelVar{\tripIndexJ}{\busIndex} \leq \EndLevelVar{\tripIndexI}{\busIndex} - \energyDeadhead{\tripIndexI}{\tripIndexJ} \,\arcVar{\tripIndexI}{\tripIndexJ}{\busIndex} + \sum\limits_{\timeIndex \in \StartStopTimeSet_{\tripIndexI, \tripIndexJ}} \electricityStartVar{\busIndex}{\tripIndexJ}{\timeIndex} + \BigMConst (1 - \arcVar{\tripIndexI}{\tripIndexJ}{\busIndex}) && \, \forall \, \busIndex \in \BusSet, \, (\tripIndexI,\tripIndexJ) \in \arcSet \\
\label{eqn:cons_joint_15}
& \electricityStartVar{\busIndex}{\tripIndexI}{\timeIndex} \leq \maxTransfer \, \StartChargeTimeVar{\tripIndexI}{\timeIndex} &&  \, \forall \, \busIndex \in \BusSet, \, \tripIndexI \in \tripSet, \,  \timeIndex \in \PreviousTimeSet{\tripIndexI} \\
\label{eqn:cons_joint_16}
& \electricityEndVar{\busIndex}{\tripIndexI}{\timeIndex} \leq \maxTransfer \, \EndChargeTimeVar{\tripIndexI}{\timeIndex} && \, \forall \, \busIndex \in \BusSet, \, \tripIndexI \in \tripSet, \, \timeIndex \in \NextTimeSet{\tripIndexI} \\
\label{eqn:cons_joint_16a}
& \StartChargeTimeVar{\tripIndexI}{\timeIndex + 1} \geq \StartChargeTimeVar{\tripIndexI}{\timeIndex} && \, \forall \, \tripIndexI \in \tripSet, \, \timeIndex, \, (\timeIndex + 1) \in \PreviousTimeSet{\tripIndexI} \\
\label{eqn:cons_joint_16b}
& \EndChargeTimeVar{\tripIndexI}{\timeIndex + 1} \leq \EndChargeTimeVar{\tripIndexI}{\timeIndex} && \, \forall \, \tripIndexI \in \tripSet, \, \timeIndex, \, (\timeIndex + 1) \in \NextTimeSet{\tripIndexI} \\
\label{eqn:cons_joint_16d}
& \StartChargeTimeVar{\tripIndexJ}{\timeIndex+\gamma_{i,j}} \leq 1 - \EndChargeTimeVar{\tripIndexI}{\timeIndex} + \BigMConst (1 - \arcVar{\tripIndexI}{\tripIndexJ}{\busIndex}) && \, \forall \, \busIndex \in \BusSet, \, (\tripIndexI,\tripIndexJ) \in \compatibleSet, \timeIndex \in \EndStopTimeSet_{\tripIndexI,\tripIndexJ} \\
\label{eqn:cons_joint_17}
& \sum\limits_{\tripIndexI \in \OriginStopsSet_{\locationIndex}} \sum\limits_{\timeIndex \in \PreviousTimeSet{\tripIndexI}} \StartChargeTimeVar{\tripIndexI}{\timeIndex} + \sum\limits_{\tripIndexI \in \DestinationStopsSet_{\locationIndex}} \sum\limits_{\timeIndex \in \NextTimeSet{\tripIndexI}} \EndChargeTimeVar{\tripIndexI}{\timeIndex} \leq \BigMConst \LocVar_{\locationIndex} && \, \forall \, \locationIndex \in \candidateLocationSet\\
\label{eqn:cons_joint_18}
& \sum\limits_{\busIndex \in \BusSet} \sum\limits_{\substack{\tripIndexI \in \tripSet: \\ \TripEndStop = \locationIndex}} \electricityEndVar{\busIndex}{\tripIndexI}{\timeIndex} + \sum\limits_{\busIndex \in \BusSet} \sum\limits_{\substack{\tripIndexI \in \tripSet: \\ \TripStartStop = \locationIndex}} \electricityStartVar{\busIndex}{\tripIndexI}{\timeIndex} \leq \contractedVar_{\locationIndex} / 60 && \forall \,\timeIndex \in \timeSet, \, \locationIndex \in \candidateLocationSet \\
\label{eqn:cons_joint_19}
&\contractedVar_{\locationIndex} \leq \BigMConst \LocVar_{\locationIndex} && \, \forall \, \locationIndex \in \candidateLocationSet \\
\label{eqn:cons_joint_20}
&\arcVar{\tripIndexI}{\tripIndexJ}{\busIndex}, \, \LocVar_{\locationIndex}, \, \StartChargeTimeVar{\tripIndexI}{\timeIndex}, \EndChargeTimeVar{\tripIndexI}{\timeAnotherIndex} \in \{0,1\} && \, \forall \, \tripIndexI, \tripIndexJ \in \nodeSet, \busIndex \in \BusSet, \timeIndex \in \PreviousTimeSet{\tripIndexI}, \timeAnotherIndex \in \NextTimeSet{\tripIndexI}, \locationIndex \in \candidateLocationSet \\
\label{eqn:cons_joint_21}
& \electricityStartVar{\busIndex}{\tripIndexI}{\timeIndex}, \electricityEndVar{\busIndex}{\tripIndexI}{\timeAnotherIndex}, 
\contractedVar_{\locationIndex}, \StartLevelVar{\tripIndexI}{\busIndex}, \EndLevelVar{\tripIndexI}{\busIndex} \in \mathbb{R^{+}} && \, \forall \, \tripIndexI \in \nodeSet, \busIndex \in \BusSet, \timeIndex \in \PreviousTimeSet{\tripIndexI}, \timeAnotherIndex \in \NextTimeSet{\tripIndexI}, \locationIndex \in \candidateLocationSet 
\end{flalign}
We ensure flow conservation at every trip node and depot node using \eqref{eqn:cons1} and \eqref{eqn:cons2a}, respectively. Constraint \eqref{eqn:cons2} indicates that each trip is served by exactly one bus. Constraint \eqref{eqn:cons3} models the scenario where a bus may remain unused and determines the buses starting at different depots. Inequalities \eqref{eqn:cons7} and \eqref{eqn:cons8} ensure that the battery can be charged up to the maximum capacity, avoiding overcharging. Used buses leave the depot fully charged according to \eqref{eqn:cons9}. Constraints \eqref{eqn:cons10}, \eqref{eqn:cons11}, and \eqref{eqn:cons12} set minimum charge-level requirements so that buses have sufficient energy to complete their assigned trips. Specifically, \eqref{eqn:cons10} handles service trips, \eqref{eqn:cons11} is applied to deadhead trips, and \eqref{eqn:cons12} deals with depot trips. Constraints \eqref{eqn:cons_joint_4} and \eqref{eqn:cons_joint_5} track the energy levels of each bus throughout its assigned trips, using the appropriate energy consumption and battery replenishment variables. Similar inequalities are defined for the deadhead trips (including the pull-in and pull-out depot trips and the compatible service trips) using \eqref{eqn:cons_joint_6} and \eqref{eqn:cons_joint_7}. Constraints  \eqref{eqn:cons_joint_15} and \eqref{eqn:cons_joint_16} restrict the amount of energy that can be transferred to a bus at every time-step. Constraints \eqref{eqn:cons_joint_16a}--\eqref{eqn:cons_joint_16d} model the split charging of buses at both ends of a layover in a manner similar to \eqref{eqn:cons_split_milp_10}--\eqref{eqn:cons_split_milp_13}. Constraint \eqref{eqn:cons_joint_17} allows charging activities to occur if a charging facility is located at a terminal stop. 
Constraint \eqref{eqn:cons_joint_18} sets $\contractedVar_{\locationIndex}$ to the maximum power supply needed for buses charging at $\locationIndex$. It is permitted to be non-negative only if a charging station is opened at $\locationIndex$, as stipulated by \eqref{eqn:cons_joint_19}. Finally, binary restrictions and limits on the decision variables are imposed by \eqref{eqn:cons_joint_20} and \eqref{eqn:cons_joint_21}.

\subsection{Surrogate CSPs}
\label{sec:surrogate_lp}
The CEE MILP presented in Section \ref{sec:csp} can be simplified using surrogate LP models. More importantly, these LPs enable the development of efficient local search heuristics, which will be discussed in Section \ref{sec:localsearch}. In these variants, we restrict buses to charge at only one of the end stops during a layover when charging stations are available at both ends. This approach eliminates the need for integer variables that determine the timing of deadheading. Two types of priority rules can be applied at such layovers: a \textit{charge-and-go} (CAG) strategy, which involves charging at the end stop of the current trip before deadheading, and a \textit{go-and-charge} (GAC) strategy, where deadheading is performed first, followed by charging at the starting stop of the next trip. Both models can be further simplified by assuming that the energy supplied to the buses either changes dynamically over time (\textit{split charging}) or remains constant (\textit{uniform charging}) during each charging opportunity. Our experiments demonstrate that these LPs are highly tractable, and there is an improvement in the overall objective and the CSP costs when they are integrated with other location and scheduling problems.

\subsubsection{Split Charging with Priority}
\label{sec:charging_lp}
\begin{figure}[H]
    \centering
    \begin{subfigure}{\textwidth}
        \centering
        \includegraphics[width=0.95\textwidth]{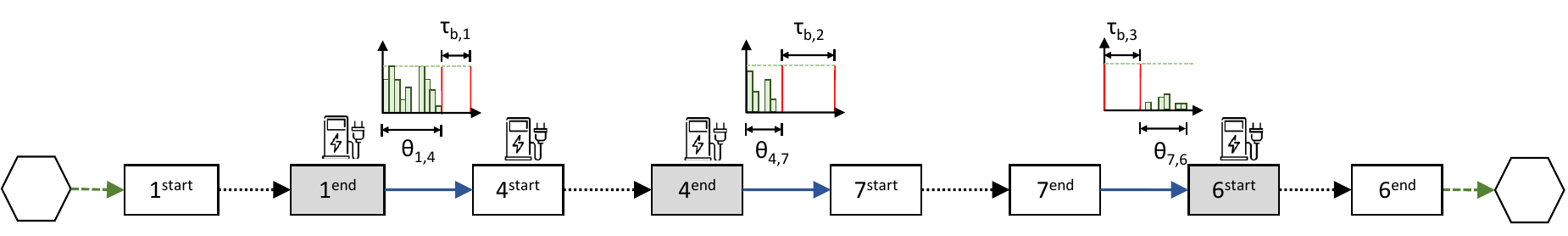}
        \caption{Charge-and-Go (CAG) strategy}
        \label{fig:diagram_cag}
    \end{subfigure}
    \\
    \begin{subfigure}{\textwidth}
        \centering
        \includegraphics[width=0.95\textwidth]{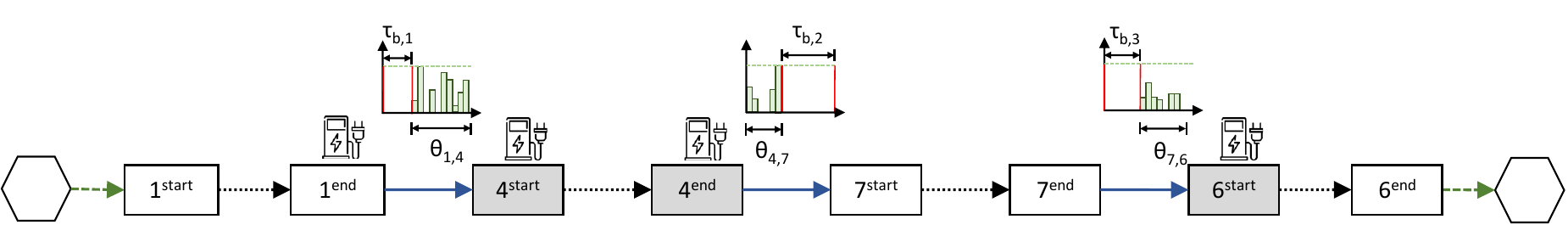}
        \caption{Go-and-Charge (GAC) strategy}
        \label{fig:diagram_gac}
    \end{subfigure}
    \caption{Illustration of split charging strategies}
    \label{fig:split}
\end{figure}

Figure \ref{fig:split} depicts an example of the split charging strategy. For the layover between trips 1 and 4, which has charging stations at both ends, the CAG strategy prioritizes charging at the ending stop of trip 1 (indicated by the gray node) as shown in Figure \ref{fig:diagram_cag}. Similarly, Figure \ref{fig:diagram_gac} demonstrates the GAC strategy, where charging at the starting stop of trip 4 is given priority. The formulations for the CAG and GAC models are similar, differing only in the input data. 
\begin{align}
\min &\sum\limits_{\periodIndex \in \timeperiodSet} \sum\limits_{\busIndex \in \BusSet} \, \sum\limits_{\substack{\locationIndex \in \LocationSet, \timeIndex \in \timeSet_{\periodIndex} \cap \timeSet_{\busIndex,\chargingIndex}: \\ \chargingIndex \in \SingleOppSet, \locationIndex_{\busIndex,\chargingIndex} = \locationIndex}} \, \electricityPrice_{\periodIndex} \, \electricityVar_{\busIndex,\locationIndex,\timeIndex} + \capacityPrice \sum\limits_{\locationIndex \in \LocationSet} \, \, \contractedVar_{\locationIndex} \label{eqn:lp_obj_split} \\
\text{s.t.} \, \, &\eqref{eqn:cons_split_milp_1}, \eqref{eqn:cons_split_milp_3a}, \eqref{eqn:cons_split_milp_4}, \eqref{eqn:cons_split_milp_6}, \eqref{eqn:cons_split_milp_7} \, \text{with} \, \DualOppSet = \emptyset \nonumber \\
&\sum\limits_{\substack{\busIndex \in \BusSet, \chargingIndex \in \SingleOppSet: \\ \locationIndex_{\busIndex,\chargingIndex} = \locationIndex, \timeIndex \in \timeSet_{\busIndex,\chargingIndex}}} \electricityVar_{\busIndex,\locationIndex,\timeIndex} \leq \contractedVar_{\locationIndex}/60 && \forall \, \locationIndex \in \LocationSet, \, \timeIndex \in \timeSet \label{eqn:cons_split_lp_6}
\end{align}
The objective \eqref{eqn:lp_obj_split} consists of two terms: the electricity consumption costs in different periods and the contracted capacity costs across all charging locations. Constraint \eqref{eqn:cons_split_lp_6} models the contracted power capacity requirements of charging stations and is a simplified version of \eqref{eqn:cons_split_milp_14}. We experimented with these strategies using a specific set of location and bus rotations and found that the objective of the CAG strategy is closer to that of the optimal CEE solution for most networks. Therefore, for the remainder of this paper, we will assume that the CSP is solved using the CAG strategy unless stated otherwise.
\subsubsection{Uniform Charging with Priority}
\begin{figure}[H]
    \centering
    \includegraphics[width=0.95\textwidth]{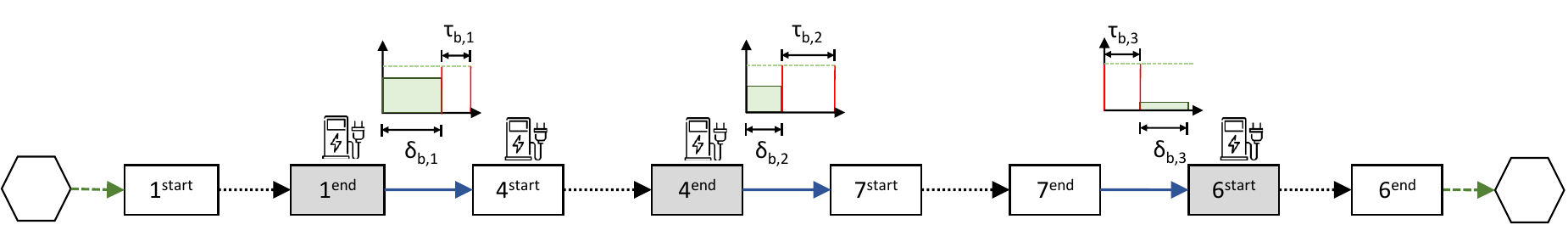}
	\caption{Illustration of the uniform CAG charging strategy}
\label{fig:csp_uniform}
\end{figure}

To further reduce the number of variables and the runtime of the local search operators, one could consider implementing a uniform charging strategy. In this formulation, we assume that buses are charged throughout their entire idle period at a charging station, thereby eliminating the need for time indices as shown in the example in Figure \ref{fig:csp_uniform}. We determine the charging levels based on the energy requirements of the buses and their idle duration at the charging locations. Let $\electricityVar_{\busIndex,\chargingIndex}$ represent the energy provided and $\delta_{\busIndex, \chargingIndex}$ denote the idle time for bus $\busIndex$ at its $\chargingIndex^\text{th}$ charging opportunity. The uniform charging formulation can be written as follows.
\begin{align}
\min &\sum\limits_{\periodIndex \in \timeperiodSet} \sum\limits_{\busIndex \in \BusSet} \sum\limits_{\substack{\chargingIndex \in \SingleOppSet: \\ \timeIndex \in \timeSet_{\periodIndex} \cap \timeSet_{\busIndex,\chargingIndex}}} \electricityPrice_{\periodIndex} \, \frac{\electricityVar_{\busIndex,\chargingIndex}}{\delta_{\busIndex, \chargingIndex}} + \capacityPrice \sum\limits_{\locationIndex \in \LocationSet} \, \,  \contractedVar_{\locationIndex} \label{eqn:lp_cont_obj}\\
\text{s.t.} \, \, &\levelVar_{ \busIndex, \chargingIndex+1} = \levelVar_{\busIndex, \chargingIndex} - \consumption_{\busIndex, \chargingIndex, \chargingIndex + 1} + \electricityVar_{\busIndex ,\chargingIndex + 1} && \forall \, \busIndex \in \BusSet, \, (\chargingIndex + 1) \in \SingleOppSet \label{eqn:cons_cont_lp_1}\\
\, \, &{\levelVar_{\busIndex, \chargingIndex} \leq \BatteryCap} && {\forall \, \busIndex \in \BusSet, \, \chargingIndex \in \SingleOppSet} \label{eqn:cons_cont_lp_2}\\
\, \, &{\levelVar_{\busIndex, \chargingIndex} - \consumption_{\busIndex, \chargingIndex, \chargingIndex + 1} \geq \MinEnergy} && {\forall \, \busIndex \in \BusSet, \, \chargingIndex \in \SingleOppSet} \label{eqn:cons_cont_lp_3}\\
&\levelVar_{\busIndex, 0} = \BatteryCap && \forall \, \busIndex \in \BusSet \label{eqn:cons_cont_lp_4}\\
&0 \leq \electricityVar_{\busIndex, \chargingIndex} \leq \maxTransfer \, \delta_{\busIndex, \chargingIndex} && \forall \, \busIndex \in \BusSet, \, \chargingIndex \in \SingleOppSet\label{eqn:cons_cont_lp_5} \\
&\sum\limits_{\substack{\busIndex \in \BusSet, \chargingIndex \in \SingleOppSet: \\ \locationIndex_{\busIndex, \chargingIndex} = \locationIndex, \timeIndex \in \timeSet_{\busIndex,\chargingIndex}}} \frac{\electricityVar_{\busIndex,\chargingIndex}}{\delta_{\busIndex, \chargingIndex}} \leq \, \contractedVar_{\locationIndex} / 60 && \forall \, \locationIndex \in \LocationSet, \, \timeIndex \in \timeSet \label{eqn:cons_cont_lp_6}
\end{align}
The objective \eqref{eqn:lp_cont_obj} again consists of two terms: the electricity costs for charging and the contracted power capacity costs across all charging locations. Constraint \eqref{eqn:cons_cont_lp_1} maintains consistency in the charging levels between two consecutive charging opportunities for each bus. Constraint \eqref{eqn:cons_cont_lp_2} prohibits buses from charging in excess of their battery limits. Constraint \eqref{eqn:cons_cont_lp_3} guarantees that the energy levels in buses are always maintained above a specified minimum threshold. Furthermore, due to \eqref{eqn:cons_cont_lp_4}, buses begin their daily activities with a battery that is fully charged. The maximum energy transferred to a bus within a single charging opportunity is captured in constraint \eqref{eqn:cons_cont_lp_5}. Finally, equation \eqref{eqn:cons_cont_lp_6} defines the station-level maximum contracted capacity, based on the energy supply rate during various charging opportunities, for buses charging simultaneously at the station.

\begin{figure}[H]
    \centering
    \includegraphics[scale=0.5]{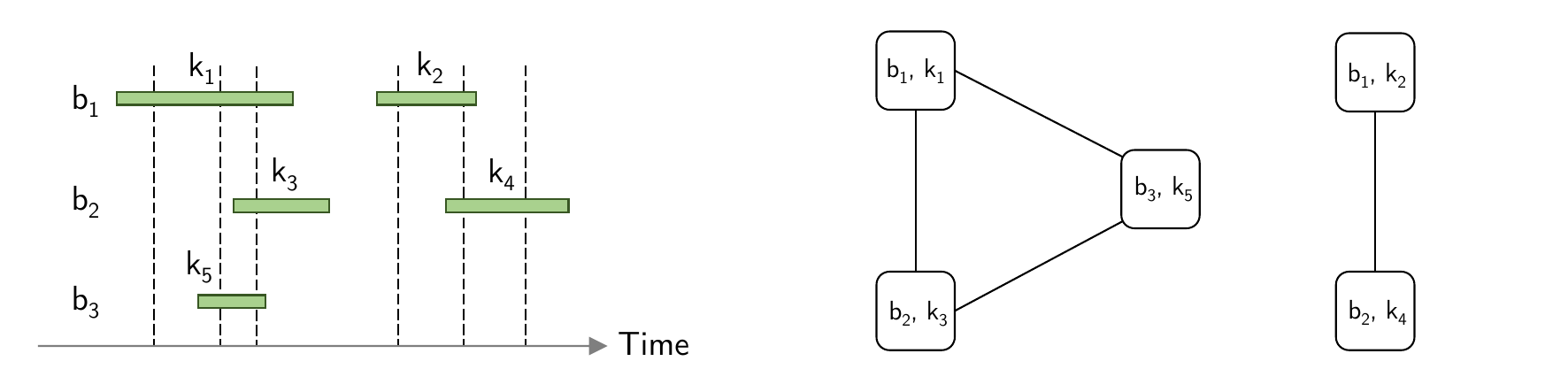}
	\caption{Clique enumeration for populating CSP constraints}
\label{fig:diagram_clique}
\end{figure}

From an implementation perspective, constraint \eqref{eqn:cons_cont_lp_6} need not be added for every time step. Alternatively, for each station, we could identify overlapping charging opportunities for different buses. This allows us to calculate the total rate at which multiple buses are charged simultaneously. To this end, we construct an undirected graph for each charging station $\locationIndex$ with pairs of $(\busIndex,\chargingIndex)$ as vertices, where $\busIndex \in \BusSet, \, \chargingIndex \in \SingleOppSet, \locationIndex_{\busIndex, \chargingIndex} = s$. We add an edge between two vertices if the charging opportunity time windows overlap. Once this graph is constructed, \eqref{eqn:cons_cont_lp_6} can be populated by enumerating all maximal cliques using the Bron-Kerbosch algorithm \citep{bron1973algorithm}. This method of finding cliques is time-efficient since it avoids duplicating constraints for different time-steps. Figure \ref{fig:diagram_clique} presents an example of three buses, $b_1, b_2,$ and $b_3$, at a station, where the charging opportunity time windows are indicated by green bars. The time windows for the pairs $(b_1, k_1)$, $(b_2, k_3)$, and $(b_3, k_5)$ overlap, as do the time windows for the pairs $(b_1, k_2)$ and $(b_2, k_4)$. These overlaps are illustrated in the right panel of the figure. Based on this graph, only two equations are required to represent \eqref{eqn:cons_cont_lp_6} at this station.

The MDVSP is known to be NP-Hard \citep{bertossi1987some}. Unsurprisingly, the CLP-EVSP-CSP, even with surrogate models, scales poorly for real-world instances. To tackle this challenge, we propose two ILS heuristic procedures: one for the CLP-EVSP, which is sequentially followed by a CSP, and another for the joint CLP-EVSP-CSP, which integrates the aforementioned LPs into the search operators. Additionally, we occasionally use the solution to a CLP-CSP where the EVSP decisions are pre-determined, and each charging station's usage is optional. This scenario is modeled by introducing additional binary variables $\LocVar_\locationIndex$ in the CAG CSP formulations, along with a constraint of the form \eqref{eqn:cons_joint_19}.

\section{Iterated Local Search}
\label{sec:localsearch}
We use an ILS framework to overcome the computational challenges of the MILP model and solve CLP-EVSP and CLP-EVSP-CSP efficiently. The procedure's initial solution is generated through the concurrent scheduler algorithm \citep{bodin1978ucost}. This algorithm organizes trips in ascending order of departure times and assigns them to available buses in a greedy manner. Additionally, it activates charging stations at terminal stops when the buses' charge levels fall below certain thresholds. Details on generating initial solutions and charge-level feasibility checks are discussed in Sections \ref{sec:conc_scheduler} and \ref{sec:feasibility}, respectively.

\begin{figure}[t]
    \centering
    \begin{subfigure}{0.42\textwidth}
        \centering
        \includegraphics[width=\textwidth]{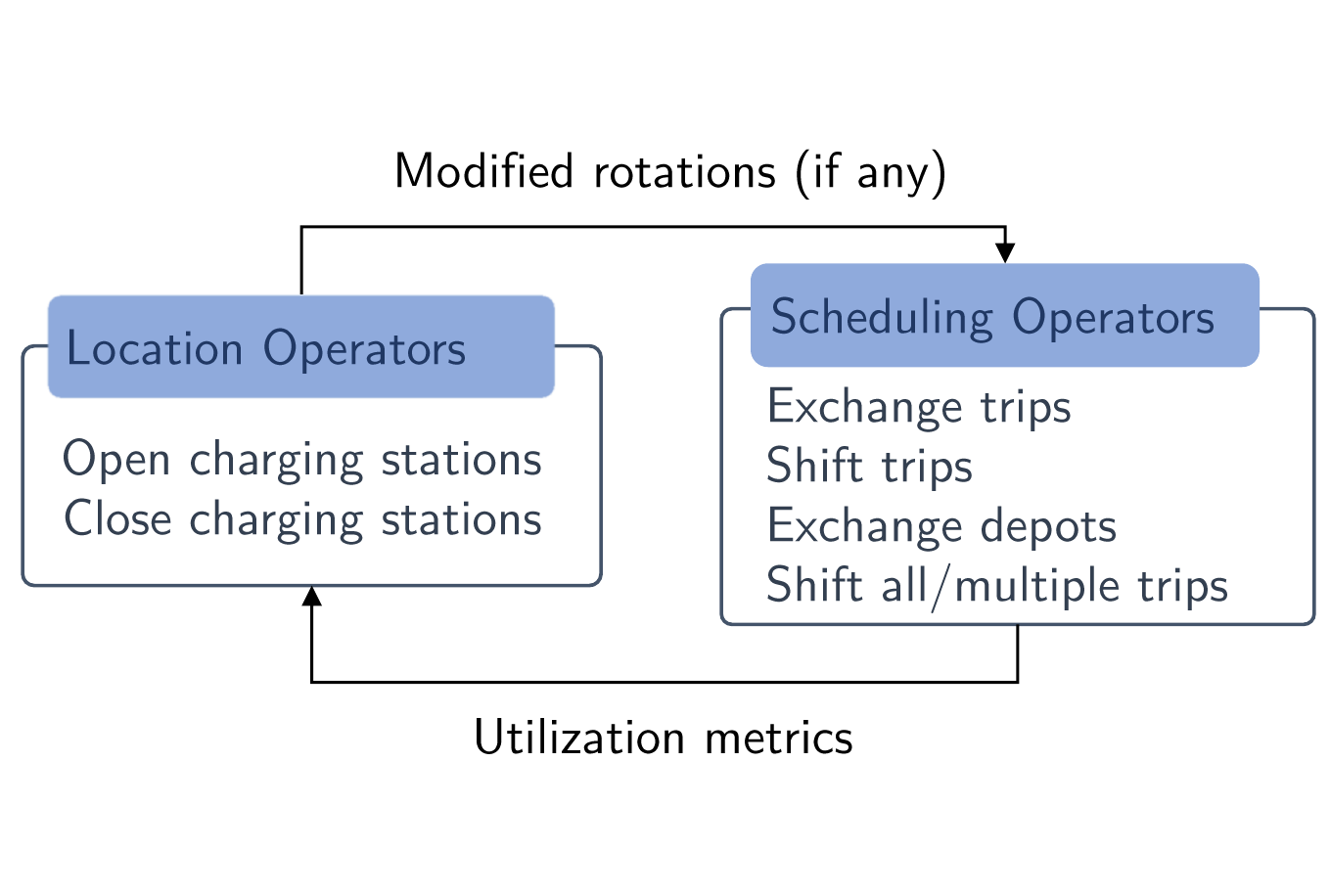}
        \caption{ILS procedure for the CLP-EVSP}
        \label{fig:local_search_summary}
    \end{subfigure}
    \hfill
    \begin{subfigure}{0.57\textwidth}
        \centering
        \includegraphics[width=\textwidth]{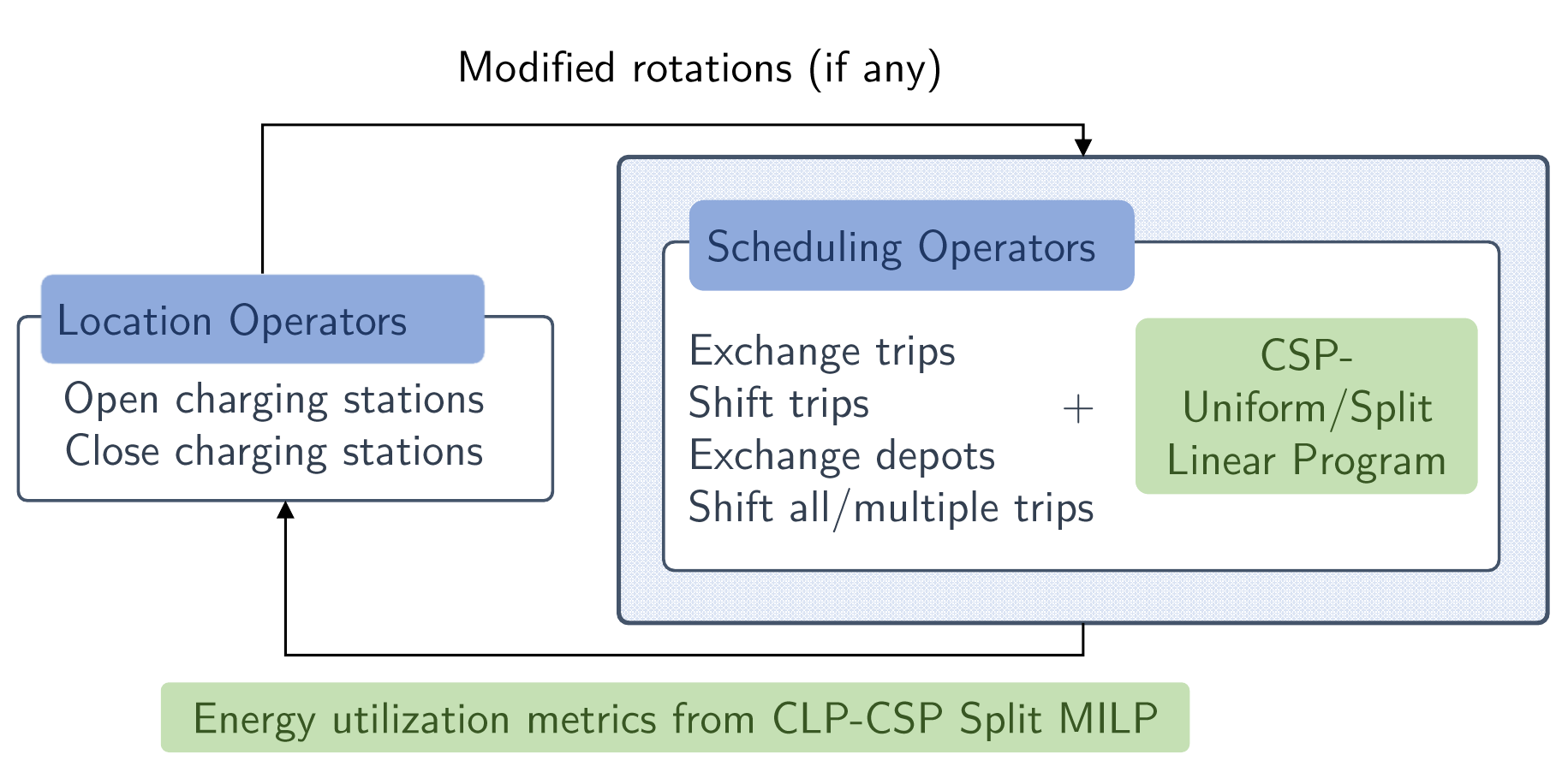}
        \caption{ILS procedure for the CLP-EVSP-CSP}
        \label{fig:local_search_joint_summary}
    \end{subfigure}
    \caption{Iterated Local Search procedures}
\end{figure}

In the sequential model, we first solve the CLP-EVSP using the ILS procedure. At a high level, the ILS for the CLP-EVSP involves two sets of operators: \textit{location operators} and \textit{scheduling operators}. The location operators (refer Section \ref{sec:location_operators}), which can be viewed as the ones making the first-stage decisions, optimize the CLP by evaluating the impact of opening and closing charging stations. While opening charging stations does not compromise the feasibility of existing vehicle rotations, closing them can render rotations charge-infeasible. In such cases, a rotation may have to be split, necessitating the creation of additional rotations. These rotations are subsequently passed to the scheduling operators (see Section \ref{sec:scheduling_operator}), which optimize the EVSP through trip/depot exchanges and shifts. Operating costs are thus reduced either due to minimized deadheading or through the removal of rotations when trips are shifted. The resulting rotations are then used to calculate utilization metrics, which track the time spent by buses at different terminal stops, helping the location operators prioritize decisions on opening and closing charging stations. Subsequently, assuming that the charging station locations and vehicle rotations are fixed, we optimize the CSP by solving the CEE CSP formulation. Figure \ref{fig:local_search_summary} summarizes the steps involved in the ILS process of the sequential model, and Table \ref{tab:notations_ls} describes the additional notation used for the local search methods. Unless stated otherwise, notations used in this section have the same meaning as those defined in the previous sections.

Integrating the MILP model for charge scheduling with the local search for the CLP-CSP, while not inconceivable, is impractical, as the MILP model need to be solved numerous times within the scheduling operators. Hence, we use our novel surrogate CSP formulations, which help find effective upper bounds to the optimal charging-related costs and assist in distinguishing the good exchanges and shifts from the bad ones. Inarguably, employing these surrogate models introduces an optimality gap when compared to the generic CEE version. However, more important than the magnitude is the consistency in the ordering of the CSP costs between the CEE MILP and the surrogate models across different rotations. Since these surrogate models are LPs, they exhibit shorter runtimes, making a strong case for their use within the ILS framework. 

\begin{table}[h]
 \small
  \centering
  \caption{Notation used in ILS methods}
    \begin{tabular}{lp{120mm}}
    \hline
    \textbf{Notation} & \textbf{Description} \\
    \hline
    \textbf{Sets/Lists:} \\
    \hline
     $\busrotationList$ & List of rotations of all the electric buses\\
     $\busrotationb$ & Rotation of $\busIndex^\text{th}$ electric bus present in a list of bus rotations \\
    \hline
    \textbf{Data/Variables:} \\
    \hline
     $\numTripsBusb$ & Number of trips performed by bus $\busIndex$\\
     $\bustrip{\busIndex}{\tripIndexJ}$ & $\tripIndexJ^\text{th}$ trip performed by bus $\busIndex$ \\
     $\bustrip{\busIndex}{0}$ & Starting depot of bus $\busIndex$ \\
     $\bustrip{\busIndex}{\numTripsBusb + 1}$ & Ending depot of bus $\busIndex$ \\
     $\nearestdepotTrip{\tripIndexI}$ & Nearest depot location from the starting bus stop location of trip $\tripIndexI \in \tripSet$ \\
     $\currentUtilization{s}$ & Current utilization of terminal $s$ \\
     $\potentialUtilization{s}$ & Potential utilization of terminal $s$ \\
     $\currentChargeCapacity{s}$ & Optimal power capacity of terminal $s$ from the CLP-CSP split charging formulation \\
     $\objectivefunction{\busrotationList}{\LocationSet}$ & Objective function with or without the CSP cost \\
    \hline 
    \end{tabular}%
  \label{tab:notations_ls}%
\end{table}%

The overarching idea guiding our decisions for integration is straightforward. Operators that frequently invoke the CSP (such as exchanges and shifts) employ the uniform CAG model due to its speed. In contrast, operators that call the CSP model sparingly, such as the objective function evaluations, utilize the split CAG version. In addition, we treat the charging locations suggested by the location operators merely as candidate locations, allowing a CLP-CSP split CAG model to optimize the actual locations within this set. This CLP-CSP solution takes into account the fixed costs of establishing charge scheduling infrastructure and the power capacity costs, and can potentially reduce the number of charging stations. 

A summary of the enhancements to the ILS procedure is shown in Figure \ref{fig:local_search_joint_summary}. The green boxes indicate the CSP-related sub-routines. The candidate locations are set as before, using the first-stage location operators. The EVSP-CSP is then addressed using the rotation scheduling operators with the CSP surrogate models. Finally, utilization metrics, which tracked the time spent by buses at different terminal stops, are now augmented with power capacity variables that guide the search for opening and closing charging stations. In the following subsections, we delve into the details of the individual components of the ILS. We summarize the overall procedure in Section \ref{sec:alg_summary}.

\subsection{Initial Solutions}
\label{sec:conc_scheduler}
The Concurrent Scheduler (CS) algorithm generates an initial feasible solution, see Algorithm \ref{alg:cs}. The solutions to the CLP-EVSP can be represented using a list of bus rotations, $\busrotationList$, and a set of charging locations $\LocationSet \subseteq \candidateLocationSet$. We denote $\busrotationList$ as $[\busrotationList_1, \busrotationList_2, \ldots, \busrotationList_b, \ldots, \busrotationList_{|\BusSet|}]$, where $\busrotationList_{\busIndex}$ is the rotation of the $\busIndex^\text{th}$ bus, containing information on the sequence of service trips, along with its initial and final depot locations. We use the indices $b$, $u$, and $v$ to refer to buses. Assume that the number of trips performed by the $\busIndex^\text{th}$ bus is $\numTripsBusb$. The list $\busrotationList_{\busIndex}$ can then be mathematically written as $\busrotationList_{\busIndex} = [\bustrip{\busIndex}{0}, \bustrip{\busIndex}{1}, \bustrip{\busIndex}{2}, \ldots , \bustrip{\busIndex}{\numTripsBusb}, \bustrip{\busIndex}{\numTripsBusb + 1}]$, where $\bustrip{\busIndex}{\tripIndexK}$ denotes the $\tripIndexK^\text{th}$ trip performed by bus $\busIndex$. The starting and ending depots for bus $\busIndex$ are represented as dummy trips $\bustrip{\busIndex}{0}$ and $\bustrip{\busIndex}{\numTripsBusb + 1}$, respectively. Note that the starting and ending depots can be different. For these rotations, the corresponding decision variables in the MILP models are $\arcVar{\bustrip{\busIndex}{0}} {\bustrip{\busIndex}{1}}{\busIndex} = 1, \arcVar{\bustrip{\busIndex}{1}} {\bustrip{\busIndex}{2}}{\busIndex} = 1, \ldots , \arcVar{\bustrip{\busIndex}{\numTripsBusb}} {\bustrip{\busIndex}{\numTripsBusb + 1}}{\busIndex} = 1$, and for the charging station locations, $\LocVar_{\locationIndex} = 1$ for all $ \locationIndex \in \LocationSet$. 

\begin{algorithm}[H]
\caption{\textsc{ConcurrentScheduler (CS)}}
\label{alg:cs}
\KwIn{$\tripSet$}
\KwOut{$\busrotationList, \LocationSet$}

$\busrotationList, \LocationSet \gets \phi, \phi$\;
Rearrange $\tripSet$ in the ascending order of the trip start times\;
Pick the first trip $\tripIndexI$ from the sorted list of $\tripSet$\;
$\busrotationList_1 \gets [\nearestdepotTrip{\tripIndexI}, \tripIndexI, \nearestdepotTrip{\tripIndexI}]$ and add $\busrotationList_1$ to $\busrotationList$\;

\tcp{\textsf{Insert trips to existing rotations}}
\For{$\tripIndexJ = 2$ \KwTo $\numTrips$}{
    $\tripIndexI_{\tripIndexJ} \gets$ The $\tripIndexJ^\text{th}$ element from the sorted list of $\tripSet$\;
    
    \For{$b = 1$ \KwTo $|\busrotationList|:$ $(\bustrip{b}{\numTripsBus{b}}, \tripIndexI_\tripIndexJ) \in \compatibleSet$}
        {
            $\busrotationAnotherList_b \leftarrow \busrotationList_b$ and insert $\tripIndexI_{\tripIndexJ}$ to the list of service trips in $\busrotationAnotherList_b$ \;
            $\textsf{insertTrip}, \LocationAnotherSet \gets \textsc{IsRotationChargeFeasible}(\busrotationAnotherList_b, \LocationSet)$\;
            \If{$\mathsf{insertTrip}$}{
                $\busrotationList_b \gets \busrotationAnotherList_b$ and $\LocationSet \gets \LocationAnotherSet$\;
                \textbf{break}
            }
        }
    \tcp{\textsf{Create new rotations}}
    \If{$\tripIndexI_\tripIndexJ$ is not assigned to any of the buses used so far}{
        Use a new bus $|\busrotationList| + 1$, update its rotation to $[\nearestdepotTrip{\tripIndexI_\tripIndexJ}, \tripIndexI_\tripIndexJ, \nearestdepotTrip{\tripIndexI_\tripIndexJ}]$, and add it to $\busrotationList$\;
    }
}
\end{algorithm}

Initially, we assume no open charging locations. The algorithm starts by sorting the service trips in the ascending order of their start times (line 2). It creates a new rotation by assigning the first trip $\tripIndexI$ to a bus by starting a bus from the depot nearest to trip $\tripIndexI$, $\nearestdepotTrip{\tripIndexI}$ (lines 3--4). The algorithm then checks if the next trip $\tripIndexI_j$ from the sorted list can be assigned to a used bus (lines 5--12). When multiple rotations can accommodate the trip, priority is given to the bus rotation, which was created first. The time compatibility for inserting the trip is checked in line 7, and a copy of the updated bus rotation is created in line 8. Additionally, the energy level feasibility of the bus rotation is checked in line 9 using the \textsc{IsRotationChargeFeasible} function. This function returns a Boolean value and an updated list of charging stations, which can be opened at any terminal stop. It is similar to Algorithm \ref{alg:feasibility_check}, described in Section \ref{sec:feasibility}. The charging locations are modified using the CAG strategy, where opening a charging facility is prioritized at the end stops of trips. If the trip insertion is feasible, the bus rotation and charging locations are updated in line 11, and the algorithm proceeds to insert the next trip. If not, the algorithm creates a new bus rotation and assigns the trip to it (lines 13--14).

\subsection{Charge Feasibility Check}
\label{sec:feasibility}
Checking the charge feasibility of a bus for a given set of charging stations is central to the local search procedure and is performed by Algorithm \ref{alg:feasibility_check}. Several sub-routines use variants of this algorithm, which are briefly described in their respective sections. The charge level at terminals is monitored using the $\levelVar$ variable. Since buses depart the depot fully charged, they are not permitted to be charged at the start location of the first trip. Lines 2--3 ensure that the bus can deadhead to the start terminal of the first trip and successfully complete it. In cases where a rotation is infeasible, it is sometimes necessary to identify the last trip in its sequence that can be completed without running out of charge. This information is tracked using the \textsf{feasibleTillTrip} variable. This variable is updated in line 5, which checks if the rotation can be terminated early and the end depot can be reached without letting the charge levels fall below the minimum threshold. 

Checking whether the charge level $\levelVar$ falls below $\MinEnergy$ must be carefully carried out after every trip, especially considering the availability of charging opportunities. We declare the bus rotation as infeasible if charging at the maximum possible rate cannot meet the energy requirements of its trips. This leads to three scenarios under the CAG strategy. In the first scenario (lines 6--8), the end terminal is a charging station, allowing the bus to either charge to the maximum possible level, i.e., recharge at the maximum rate possible during the idle time between trips, denoted as $\chargingRate \delta_{\bustrip{\busIndex}{j},\bustrip{\busIndex}{j+1}}$. Subsequently, the bus may need to deadhead to the start location of the next trip, after which the charge level is compared with $\MinEnergy$. 

In the second scenario (lines 9--12), the bus can only charge after deadheading. Therefore, the charge levels are first checked for infeasibility and then updated with the maximum possible energy that can be provided within its available idle time. In the third scenario (lines 13--15), where there is no charging opportunity, the charge levels are updated following the deadheading trip. Line 16 accounts for the energy required for the next trip, and line 20 addresses the scenario of deadheading to the end depot after the last trip. 

\begin{algorithm}[H]
\caption{\textsc{IsRotationChargeFeasible}}\label{alg:feasibility_check}
\KwIn{$\busrotationList_b, \LocationSet$}
\KwOut{True/False}

\textsf{feasibleTillTrip} $\gets -1$\;
    $\levelVar \gets \BatteryCap - (\energyDeadhead{\bustrip{\busIndex}{0}}{\bustrip{\busIndex}{1}} - \energyTrip{\bustrip{\busIndex}{1}})$\;
    \lIf{$\levelVar < \MinEnergy$}{
        \Return \textbf{false}}
    \For{$j = 1$ \KwTo $(\numTripsBusb - 1)$}{
        \lIf{$\levelVar - \energyDeadhead{\bustrip{\busIndex}{j}}{\bustrip{\busIndex}{\numTripsBusb+1}} \geq \MinEnergy$}{$\textsf{feasibleTillTrip} \gets j$}
\tcp{\textsf{Current trip end is a charging station}}
        \If{$\bustrip{\busIndex}{j, end} \in \LocationSet$}{
        $\levelVar \gets \min\{\BatteryCap, l + \chargingRate \delta_{\bustrip{\busIndex}{j}, \bustrip{\busIndex}{j+1}}\} - \energyDeadhead{\bustrip{\busIndex}{j}}{\bustrip{\busIndex}{j+1}}$\;
        \lIf{$\levelVar < \MinEnergy$}{$\textsf{isNextTripFeasible} \gets$ \textbf{false}}
            }
        \tcp{\textsf{Only next trip start is a charging station}} 
        \If{$\bustrip{\busIndex}{j, end} \notin \LocationSet$, $\bustrip{\busIndex}{j+1, start} \in \LocationSet$}{
        $\levelVar \gets \levelVar - \energyDeadhead{\bustrip{\busIndex}{j}}{\bustrip{\busIndex}{j+1}}$\;
        \lIf{$\levelVar < \MinEnergy$}{$\textsf{isNextTripFeasible} \gets$ \textbf{false}}
        $\levelVar \gets \min\{\BatteryCap, l + \chargingRate \delta_{\bustrip{\busIndex}{j}, \bustrip{\busIndex}{j+1}}\}$\; 
        }
        \tcp{\textsf{Current trip end and next trip start are not charging stations}}
        \If{$\bustrip{\busIndex}{j, end} \notin \LocationSet, \bustrip{\busIndex}{j+1, start} \notin \LocationSet$}{
        $\levelVar \gets \levelVar - \energyDeadhead{\bustrip{\busIndex}{j}}{\bustrip{\busIndex}{j+1}}$\;
        \lIf{$\levelVar < \MinEnergy$}{$\textsf{isNextTripFeasible} \gets$ \textbf{false}}
        }

$\levelVar \gets \levelVar - \energyTrip{\bustrip{\busIndex}{j+1}}$\;

\lIf{$\levelVar \geq \MinEnergy$}{$\textsf{isNextTripFeasible} \gets$ \textbf{true}}

\If{$\mathbf{not}~\mathsf{isNextTripFeasible}$}{

            \Return \textbf{false}
        }
        }

    $\levelVar \gets \levelVar - \energyDeadhead{\bustrip{\busIndex}{\numTripsBusb}}{\bustrip{\busIndex}{\numTripsBusb+1}}$\;
    
    \If{$\levelVar < \MinEnergy$}{
        \Return \textbf{false}
    }
    {$\textsf{feasibleTillTrip} \gets \numTripsBusb$}\;
\Return \textbf{true}\;
\end{algorithm}

\begin{figure}[H]
\centering
\begin{subfigure}{\textwidth}
  \centering
  \includegraphics[width=0.7\linewidth]{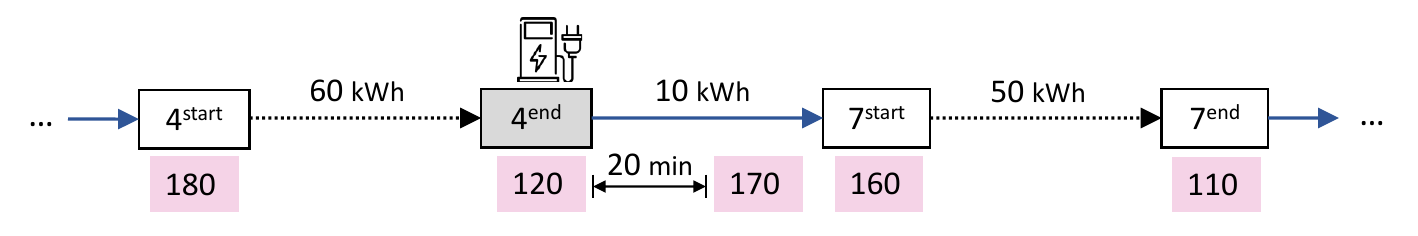}  
  \caption{Charging station is at the end terminal of the current trip}
  \label{fig:sub1}
\end{subfigure}

\begin{subfigure}{\textwidth}
  \centering
  \includegraphics[width=0.7\linewidth]{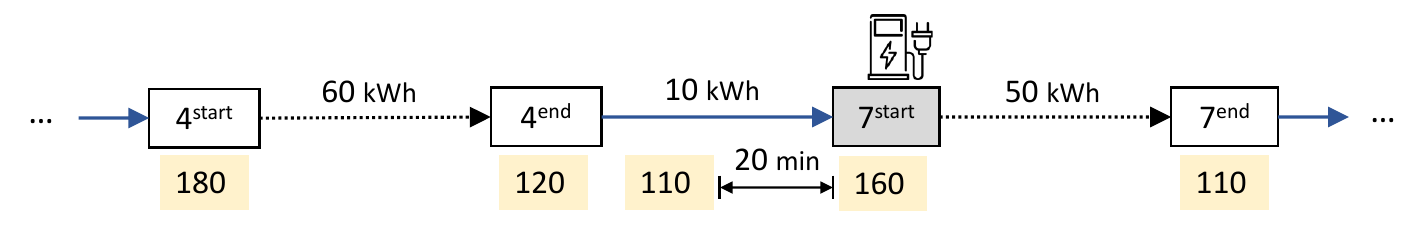}
  \caption{Charging station is at the start terminal of the next trip}
  \label{fig:sub2}
\end{subfigure}

\begin{subfigure}{\textwidth}
  \centering
  \includegraphics[width=0.7\linewidth]{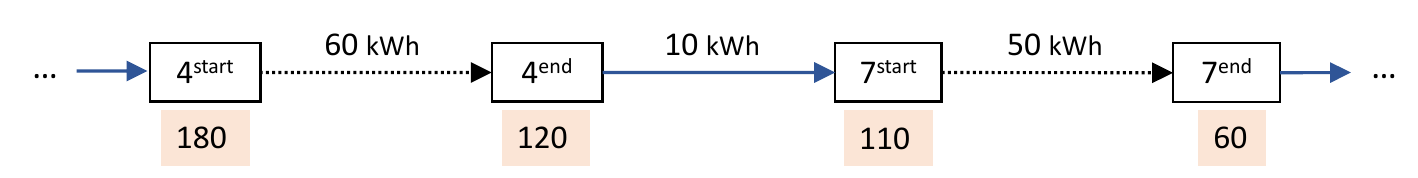}
  \caption{There are no charging stations at the end terminal of the current trip and the start terminal of the next trip}
  \label{fig:sub3}
\end{subfigure}
\caption{Scenarios for charging between trips with a maximum charging rate $\chargingRate = 2.5$ kWh/min}
\label{fig:charging_scenarios}
\end{figure}

Figure \ref{fig:charging_scenarios} illustrates the charge level updates between two consecutive trips ($4$ and $7$) of a bus rotation. The numbers on the arcs indicate the energy required for each trip or for deadheading. The values displayed in the colored boxes represent the charge levels under different scenarios. Assuming that charging occurs throughout the entire idle time at a maximum rate of 2.5 kWh/min, the charge levels increase by 50 kWh at each charging opportunity. In the following sections, we use variants of this function, such as the \textsc{AreRotationChargeFeasible}, in which this function is applied in a loop across different rotations. The notations for these functions are slightly abused to spare the reader from a description of overloaded versions of the function.

\subsection{Optimizing Rotations}
\label{sec:scheduling_operator}
We optimize rotations through simple local search moves that involve exchanging and shifting service trips or exchanging depots. Each iteration of this procedure (see Algorithm \ref{alg:opt_rotations})  alters the rotations in a greedy way by selecting a move that yields the maximum objective function improvement. This iterative process yields solutions that we call \textit{exchange optimal} and \textit{shift optimal}, i.e., solutions that cannot be improved through exchanging two trips or shifting a trip. One can view this as being similar to 2-opt moves in the traveling salesman problem. Exchanging and shifting multiple trips can also be performed, but this comes at a higher computational cost as the number of possibilities explodes. In this context, methods such as improvement graphs \citep{ahuja2002survey} are not directly applicable since the savings are not additive. Instead, we use an \textsc{OptimizeMultipleShifts} subroutine on a limited set of rotations involving few trips (based on a user-defined threshold) to determine if multiple trips can be shifted elsewhere.

\begin{algorithm}[h]
\caption{\textsc{OptimizeRotations}}
\label{alg:opt_rotations}

\KwIn{$\busrotationList, \LocationSet$}
\KwOut{$\busrotationList$}

$\textsf{improvement} \gets \infty$\;
$\bestbusrotationList \gets \busrotationList$\;
\While{$\mathsf{improvement} > 0$}{
    $\busrotationList \gets \bestbusrotationList$\;
    $\bestbusrotationList \gets \textsc{ApplyBestImprovement}(\busrotationList,\LocationSet)$\;
    $\textsf{improvement} \gets f(\busrotationList, \LocationSet) - f(\bestbusrotationList, \LocationSet)$\;
}
$\busrotationList \gets \textsc{OptimizeMultipleShifts}(\busrotationList, \LocationSet)$\;
\end{algorithm}

Rotations are updated within each iteration using the \textsc{ApplyBestImprovement} sub-routine (Algorithm \ref{alg:apply_best_improvement}). The best exchange and shift across all service trips are first evaluated (lines 1--2), and the better of the two is chosen (lines 8--11). The subroutines \textsc{Exst} and \textsc{Sst} are described in \ref{sec:appendix_optimize_rotations}. When service trip exchanges and shifts do not yield any benefit, we try exchanging depots (lines 3--7). Note that buses do not need to return to the same depot at the end of the day as long as the distribution of buses across depots remains the same at the start and end of operations. This allows the service provider to repeat schedules without needing to re-balance the buses. 
Trip exchanges and shifts are the most time-consuming parts of the local search process, particularly when the CSP is integrated. Hence, when dealing with networks with a large number of trips, we use a hybrid version of this algorithm to keep the computation tractable. In this hybrid method, a candidate pool of exchanges and shifts is shortlisted based on the savings in deadheading, and the CSP is jointly solved only for new rotations derived from this pool.

\begin{algorithm}[h]
\caption{\textsc{ApplyBestImprovement}}
\label{alg:apply_best_improvement}
\KwIn{$\busrotationList, \LocationSet$}
\KwOut{$\busrotationList$}

$\ExchangetriprotationList, \textsf{exchangeSavings} \gets \textsc{Exst}(\busrotationList, \LocationSet)$\;
$\ShifttriprotationList, \textsf{shiftSavings} \gets \textsc{Sst}(\busrotationList, \LocationSet)$\;

\tcp{\textsf{Exchange depots}}
\If{$\mathsf{exchangeSavings} \leq 0 \mathbf{~and~} \mathsf{shiftSavings} \leq 0$}{
    $\ExchangedepotrotationList, \textsf{depotExchangeSavings} \gets  \textsc{Exd}(\busrotationList, \LocationSet)$\;
    \If{$\mathsf{depotExchangeSavings} > 0$}{
        $\busrotationList \gets \ExchangedepotrotationList$\;
    }
    \Return\;
}
{
\tcp{\textsf{Exchange or shift service trips}}
    \eIf{$\mathsf{exchangeSavings} > \mathsf{shiftSavings}$}{
        $\busrotationList \gets \ExchangetriprotationList$\;
    }{
        $\busrotationList \gets \ShifttriprotationList$\;
    }
}
\end{algorithm}

\subsection{Optimizing Stations}
\label{sec:location_operators}
\subsubsection{Utilization Metrics}
To prioritize the opening and closing of charging stations, we track the time spent by buses at different charging stations. If several buses are idle at stops that currently are not charging stations, opening a station at that location can potentially reduce charging costs. Likewise, if buses do not spend much time at locations where charging stations are currently open, closing them would help recover the fixed costs associated with charging infrastructure. 

\begin{algorithm}[H]
\caption{\textsc{UpdateUtilization}}
\label{alg:update_utilization}

\KwIn{$\busrotationList, \LocationSet$}
\KwOut{$\LocationSet$}
\tcp{\textsf{Initialize current and potential utilizations}}
{
    $\currentUtilization{s}, \potentialUtilization{s} \gets 0, 0 \, \forall \, s \in \LocationSet$\;
}

\For{$\{b \in \busrotationList : b \text{ requires charging}\}$}{
    \For{$j = 1$ \KwTo $\numTripsBusb - 1$}{
    \tcp{\textsf{The current trip's end terminal is a charging station}}
        \If{$\bustrip{\busIndex}{j, end} \in \LocationSet$}{
            $\currentUtilization{\bustrip{\busIndex}{j, end}} \gets \currentUtilization{\bustrip{\busIndex}{j, end}} + \idleTime_{\bustrip{\busIndex}{j}, \bustrip{\busIndex}{j+1}}$ 
            }
        \tcp{\textsf{Only the next trip's start terminal is a charging station}} 
        \If{$\bustrip{\busIndex}{j, end} \notin \LocationSet$, $\bustrip{\busIndex}{j+1, start} \in \LocationSet$}{
            $\currentUtilization{\bustrip{\busIndex}{j+1, start}} \gets \currentUtilization{\bustrip{\busIndex}{j+1, start}} + \idleTime_{\bustrip{\busIndex}{j}, \bustrip{\busIndex}{j+1}}$  \\
            $\potentialUtilization{\bustrip{\busIndex}{j, end}} \gets \potentialUtilization{\bustrip{\busIndex}{j, end}} + \idleTime_{\bustrip{\busIndex}{j}, \bustrip{\busIndex}{j+1}}$ 
        }
        \tcp{\textsf{The current trip's end and next trip's start terminals are not charging stations}}
        \If{$\bustrip{\busIndex}{j, end} \notin \LocationSet, \bustrip{\busIndex}{j+1, start} \notin \LocationSet$}{
            $\potentialUtilization{\bustrip{\busIndex}{j+1, start}} \gets \potentialUtilization{\bustrip{\busIndex}{j+1, start}} + \idleTime_{\bustrip{\busIndex}{j}, \bustrip{\busIndex}{j+1}}$\;
            $\potentialUtilization{\bustrip{\busIndex}{j, end}} \gets \potentialUtilization{\bustrip{\busIndex}{j, end}} + \idleTime_{\bustrip{\busIndex}{j}, \bustrip{\busIndex}{j+1}}$\;  
        }
    }
}

\tcp{\textsf{Close charging terminals with zero utilization}}
    \For{$s \in \LocationSet : \currentUtilization{s} = 0$}{
    $\LocationSet \gets \LocationSet \setminus \{s\}$ \;}
\end{algorithm}

Algorithm \ref{alg:update_utilization} is designed to update the utilization of various locations. We keep track of two types of utilization metrics at each terminal $s$: \textit{current utilization} $\currentUtilization{s}$ and \textit{potential utilization} $\potentialUtilization{s}$. The current utilization simply stores the cumulative time spent by all buses at a charging station. Potential utilization, on the other hand, captures the total idle time spent by all the buses should a charging station be opened at $s$. Finally, the algorithm also removes any charging stations from the set of locations that have zero current utilization in lines 12--13.  

Initially, the algorithm sets both the current and potential utilizations of each location to zero. For each trip made by a bus, it checks if the ends of a layover have a charging station and works in three scenarios just like Algorithm \ref{alg:feasibility_check}. In case 1, if the end terminal of the current trip is a charging station, the algorithm increments the current utilization of that station by the idle time between the current and the next trip. Note that the utilization stats are dependent on the current charging station configuration since we assume that operations are carried out under the CAG policy. Hence, an open charging station at the start of a trip will not be utilized if there is already a charging station at the end terminal of its previous trip. In scenario 2, if the end terminal is not a charging station but the start terminal of the next trip is, it increments the current utilization of the next trip's start terminal and the potential utilization of the current trip's end terminal by the idle time. Finally, in scenario 3, if neither end of a layover is a charging station, only the potential utilization for both terminals is incremented by the idle time.

\subsubsection{Opening Charging Stations}
Algorithm \ref{alg:open_stations} evaluates the impact of adding new charging stations on the overall cost. It takes as input a list of charging stations to open $\LocationSet^{open}$ and checks for improvement in the rotations (and charge scheduling in the case of the joint model). That is, the fixed costs of opening the station must be offset by a reduction in the vehicle and charge scheduling costs. Charging stations that are not needed are then removed using the \textsc{UpdateUtilization} function (in the case of the sequential model) and the CLP-CSP model (in the joint model case). If the new objective value is not better, no changes are made to the set of charging stations and bus rotations (lines 7--9).

\begin{algorithm}[h]
\caption{\textsc{OpenChargingStation}}
\label{alg:open_stations}

\KwIn{$\busrotationList, \LocationSet, \LocationSet^{open}$}
\KwOut{$\busrotationList, \LocationSet$}
$\LocationSet^{temp} \gets \LocationSet \cup \{\LocationSet^{open}\} $\;
$\busrotationAnotherList \gets \textsc{OptimizeRotations}(\busrotationList, \LocationSet^{temp})$\;
\If{the sequential model is solved}{
$\LocationSet^{temp} \gets \textsc{UpdateUtilization}(\busrotationList, \LocationSet^{temp})$}

\If{if the joint model is solved}{Solve CLP-CSP model (Section \ref{sec:surrogate_lp}), remove unused stations from $\LocationSet^{temp}$}
\If{$f(\busrotationList^{temp}, \LocationSet^{temp}) < f(\busrotationList, \LocationSet)$}{
    $\busrotationList \gets \busrotationAnotherList$\;
    $\LocationSet \gets \LocationSet^{temp}$
}
\end{algorithm}
\subsubsection{Closing Charging Stations and Splitting Trips}
We also test the effect of closing charging stations in a similar way, except that closing charging stations affects charge feasibility, and hence, we allow rotations to be split. Algorithm \ref{alg:close_charge_station} takes as input the specific charging station that is considered for closure and identifies a subset of buses that are affected by the closure of the specified charging station under the CAG strategy (lines 2--6). The algorithm then temporarily removes the specified charging station from the location set and checks if the bus rotations are still feasible. If the rotations are not feasible, the algorithm splits trips at locations from where the bus can return to the depot without falling below the lower threshold of charge using the \textsc{AreRotationsChargeFeasible} function. This function, a variant of Algorithm \ref{alg:feasibility_check}, returns a Boolean value and a set of potentially modified rotations. It yields `true' if closing the station does not result in a bus being stranded due to inadequate charge. This is determined after recursively splitting rotations and checking their charge feasibility. Splitting trips will require new buses and can increase the fixed costs of buses. However, this can be compensated by better bus-to-trip assignments. Hence, the bus rotations are optimized using the \textsc{OptimizeRotations} sub-routine in line 11. If there is an improvement, the station is kept closed, and the algorithm terminates.

\begin{algorithm}[h]
\caption{\textsc{CloseChargingStation}}
\label{alg:close_charge_station}

\KwIn{$\busrotationList, \LocationSet, s^{close}$}
\KwOut{$\busrotationList, \LocationSet$}

\tcp{\textsf{Identify the rotations affected by the terminal closure}}
$\textsf{affectedRotations} \gets \emptyset$\;
\For{$\{b \in \busrotationList : b \text{ requires charging}\}$}{
\For{$j = 1$ \KwTo $\numTripsBusb - 1$}{
    \If{$(\bustrip{\busIndex}{j, end} = s^{close}) \mathbf{~or~} (\bustrip{\busIndex}{j+1, start} = s^{close} \mathbf{~and~} \bustrip{\busIndex}{j, end} \notin \LocationSet)$}{
        $\textsf{affectedRotations} \gets \textsf{affectedRotations} \cup \{\busIndex \}$\;
        \textbf{break}\;
    }
}
}

$\LocationSet^{temp} \gets \LocationSet \setminus \{s^{close}\} $\;

\tcp{\textsf{Check for charge feasibility after terminal closure and split trips if necessary}}
$\textsf{closureFeasible}, \busrotationAnotherList \gets \textsc{AreRotationsChargeFeasible}(\busrotationList, \LocationSet^{temp}, \textsf{affectedRotations})$ \;
\If{$\mathbf{not~} \mathsf{closureFeasible}$}{
\Return\;
}

\tcp{\textsf{Optimize rotations after terminal closure}}
$\busrotationAnotherList \gets \textsc{OptimizeRotations}(\busrotationAnotherList, \LocationSet^{temp})$\;

\If{$f(\busrotationList^{temp}, \LocationSet^{temp}) < f(\busrotationList, \LocationSet)$}{
    $\busrotationList \gets \busrotationAnotherList$\;
    $\LocationSet \gets \LocationSet^{temp}$
}
\end{algorithm}

\subsection{Summary}
\label{sec:alg_summary}
The overall ILS procedure is summarized in Algorithm \ref{alg:fleet_optimization}. The CS algorithm initalizes bus rotations (line 1). The \textsc{OptimizeMultipleShifts} operator (line 2) then optimizes these rotations by efficiently reallocate trips, focusing particularly on rotations with fewer trips. We present the pseudocode for this operator in \ref{sec:appendix_multiple_shift}.

\begin{algorithm}[H]
\caption{\textsc{OptimizeEBusFleet}}
\label{alg:fleet_optimization}

\KwIn{$\tripSet$}
\KwOut{$\busrotationList, \LocationSet$}

\tcp{\textsf{Initialize rotations and locations}}
$\busrotationList, \LocationSet \gets \textsc{ConcurrentScheduler}(\tripSet)$\;

\tcp{\textsf{Optimize rotations that serve few trips}}
$\busrotationList \gets \textsc{OptimizeMultipleShifts}(\busrotationList, \LocationSet)$

\tcp{\textsf{Open charging stations and optimize rotations}}
$\LocationSet \gets \textsc{UpdateUtilization}(\busrotationList, \LocationSet)$\;
$\LocationSet^{open} \gets$ Sort terminals not in $\LocationSet$ for which $\potentialUtilization{\locationIndex} > 0$ in ascending order \; 
$\busrotationList, \LocationSet \gets \textsc{OpenChargingStations}(\busrotationList, \LocationSet, \LocationSet^{open})$\;
\tcp{\textsf{Close charging stations and optimize rotations}}

$\LocationSet^{close} \gets$ Select all terminals in $\LocationSet$\; 

\SetKwRepeat{Do}{do}{while} 
\Do{$\LocationSet^{close} \neq \emptyset$}{

\If{the sequential model is solved}{
$\LocationSet \gets \textsc{UpdateUtilization}(\busrotationList, \LocationSet)$}

\If{if the joint model is solved}{Solve CLP-CSP model (Section \ref{sec:surrogate_lp}), remove unused stations from $\LocationSet$, and update the optimal power capacities $\currentChargeCapacity{\locationIndex}$}

Remove terminal $s$ from $\LocationSet^{close}$ for which $\currentUtilization{\locationIndex}=0$ (sequential model) or $\currentChargeCapacity{\locationIndex}=0$ (joint model)
    
Find terminal $s^{close}$ with minimum $\currentUtilization{\locationIndex}$ (sequential model) or $\currentChargeCapacity{\locationIndex}$ (joint model) and remove it from $\LocationSet^{close}$\;

$\busrotationList, \LocationSet \gets \textsc{CloseChargingStation}(\busrotationList, \LocationSet, s^{close})$ \;
}
Find the optimal charge schedule using the CSP MILP (Section \ref{sec:csp})
\end{algorithm}

Subsequently, the algorithm calculates utilization metrics and opens charging stations at locations that exhibit positive potential utilization (see Algorithm \ref{alg:open_stations}). This process can also be executed sequentially for each station, selecting a station with the highest potential utilization $\potentialUtilization{s}$. These utilization metrics can then be updated after each iteration after further optimizing the rotations. 

From lines 7 to 15, the algorithm evaluates the impact of sequentially closing charging stations, optimizing rotations during this process, and addressing any charge infeasibilities by splitting rotations. The stations currently in operation are tracked in $\LocationSet^{close}$. In the sequential model, the algorithm updates utilizations and identifies the least utilized station for closure. For joint models, the CLP-CSP split formulation is solved, and station power capacities are utilized to prioritize closures. Stations exhibiting zero utilization or power capacity are also removed from the potential closure list.

Ultimately, after determining the optimal configuration of charging stations and rotations, the algorithm solves the MILP CSP formulation described in Section \ref{sec:csp} to optimize charge scheduling costs. Although we explored the possibility of swapping an open charging station with a closed one, the benefits were negligible in our tested datasets and did not justify the additional computational time.

\section{Experiments and Results}
\label{sec:results}
We demonstrate the benefits of our framework on 14 real-world bus transit networks by solving both the joint CLP-EVSP-CSP and the sequential CLP-EVSP with a subsequent CSP model. The MILP model, outlined in Section \ref{sec:cl_evsp_csp}, could only be solved for small toy networks using CPLEX. The local search subroutines were implemented in C++ using the GCC compiler version 11.4.0. LPs were solved using CPLEX, and the Boost library was utilized for the Bron-Kerbosch algorithm within the uniform charging LPs as CSP surrogates. 
All experiments were conducted on a Dell Inc. Precision 7875 Tower equipped with an AMD Ryzen Threadripper Pro 7975wx CPU @ 4.0 GHz and 128 GB of RAM. OpenMP was employed with $32$ cores to run the exchange and shift operators in parallel. The source codes are available at \href{https://github.com/transnetlab/clp-evsp-csp.git}{github.com/transetlab/clp-evsp-csp}.

\subsection{Performance Evaluation of ILS}

\textbf{Datasets and Parameters:}
The GTFS data for all networks were collected from \href{https://transitfeeds.com/}{transitfeeds.com}. For deadheading trips, the shortest path distances between stops are calculated using Open Street Maps (\href{https://www.openstreetmap.org/}{OSM}) and converted into time, assuming a constant bus speed of 30 km/h. This speed is consistent with the bus speeds implied by the trip schedules for various networks. Table \ref{tab:datasets} provides a summary of the network data. The origins and destinations of different routes determine candidate charging locations. Considering the number of trips originating or terminating at each location, a maximum of five depots are chosen from this subset. The initial CS solution determines the column ‘Depots used’ based on the number of depots assigned to rotations.

\begin{figure}[H]
  \centering
  \begin{subfigure}{0.195\textwidth}
    \includegraphics[width=\linewidth]{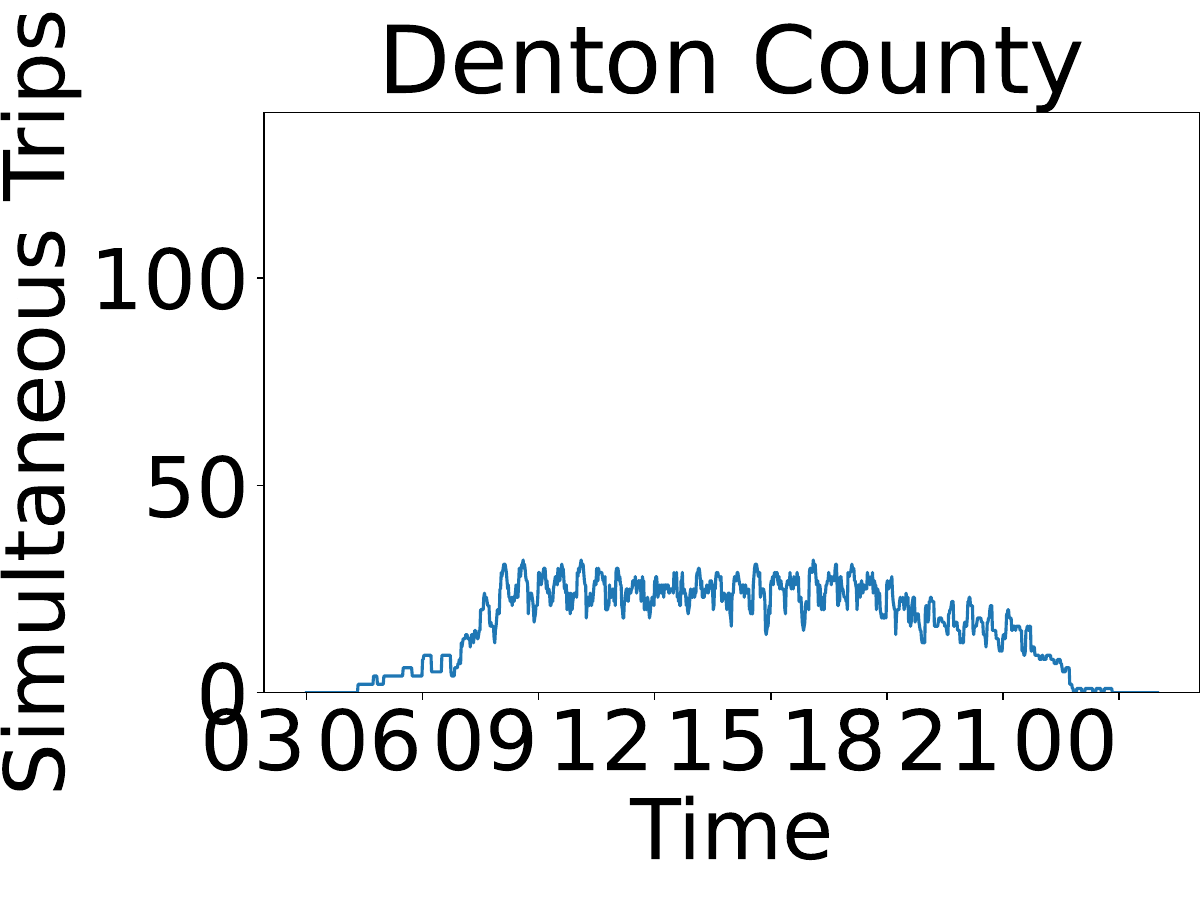}
  \end{subfigure}
  \hfill
  \begin{subfigure}{0.195\textwidth}
    \includegraphics[width=\linewidth]{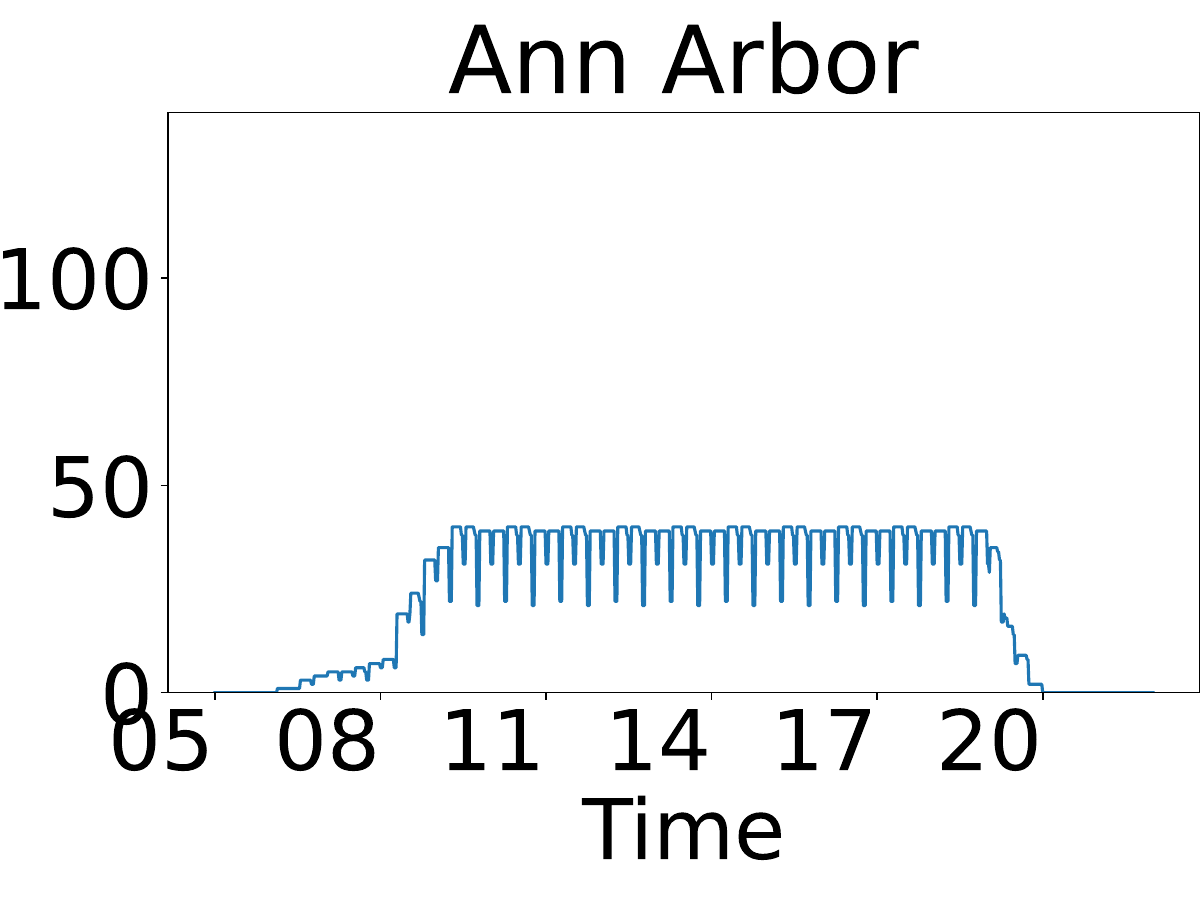}
  \end{subfigure}
    \hfill
  \begin{subfigure}{0.195\textwidth}
    \includegraphics[width=\linewidth]{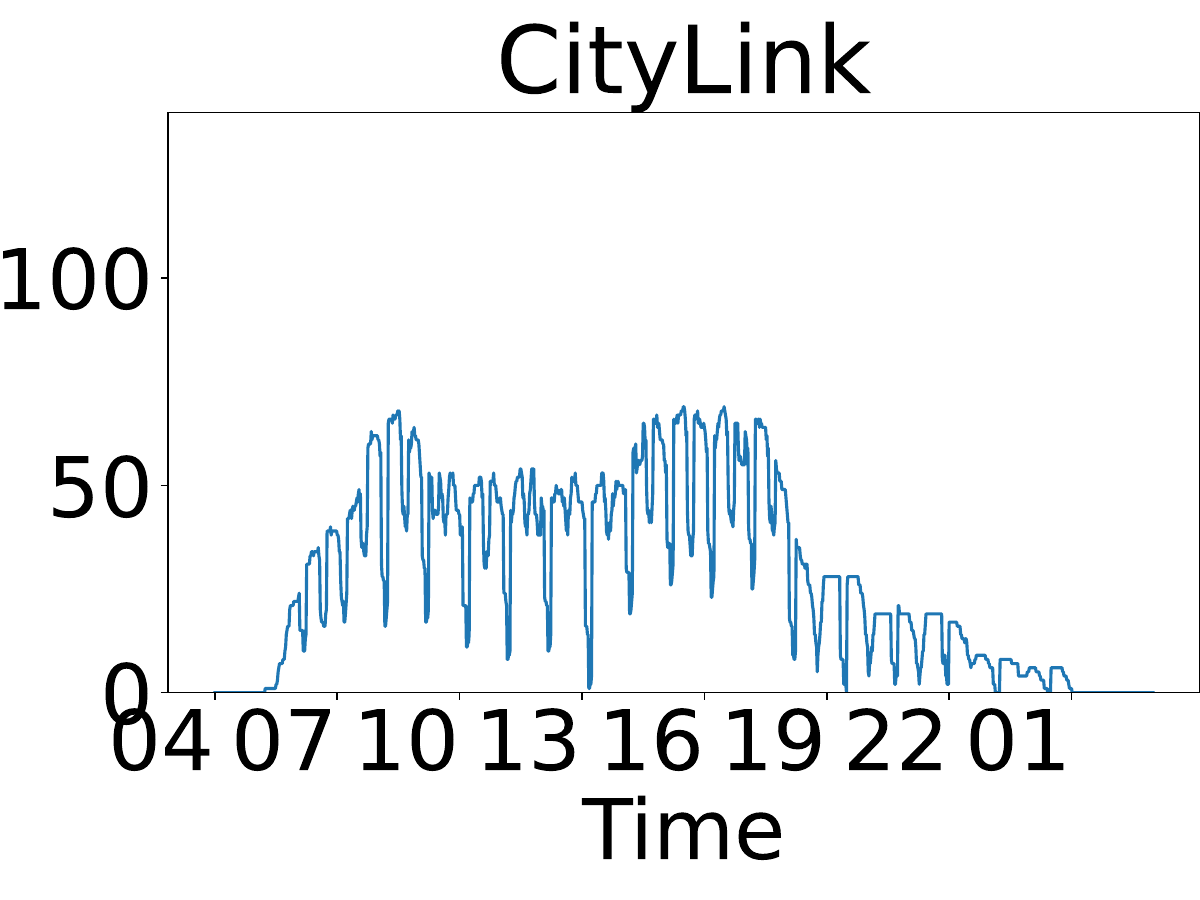}
  \end{subfigure}
    \hfill
  \begin{subfigure}{0.195\textwidth}
    \includegraphics[width=\linewidth]{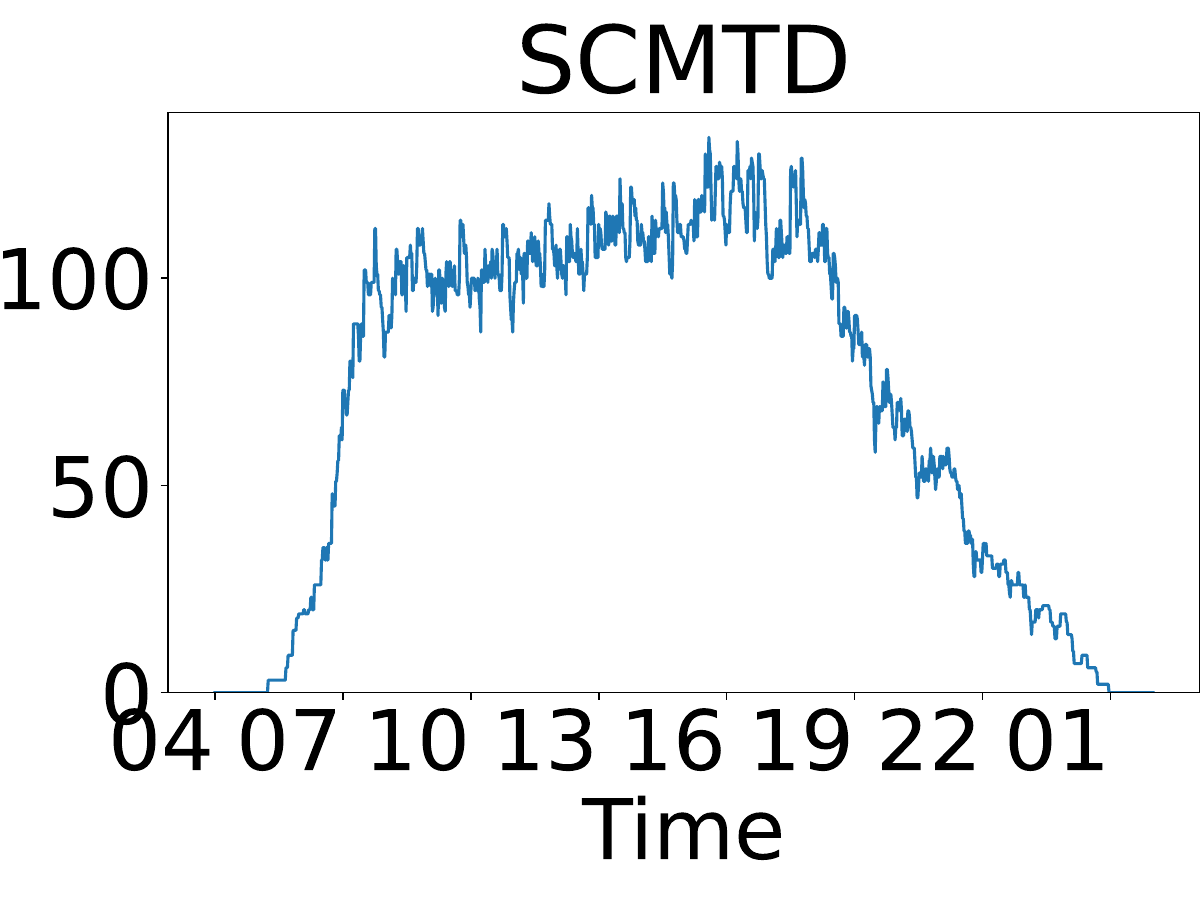}
  \end{subfigure}
      \hfill
  \begin{subfigure}{0.195\textwidth}
    \includegraphics[width=\linewidth]{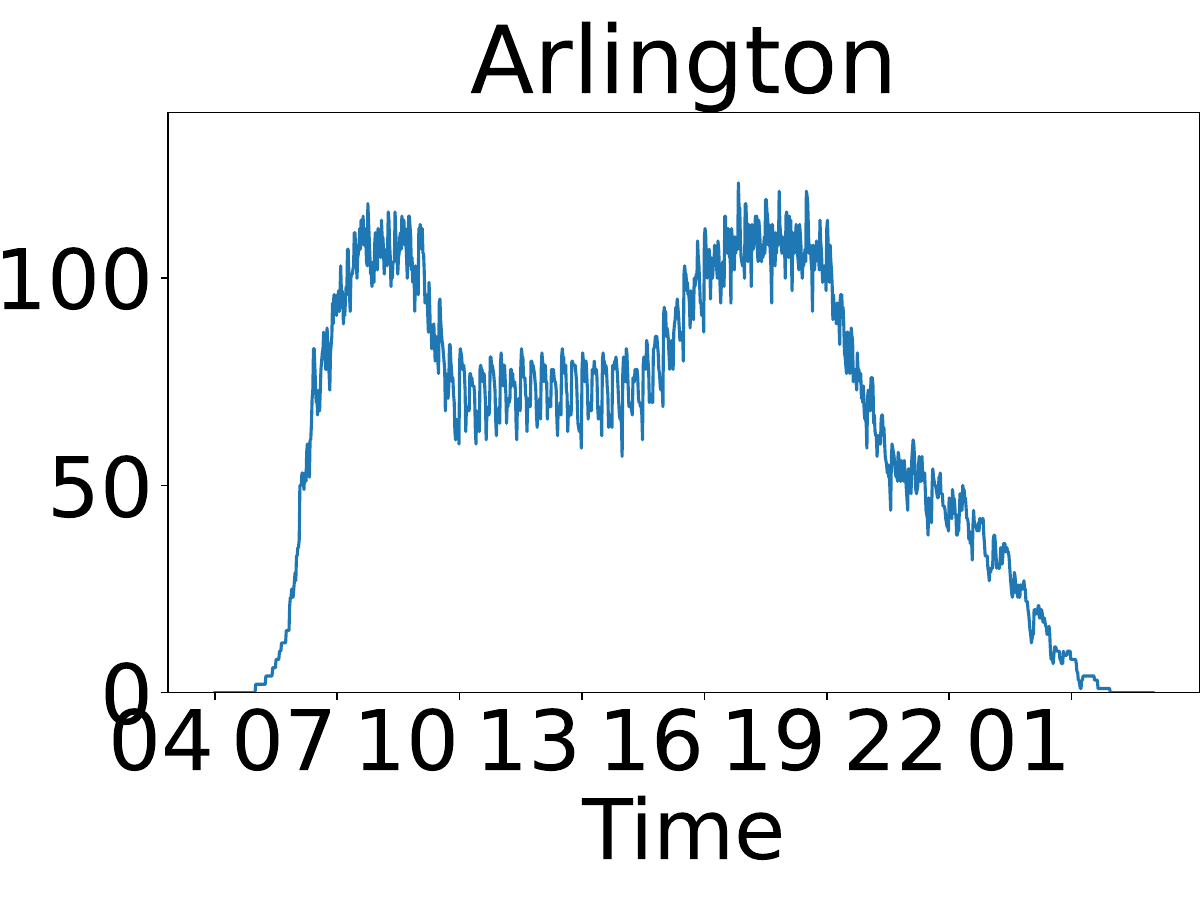}
  \end{subfigure}
  \caption{Distribution of simultaneous trips over time for sample networks}
  \label{fig:simultaneous_trips}
\end{figure}

To establish a lower limit on the required number of buses, we track the number of simultaneous trips across time for all the networks. For example, see Figure \ref{fig:simultaneous_trips}, where the SCMTD network requires a minimum of 134 buses to accommodate simultaneous afternoon trips. Table \ref{tab:datasets} also presents the lower limit on the number of buses for all networks. 

\begin{table}[htbp]
  \centering
  \caption{Summary of transit network characteristics}
    \begin{tabular}{p{40mm}p{12mm}p{12mm}p{12mm}p{12mm}p{12mm}p{18mm}}
    \hline
    \textbf{Network name} & \textbf{Trips} & \textbf{Routes} & \textbf{Stops} & \textbf{$|\candidateLocationSet|$} & \textbf{Depots used} & \textbf{Min buses required} \\
    \hline
Cornwall, Canada & 432 & 13 & 308 & 23 & 4 & 18\\
Milton, Canada & 501 & 9 & 312 & 5 & 3 & 22\\
Mountain Line, US & 716 & 9 & 181 & 11 & 5 & 42\\
Ann Arbor, US & 858 & 22 & 966 & 29 & 3 & 40\\
LA Go, US & 928 & 20 & 909 & 28 & 3 & 57\\
Cascades East, US & 1093 & 20 & 310 & 27 & 3 & 42\\
Denton County,US & 1148 & 32 & 560 & 25 & 5 & 32\\
Gold Coast, US & 1292 & 23 & 651 & 35 & 4 & 78\\
CityLink, US & 1364 & 20 & 1026 & 30 & 5 & 69\\
Intercity, US & 1591 & 22 & 954 & 22 & 4 & 85\\
SCMTD, US & 1830 & 24 & 786 & 26 & 4 & 134\\
Strathcona, Canada & 2355 & 28 & 373 & 30 & 3 & 117\\
Embark, US & 2830 & 34 & 1508 & 58 & 3 & 125\\
Arlington, US & 3337 & 19 & 635 & 32 & 5 & 123\\
    \hline 
    \end{tabular}%
  \label{tab:datasets}%
\end{table}%

Table \ref{tab:parameter_val} presents the parameter values used in this study. We apply an exchange rate factor of 1.09 to convert cost parameters from euros (\euro) to dollars (\$) based on some of the original sources mentioned in Table \ref{tab:parameter_val}. We assume dynamic electricity pricing with the following schedule: 444 \$/kWh from 9 AM to 2 PM, 555 \$/kWh from 2 PM to 4 PM and 9 PM to 9 AM, with a peak rate of 1355 \$/kWh from 4 PM to 9 PM \citep{electricschedule}. The table includes costs (in \$) and energy-related parameters sourced from existing literature. Some components of the objective are strategic, while other costs are from tactical and operational decisions. Therefore, these values are adjusted to account for a 12-year life cycle operational period for buses and charging infrastructure. We do not consider discounting and projected future energy prices, but these modifications can be easily incorporated with good-quality estimates. In the following sections, we detail the results of various heuristics for one network, Ann Arbor, US. Similar trends were observed in other cases.

\begin{table}[H]
 \small
  \centering
  \caption{Parameter values (operating costs are multiplied by $12 \times 365$ to account for a 12-year operational period)}
  \label{tab:parameter_val}
    \centering
    \begin{tabular}{lll}
    \hline
    \textbf{Parameters} & \textbf{Values} & \textbf{References}\\
    \hline
    $\BusCost$ & \$381,500 & \cite{dirks2022integration}\\
    $\LocationCost$ & \$218,000 & \cite{olsen2022location}\\
    $\text{per km travel cost}$ & \$2,100 & \cite{dirks2022integration} \\
    $\capacityPrice$ & 654 (\$/kW) & \cite{dirks2022integration}\\
    $\chargingRate$ & 2.505 (kWh/min) & \cite{jahic2019charging}\\
    $\BatteryCap$ & 300 (kWh) & \cite{stumpe2021study} \\
    $\MinEnergy$ & 45 (kWh) & \cite{sadati2019operational}\\
    $\maxTransfer$ & 2.505 (kWh/min) & \cite{jahic2019charging} \\
    \hline
    \end{tabular}%
  \label{tab:cost_parameters}
\end{table}%

\begin{figure}[H]
\centering
\begin{subfigure}{\textwidth}
\centering
\includegraphics[scale=0.145]{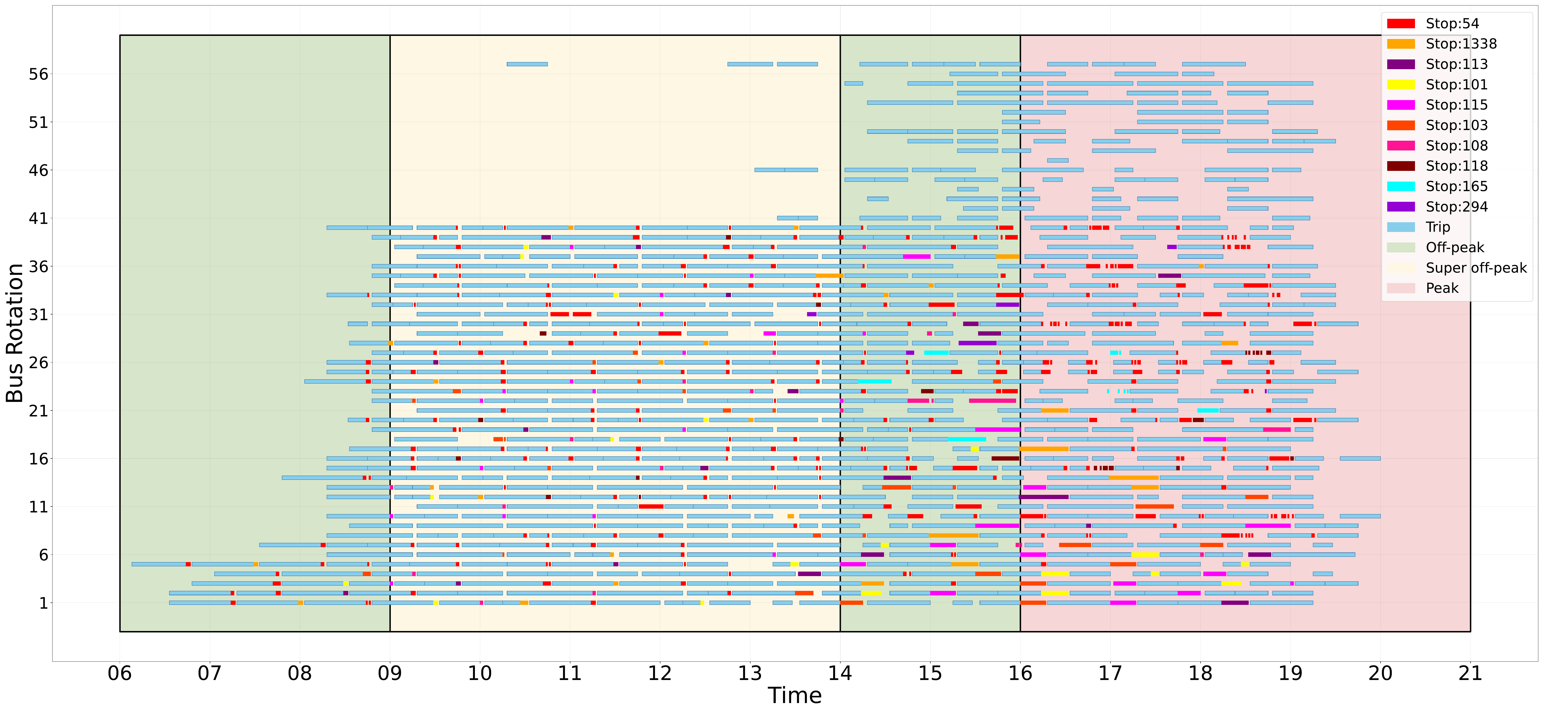}
\caption{Concurrent Scheduler model}
\label{fig:cs_gantt}
\end{subfigure}
\begin{subfigure}{\textwidth}
\centering
 \includegraphics[scale=0.145]{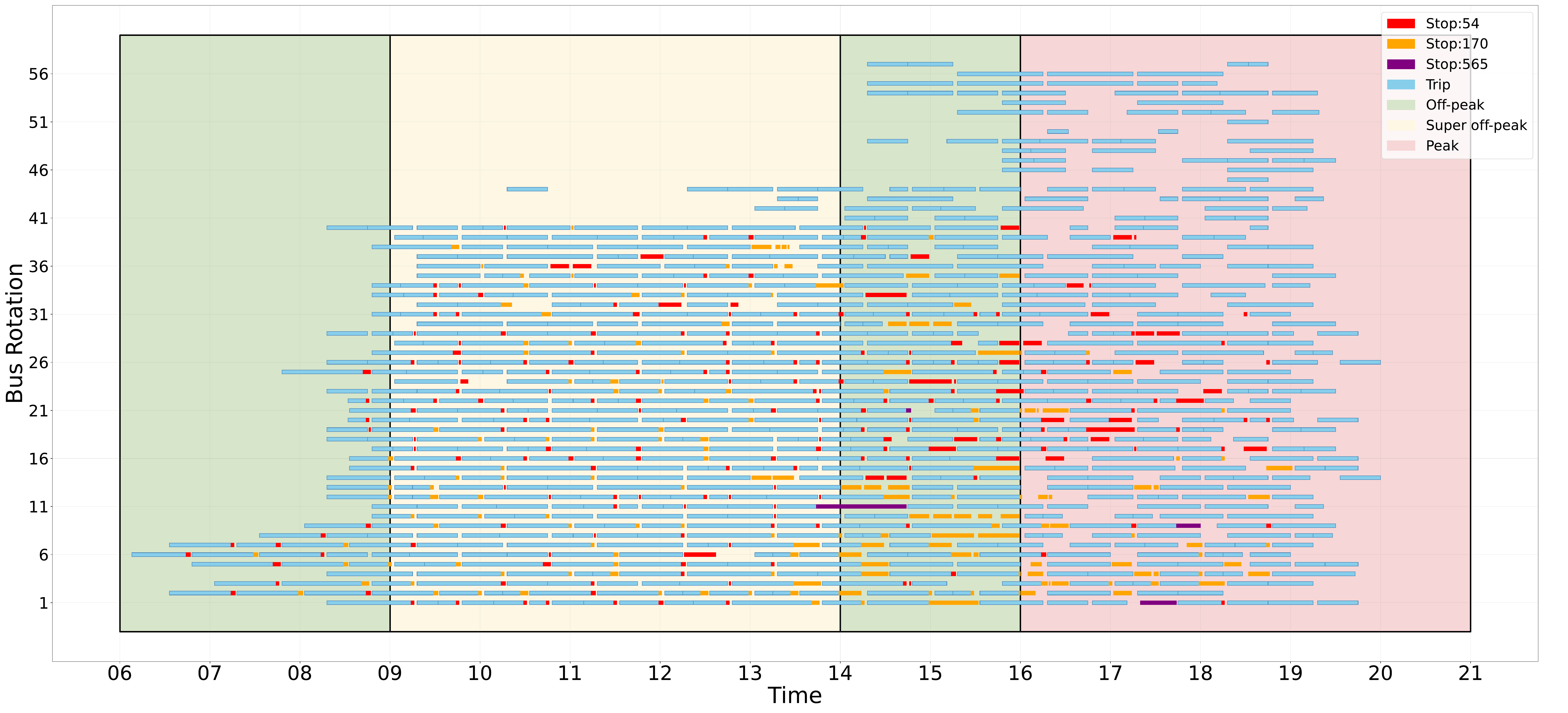}
\caption{Sequential model}
\label{fig:sequential_gantt}
\end{subfigure}
\begin{subfigure}{\textwidth}
\centering
  \includegraphics[scale=0.145]{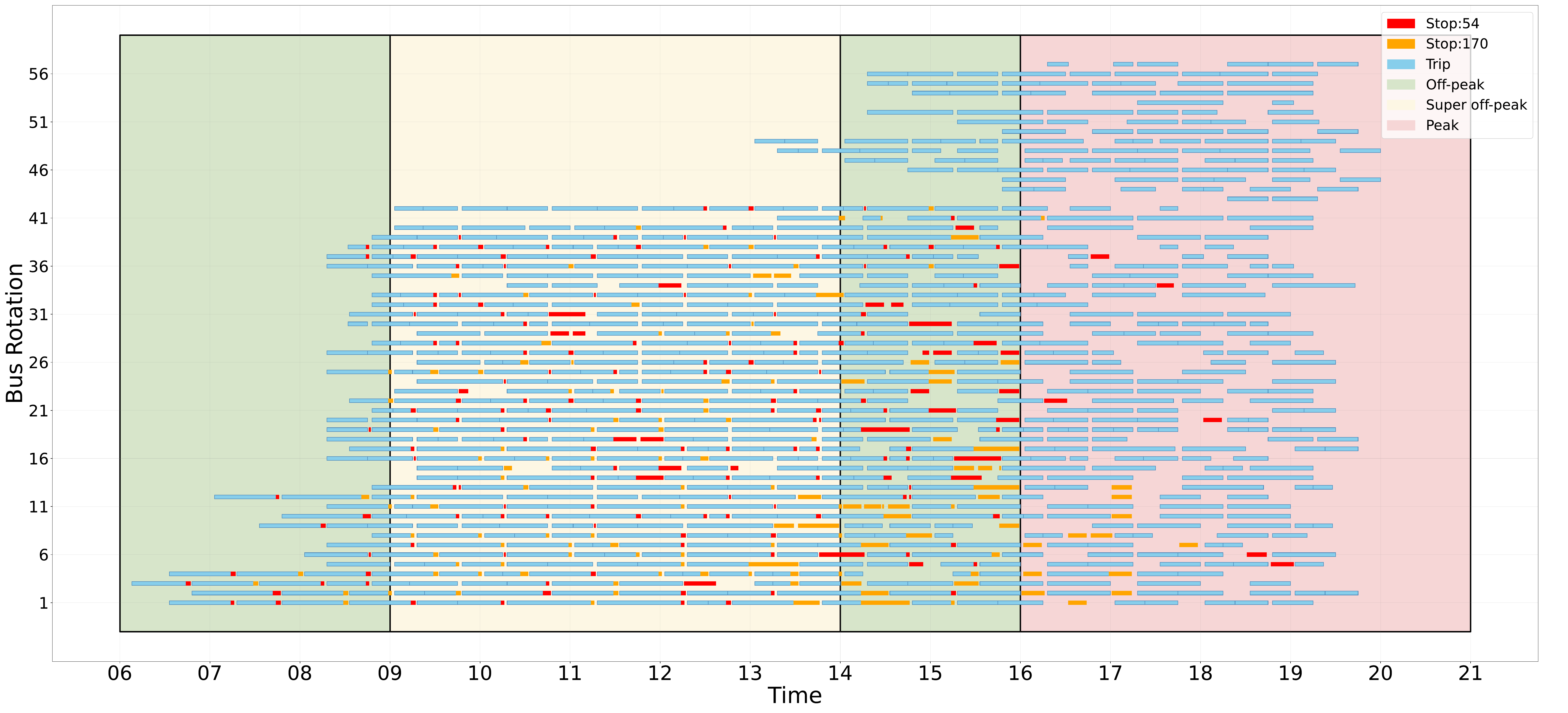}
\caption{Joint model}
\label{fig:joint_gantt}
\end{subfigure}
\caption{Bus schedules for the Ann Arbor network}
\end{figure}

\textbf{Initial Solution Generation:}
We generate an initial feasible solution for both the sequential and joint models using the CS algorithm. The initial solution often recommends opening a large number of charging locations. However, only a subset of these locations might be used when we solve the CSP using the MILP model. The fixed charging location cost is calculated based on the stations utilized by the post-CSP solution. Table \ref{tab:concurrent} displays the initial solutions for all networks, revealing that bus acquisition costs are a significant share of the total cost. The \% share of charge scheduling costs typically surpasses facility opening and deadheading costs. In Figure \ref{fig:cs_gantt}, the Gantt chart illustrates the bus schedules, featuring trips (indicated by blue bars) and charging events. Charging gaps result from split charging, while gaps in bus schedules include deadheading trips. The schedules are densely packed, aligning with the CS algorithm's logic. Colored backgrounds and borders help distinguish the energy price periods. The Ann Arbor network's initial solution employs 10 charging locations, with many charging activities scheduled during peak periods. Subsequent refinements of bus and charging schedules are made through ILS and CSP models.

\textbf{ILS -- Sequential Model Results:}
The ILS for the sequential model addresses the CLP-EVSP by refining the CS algorithm solution. Subsequently, the CEE model is employed to solve the CSP. Table \ref{tab:sequential} presents outcomes for various networks, revealing substantial cost savings over the CS solution. The reduction is primarily due to fewer charging locations compared to the CS solution, accompanied by notable deadheading cost reductions via exchange and shift operators. Charge scheduling costs are also consistently lower across transit networks. Figure \ref{fig:sequential_gantt} depicts the resulting Gantt chart of bus schedules, featuring structural changes due to trip exchanges and shifts. 

\begin{table}[H]
  \centering
  \caption{Concurrent Scheduler results}
    \begin{tabular}{p{24mm}p{10mm}p{16mm}p{20mm}p{18mm}p{18mm}p{18mm}p{18mm}}
    \toprule
    \multirow{2}[4]{24mm}{\parbox{24mm}{\textbf{Network\\ Name}}} & \multirow{2}[4]{10mm}{\textbf{Buses\\ Used}} & \multirow{2}[4]{16mm}{\parbox{16mm}{\textbf{Charging\\ Locations\\ Opened}}} & \multirow{2}[4]{20mm}{\parbox{20mm}{\textbf{Total Cost\\ (Million\\ \$)}}} & \multicolumn{4}{c}{\textbf{\% Share}} \\
\cmidrule{5-8}          &       &       &       & \textbf{Bus acquisition} & \textbf{Facility Opening} & \textbf{Deadhead} & \textbf{CSP} \\
    \bottomrule
    Cornwall & 20 & 1 & 10.20 & 74.81 & 2.14 & 2.36 & 20.69 \\
Milton & 25 & 4 & 14.31 & 66.64 & 6.09 & 2.39 & 24.88 \\
Mountain Line & 45 & 1 & 22.36 & 76.76 & 0.97 & 1.99 & 20.28 \\
Ann Arbor & 57 & 10 & 31.93 & 68.11 & 6.83 & 7.30 & 17.76 \\
LA Go & 87 & 13 & 57.67 & 57.55 & 4.91 & 25.69 & 11.85 \\
Cascades East & 71 & 11 & 41.67 & 65.00 & 5.75 & 17.13 & 12.12 \\
Denton County & 46 & 13 & 32.23 & 54.46 & 8.79 & 14.96 & 21.79 \\
Gold Coast & 97 & 12 & 57.52 & 64.33 & 4.55 & 14.34 & 16.78 \\
CityLink & 87 & 13 & 56.69 & 58.55 & 5.00 & 15.16 & 21.29 \\
Intercity & 195 & 11 & 95.70 & 77.73 & 2.51 & 6.99 & 12.77 \\
SCMTD & 187 & 11 & 134.62 & 52.99 & 1.78 & 19.60 & 25.63 \\
Strathcona & 146 & 4 & 82.98 & 67.13 & 1.05 & 7.86 & 23.96 \\
Embark & 183 & 31 & 116.84 & 59.75 & 5.78 & 15.07 & 19.40 \\
Arlington & 145 & 20 & 84.43 & 65.52 & 5.16 & 11.22 & 18.10 \\
    \hline 
    \end{tabular}%
  \label{tab:concurrent}%
\end{table}%

\begin{figure}[H]
  \centering
    \includegraphics[scale=0.3]{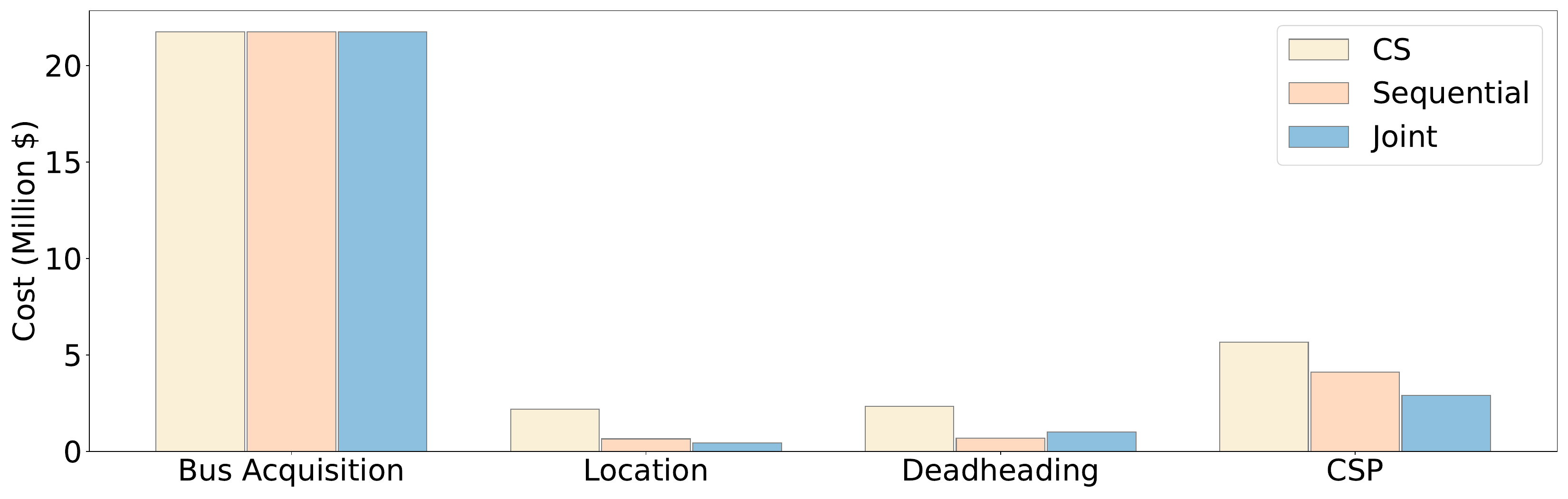}
    \caption{Cost components of the objective for the Ann Arbor network}
  \label{fig:cost_components}
\end{figure}

\textbf{ILS -- Joint Model Results:}
We then applied the ILS technique to the joint CLP-EVSP-CSP model, using the CS algorithm solution as a starting point. This approach yields noticeable savings in charging costs, even though deadheading costs increase compared to the sequential model for most networks. The joint model's results, summarized in Table \ref{tab:joint}, show overall cost savings across all networks. The values in the table indicate savings relative to the solutions of the sequential model. 

Figure \ref{fig:cost_components} breaks down the costs associated with using the CS, sequential, and joint models for the Ann Arbor network. The joint model utilizes fewer charging locations than the sequential model. Despite the increased deadheading costs in the joint model, they remain lower than those incurred using the CS algorithm. Figure \ref{fig:joint_gantt} displays schedules from the joint model, where one can prominently notice significant cost savings in electricity consumption during peak periods, as only a few buses are charged. Similar trends were observed in other networks.

A majority of the cost is incurred in bus acquisition. For agencies where bus procurement happens through subsidies or when gasoline/diesel buses are exchanged/replaced, these savings estimates are conservative and can actually be closer to the CSP savings.

\begin{table}[H]
  \centering
  \caption{ILS model results (Instances marked with a `*' were solved using a hybrid method where the top $400$ shifts/exchanges were considered based on savings in deadheading costs.)}
  \begin{subtable}{\linewidth}
  \caption{Sequential model}
    \begin{tabular}{p{24mm}p{10mm}p{16mm}p{15mm}p{15mm}p{15mm}p{15mm}p{12mm}p{13mm}}
    \toprule
    \multirow{2}[4]{24mm}{\parbox{24mm}{\textbf{Network\\ Name}}} & \multirow{2}[4]{10mm}{\textbf{Buses\\ Used}} & \multirow{2}[4]{16mm}{\parbox{16mm}{\textbf{Charging\\ Locations\\ Opened}}} & \multicolumn{5}{c}{\textbf{\% Savings w.r.t. CS Model}} \\
\cmidrule{4-8}          &       &       &  \textbf{Overall} & \textbf{Bus acquisition} & \textbf{Facility Opening} & \textbf{Deadhead} & \textbf{CSP} & \textbf{Runtime (s)} \\
    \bottomrule
    Cornwall & 20 & 1 & $0.51$ & $0.00$ & $0.00$ & $15.93$ & $0.67$ & $0.08$\\
Milton & 25 & 2 & $7.48$ & $0.00$ & $50.00$ & $67.26$ & $11.38$ & $0.10$\\
Mountain Line & 45 & 1 & $1.76$ & $0.00$ & $0.00$ & $48.89$ & $3.90$ & $0.14$\\
Ann Arbor & 57 & 3 & $14.84$ & $0.00$ & $70.00$ & $70.49$ & $27.69$ & $0.28$\\
LA Go & 87 & 11 & $8.13$ & $0.00$ & $15.38$ & $20.42$ & $18.00$ & $0.45$\\
Cascades East & 71 & 8 & $11.24$ & $0.00$ & $27.27$ & $46.23$ & $14.46$ & $0.46$\\
Denton County & 46 & 7 & $17.60$ & $0.00$ & $46.15$ & $60.26$ & $20.78$ & $0.54$\\
Gold Coast & 97 & 7 & $13.35$ & $0.00$ & $41.67$ & $45.16$ & $29.67$ & $0.73$\\
CityLink & 87 & 6 & $17.05$ & $0.00$ & $53.85$ & $67.85$ & $19.11$ & $0.82$\\
Intercity & 195 & 7 & $9.99$ & $0.00$ & $36.36$ & $66.98$ & $34.47$ & $1.28$\\
SCMTD & 187 & 9 & $14.77$ & $0.00$ & $18.18$ & $45.66$ & $21.44$ & $1.76$\\
Strathcona & 146 & 3 & $6.28$ & $0.00$ & $25.00$ & $52.96$ & $7.73$ & $1.79$\\
Embark & 183 & 28 & $15.16$ & $0.00$ & $9.68$ & $57.39$ & $30.69$ & $4.07$\\
Arlington & 145 & 12 & $11.49$ & $0.00$ & $40.00$ & $52.19$ & $19.73$ & $6.25$\\
    \hline 
    \end{tabular}%
  \label{tab:sequential}%
  \end{subtable}
  
  \vspace{1em}  

  \begin{subtable}{\linewidth}
  \caption{Joint model} 
    \begin{tabular}{p{24mm}p{10mm}p{16mm}p{15mm}p{15mm}p{15mm}p{15mm}p{12mm}p{13mm}}
    \toprule
    \multirow{2}[4]{24mm}{\parbox{24mm}{\textbf{Network\\ Name}}} & \multirow{2}[4]{10mm}{\textbf{Buses\\ Used}} & \multirow{2}[4]{16mm}{\parbox{16mm}{\textbf{Charging\\ Locations\\ Opened}}} & \multicolumn{5}{c}{\textbf{\% Savings w.r.t. Sequential Model}} \\
\cmidrule{4-8}          &       &       &      \textbf{Overall} & \textbf{Bus acquisition} & \textbf{Facility Opening} & \textbf{Deadhead} & \textbf{CSP} & \textbf{Runtime (s)} \\
    \bottomrule
    Cornwall & 20 & 1 & $7.69$ & $0.00$ & $0.00$ & $-3.88$ & $37.58$ & $41.58$\\
Milton & 25 & 2 & $6.52$ & $0.00$ & $0.00$ & $-46.93$ & $29.01$ & $59.17$\\
Mountain Line & 45 & 1 & $5.11$ & $0.00$ & $0.00$ & $10.03$ & $25.25$ & $90.65$\\
Ann Arbor & 57 & 2 & $4.04$ & $0.00$ & $33.33$ & $-47.11$ & $29.34$ & $102.70$\\
LA Go & 89 & 8 & $0.99$ & $-2.30$ & $27.27$ & $-4.09$ & $19.90$ & $865.11$\\
Cascades East & 71 & 5 & $5.17$ & $0.00$ & $37.50$ & $-7.47$ & $35.74$ & $702.00$\\
Denton County & 46 & 4 & $4.84$ & $0.00$ & $42.86$ & $-9.71$ & $14.67$ & $509.89$\\
Gold Coast & 98 & 6 & $0.94$ & $-1.03$ & $14.29$ & $-10.68$ & $16.45$ & $1420.18$\\
CityLink & 88 & 2 & $6.66$ & $-1.15$ & $66.67$ & $-5.43$ & $28.57$ & $4631.02$\\
Intercity & 195 & 2 & $7.70$ & $0.00$ & $71.43$ & $-16.36$ & $73.65$ & $11897.79$\\
SCMTD* & 193 & 7 & $0.13$ & $-3.21$ & $22.22$ & $-11.66$ & $13.54$ & $7024.60$\\
Strathcona* & 147 & 1 & $3.70$ & $-0.68$ & $66.67$ & $-27.41$ & $19.97$ & $1265.62$\\
Embark* & 188 & 19 & $2.95$ & $-2.73$ & $32.14$ & $-13.54$ & $24.74$ & $3135.22$\\
Arlington* & 146 & 7 & $1.55$ & $-0.69$ & $41.67$ & $-29.47$ & $14.56$ & $142177.98$\\
    \hline 
    \end{tabular}%
  \label{tab:joint}%
  \end{subtable}
\end{table}%

\subsection{Extended Analysis}

\textbf{Power Requirements:} The predominant objective component in the CSP is electricity consumption costs. Table \ref{tab:power} summarizes the maximum and average power requirements for charging locations opened using sequential and joint models. In most instances, the maximum power required by the joint model is lower than that of the sequential model. In cases where the sequential model exhibits lower average power requirements, the difference is relatively modest. Given that we model these decisions as continuous variables, the overall power requirement can also be translated into the equivalent number of chargers to be stationed at a charging facility. For example, a power requirement of 300 kW will warrant installing two 150 kW chargers at the charging location. 

\begin{table}[H]
  \centering
  \caption{Power requirements of charging stations in kW (The lower value of the two models are shown in bold)} 
  \begin{tabular}{p{28mm}p{30mm}p{30mm}p{30mm}p{30mm}}
    \hline
    \textbf{Network Name} & \multicolumn{2}{c}{\textbf{Maximum Power Required}} & \multicolumn{2}{c}{\textbf{Average Power Required}} \\
    \cmidrule(lr){2-3} \cmidrule(lr){4-5}
    & \textbf{Sequential} & \textbf{Joint} & \textbf{Sequential} & \textbf{Joint} \\
    \hline
    Cornwall & 811.22 & \textcolor{black}{\textbf{262.29}} & 811.22 & \textcolor{black}{\textbf{262.29}} \\
Milton & 543.13 & \textcolor{black}{\textbf{300.60}} & 421.86 & \textcolor{black}{\textbf{250.31}} \\
Mountain Line & 1032.66 & \textcolor{black}{\textbf{681.01}} & 1032.66 & \textcolor{black}{\textbf{681.01}} \\
Ann Arbor & 694.65 & \textcolor{black}{\textbf{601.65}} & \textcolor{black}{\textbf{450.49}} & 601.43 \\
LA Go & 377.14 & \textcolor{black}{\textbf{300.60}} & \textcolor{black}{\textbf{148.56}} & 152.02 \\
Cascades East & 314.09 & \textcolor{black}{\textbf{219.98}} & 179.34 & \textcolor{black}{\textbf{144.06}} \\
Denton County & \textcolor{black}{\textbf{302.00}} & 315.03 & \textcolor{black}{\textbf{198.84}} & 274.50 \\
Gold Coast & \textcolor{black}{\textbf{536.57}} & 549.97 & 246.42 & \textcolor{black}{\textbf{241.56}} \\
CityLink & 2004.53 & \textcolor{black}{\textbf{1106.46}} & \textcolor{black}{\textbf{549.82}} & 778.68 \\
Intercity & 1503.00 & \textcolor{black}{\textbf{498.72}} & 318.35 & \textcolor{black}{\textbf{278.07}} \\
SCMTD & \textcolor{black}{\textbf{1653.30}} & 1675.09 & \textcolor{black}{\textbf{600.01}} & 730.78 \\
Strathcona & 4401.20 & \textcolor{black}{\textbf{3007.14}} & \textcolor{black}{\textbf{2068.27}} & 3007.14 \\
Embark & \textcolor{black}{\textbf{450.90}} & \textcolor{black}{\textbf{450.90}} & \textcolor{black}{\textbf{183.29}} & 191.17 \\
Arlington & \textcolor{black}{\textbf{300.60}} & 332.96 & \textcolor{black}{\textbf{237.97}} & 298.11 \\
    \hline
  \end{tabular}%
  \label{tab:power}%
\end{table}%

\begin{figure}[H]
\centering
  \includegraphics[width=0.95\textwidth]{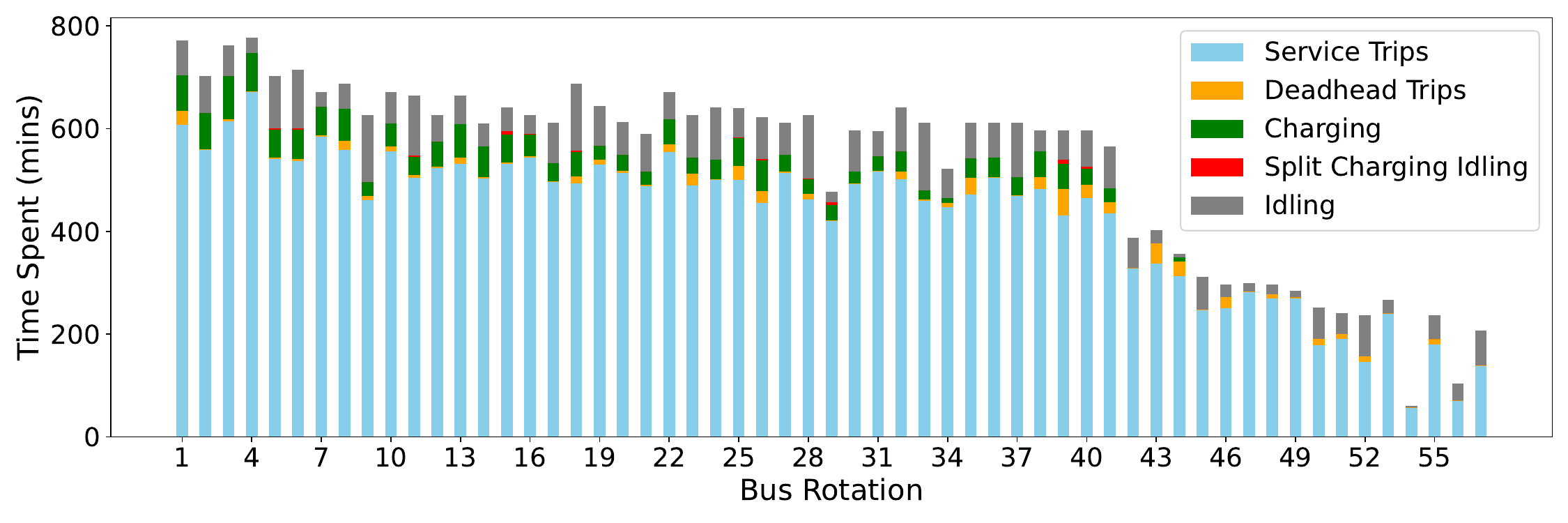}
\caption{Activity-based duration for buses in the Ann Arbor network using the joint model}
\label{fig:time_spent_2}
\end{figure}

\textbf{Activity-Based Time Allocation:} Figure \ref{fig:time_spent_2} illustrates the time allocated to various activities, as determined by the joint model for the Ann Arbor network. Buses engage in service trips, deadhead trips, recharging activities, idling during recharging due to split charging, and true idling when not plugged into the grid. On average, buses spend approximately 8 hours on trips and charging activities. Some buses covering short distances may not require opportunity charging but require overnight charging. Periodic bus swapping between different rotations/itineraries might be necessary to balance overall usage and mitigate potential long-term maintenance issues.

\begin{figure}[h]
  \centering
  \begin{subfigure}{0.32\textwidth}
    \includegraphics[width=\linewidth]{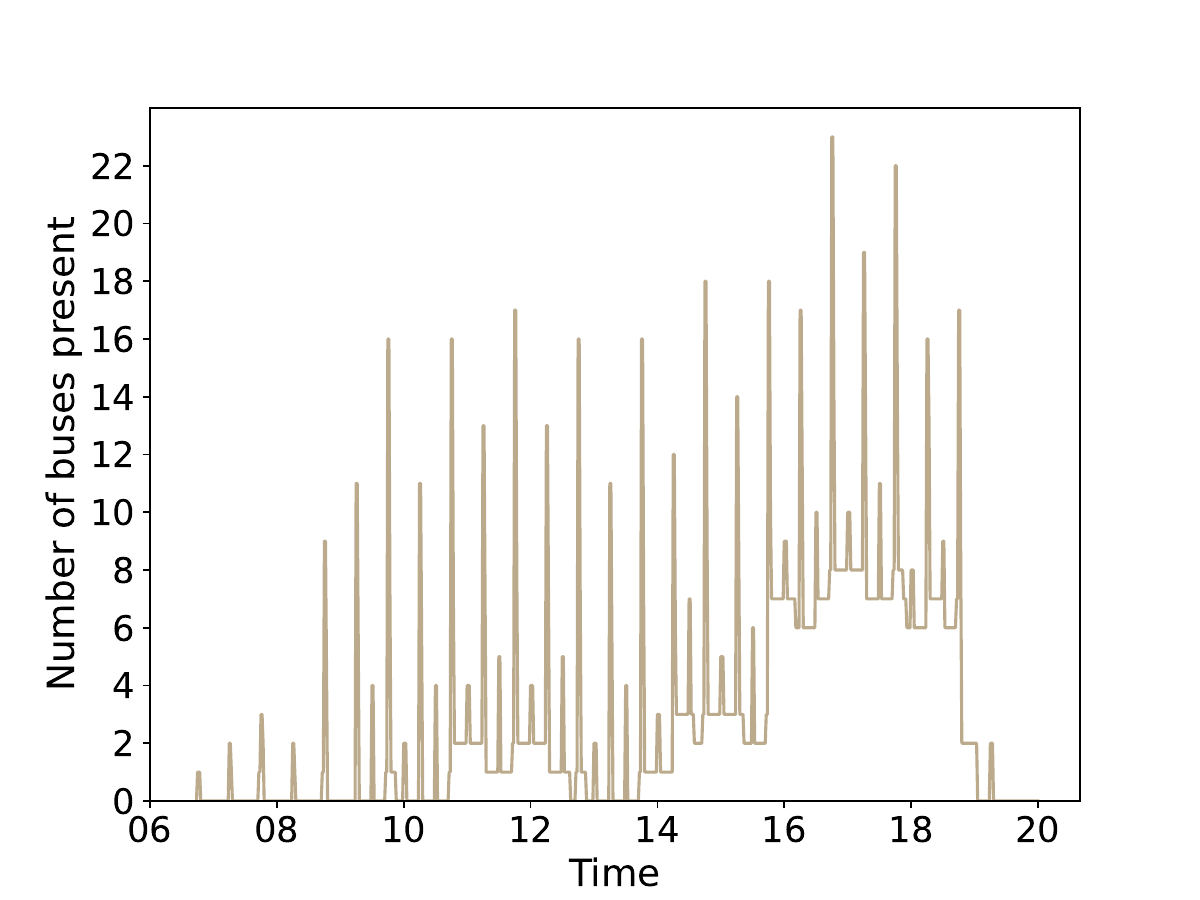}
    \caption{Buses present (Sequential)}
    \label{fig:buses_present_sequential}
  \end{subfigure}
  \hfill
  \begin{subfigure}{0.32\textwidth}
    \includegraphics[width=\linewidth]{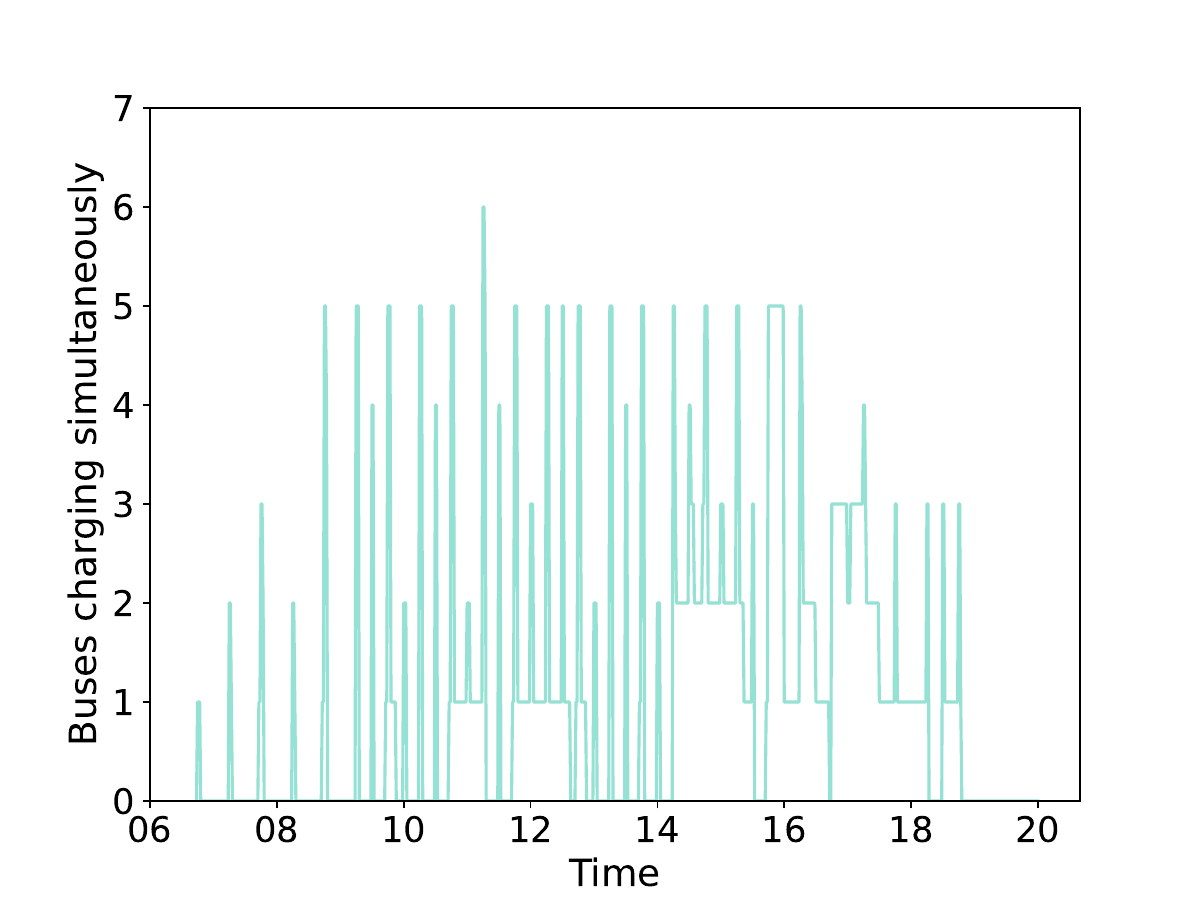}
    \caption{Concurrent charging (Sequential)}
    \label{fig:simultaneous_sequential}
  \end{subfigure}
  \hfill
  \begin{subfigure}{0.32\textwidth}
    \includegraphics[width=\linewidth]{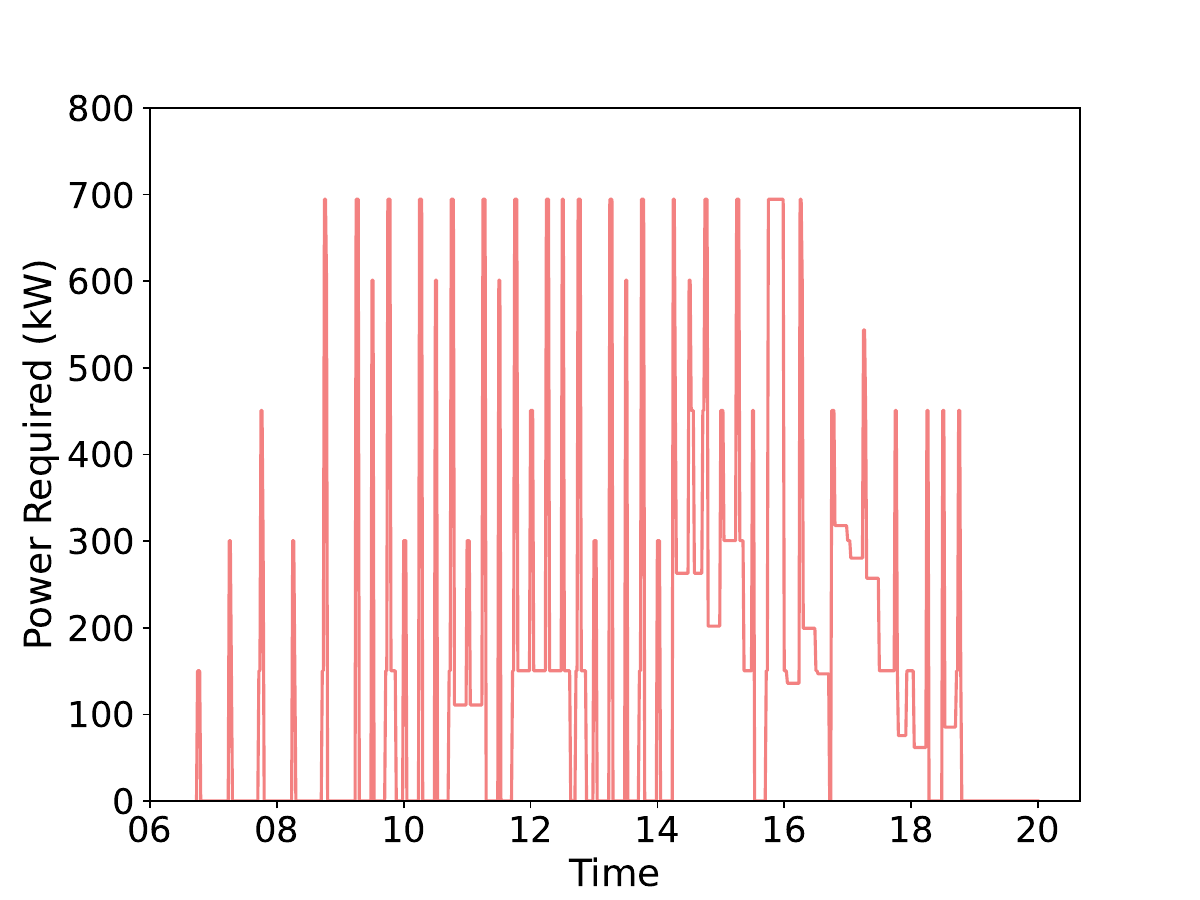}
    \caption{Power requirement (Sequential)}
    \label{fig:power_sequential}
  \end{subfigure}

  \begin{subfigure}{0.32\textwidth}
    \includegraphics[width=\linewidth]{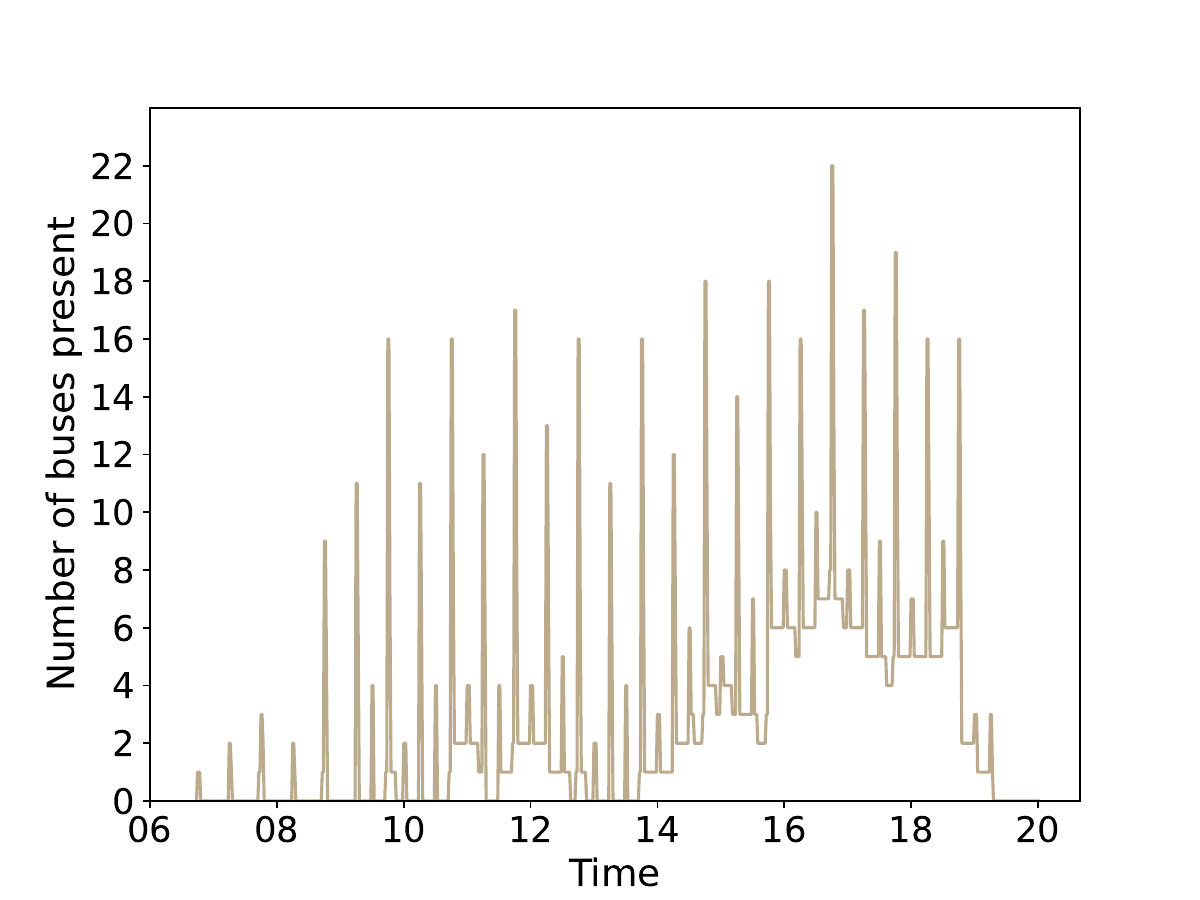}
    \caption{Buses present (Joint)}
    \label{fig:buses_present_joint}
  \end{subfigure}
  \hfill
  \begin{subfigure}{0.32\textwidth}
    \includegraphics[width=\linewidth]{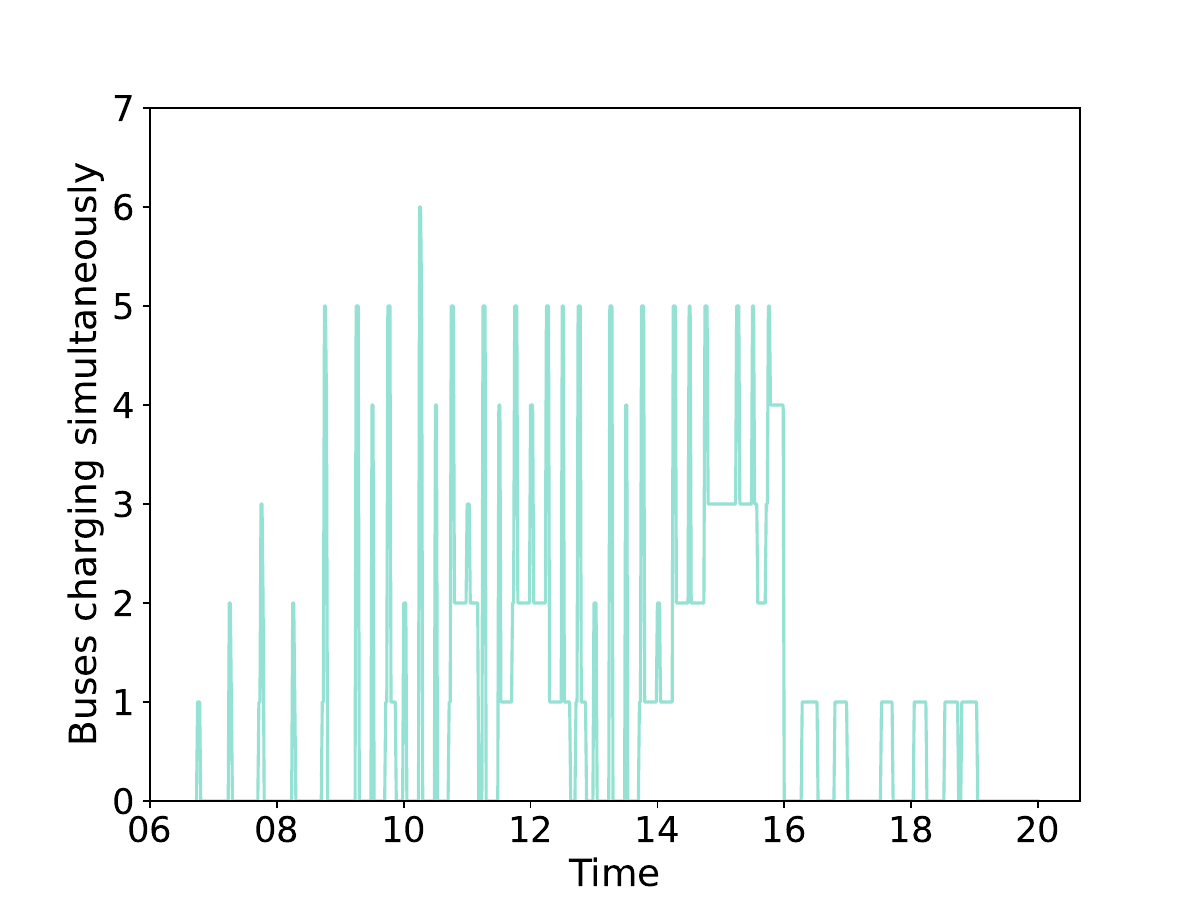}
    \caption{Concurrent charging (Joint)}
    \label{fig:simultaneous_joint}
  \end{subfigure}
  \hfill
  \begin{subfigure}{0.32\textwidth}
    \includegraphics[width=\linewidth]{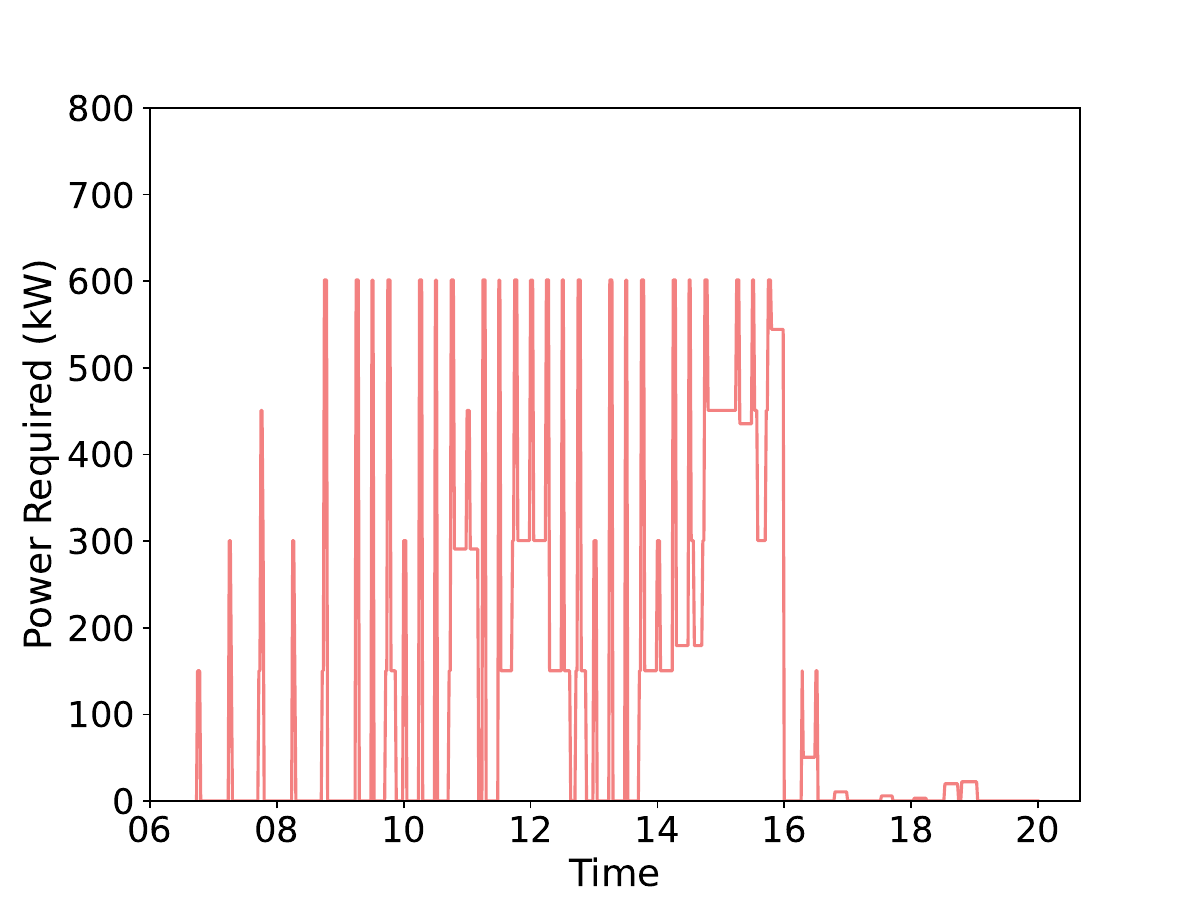}
    \caption{Power requirement (Joint)}
    \label{fig:power_joint}
  \end{subfigure}

  \caption{Station-level comparison between the sequential and joint models for the charging location at the Blake Transit Center of the Ann Arbor network}
  \label{fig:combined_figure_location}
\end{figure}

\textbf{Station-level Analysis:}
Figure \ref{fig:combined_figure_location} displays station-level metrics for a sample charging station situated at the Blake Transit Center (bus stop 54), Ann Arbor. The metrics plotted include the number of buses present at the station, the number of buses charging simultaneously, and the power demand at the charging station over time. 
Under sequential and joint models, a maximum of 23 and 22 buses are observed to be simultaneously present at the charging location, respectively. There are slight variations in the number of vehicles present (see Figures \ref{fig:buses_present_sequential} \& \ref{fig:buses_present_joint}) between the sequential and joint models due to variations in vehicle scheduling. Both sequential (Figure \ref{fig:simultaneous_sequential}) and joint models (Figure \ref{fig:simultaneous_joint}) accommodate a maximum of six simultaneously charging buses. However, 
the overall power demand at the charging station remains lower for the joint model than that observed in the sequential model (see Figures \ref{fig:power_sequential} \& \ref{fig:power_joint}) at any given point in time. We notice this pattern despite a nearly equal number of buses being present at the station over time in both the sequential and joint models. Furthermore, very few charging events take place during the peak periods in the joint case. 

\begin{figure}[h]
  \centering
  \begin{subfigure}{0.32\textwidth}
    \centering
    \includegraphics[width=0.95\linewidth]{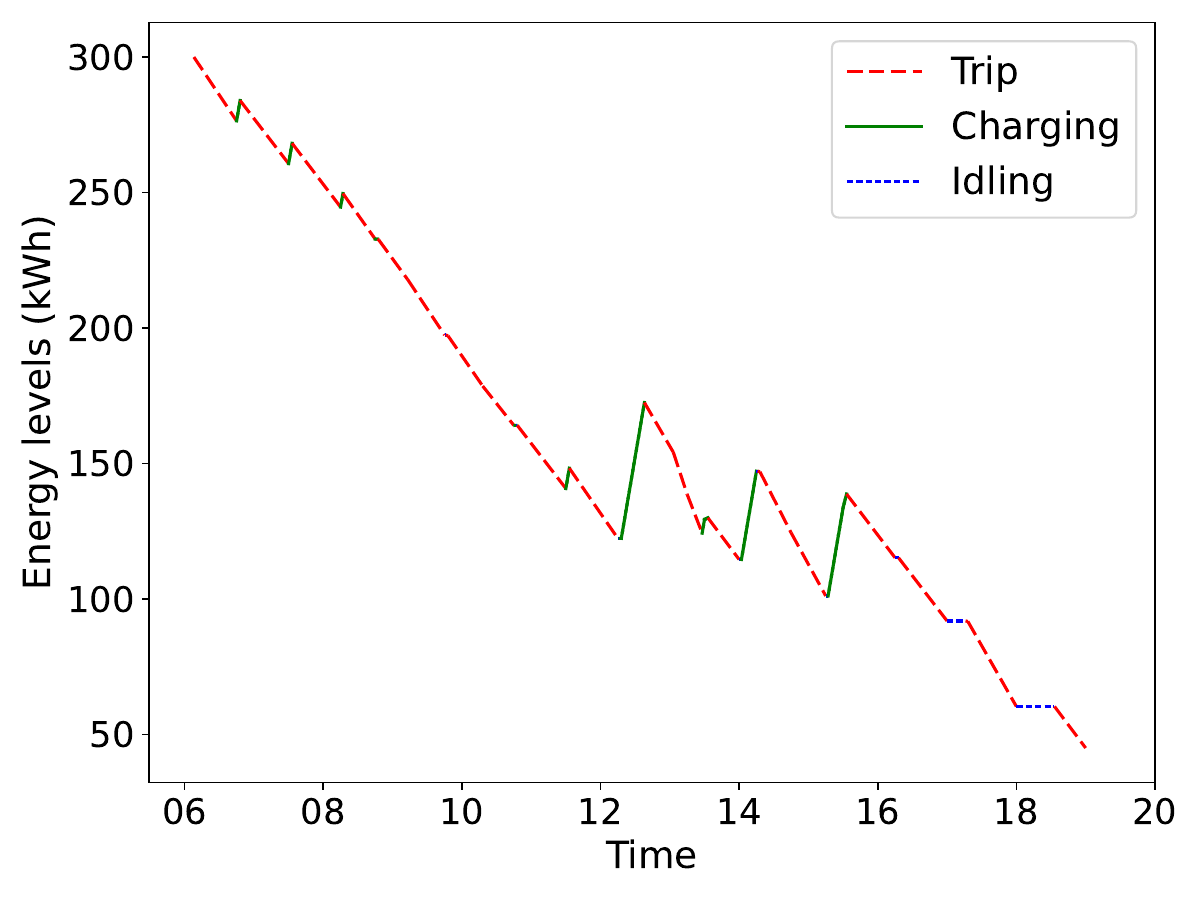}
    \caption{Bus 1}
  \end{subfigure}
  \hfill
  \begin{subfigure}{0.32\textwidth}
    \centering
    \includegraphics[width=0.95\linewidth]{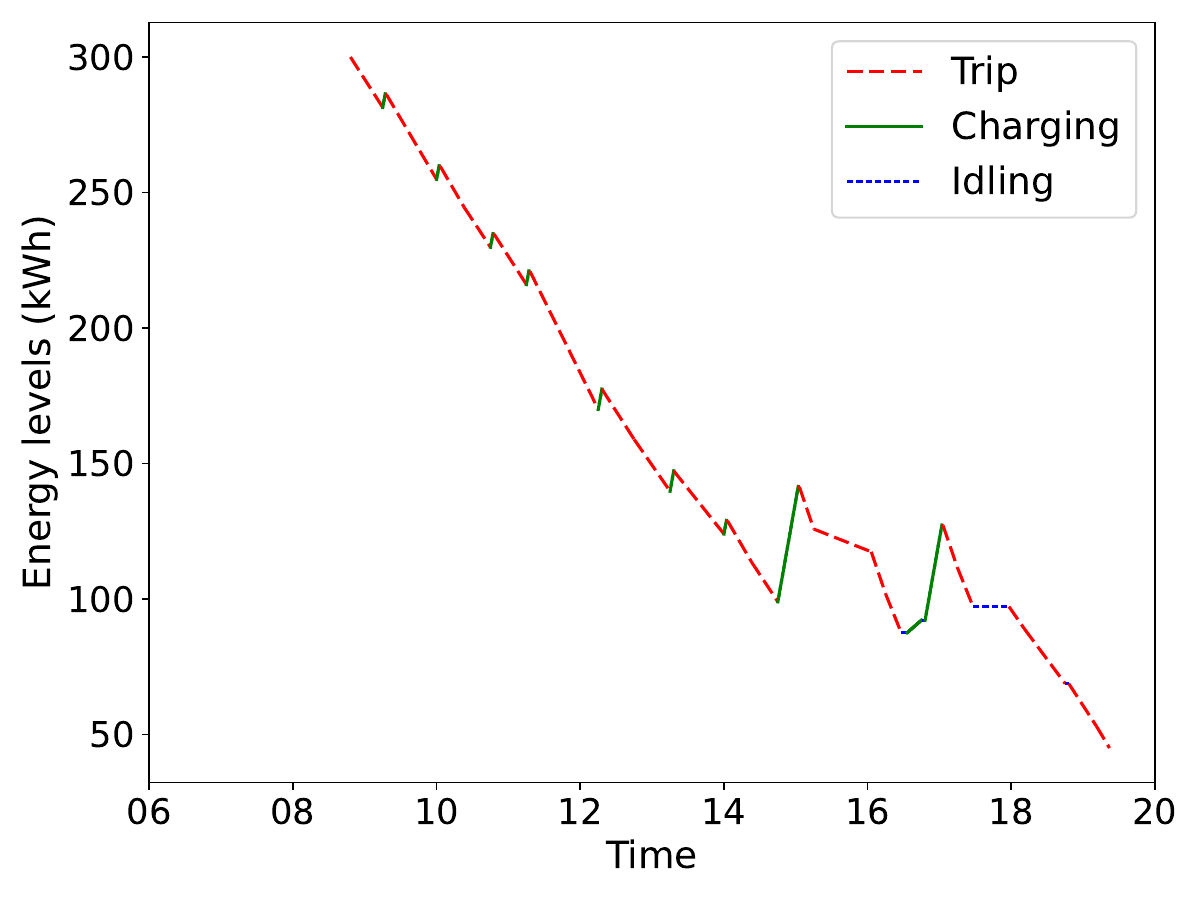}
    \caption{Bus 26}
  \end{subfigure}
  \hfill
  \begin{subfigure}{0.32\textwidth}
    \centering
    \includegraphics[width=0.95\linewidth]{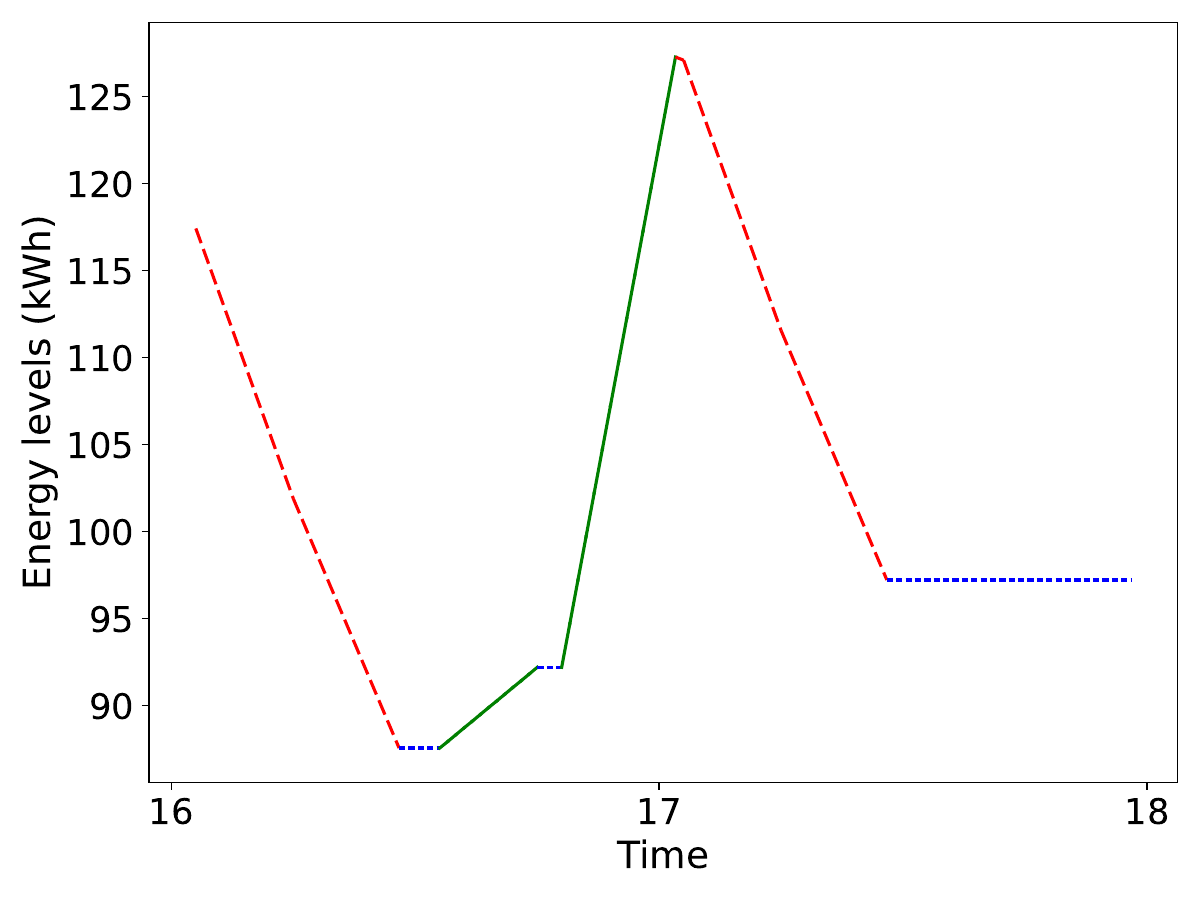}
    \caption{Bus 26 (16:00 -- 18:00)}
    \label{fig:bus_26_zoomed}
  \end{subfigure}
  \caption{Bus energy levels in the Ann Arbor network using the joint model}
  \label{fig:bus_energy_level}
\end{figure}

\textbf{Bus-level Analysis:}
We can also monitor the energy levels of each bus during operational hours using the results obtained from the ILS methods. Figure \ref{fig:bus_energy_level} traces the energy levels of buses 1 and 26 using the joint model solution. Buses commence daily operations with a fully charged battery capacity of 300 kWh. Buses that require charging during the day finish daily operations with an energy level of 45 kWh, indicating no overcharging. The red lines denote energy dissipation during service or deadheading trips. Green lines represent charging activities, while blue lines indicate bus idling. Varying slopes of the red lines result from different speeds during service and deadheading trips. Bus speeds fluctuate based on the schedule and the assumed deadheading speed. On the other hand, varying slopes of the green lines result from different charging rates and the split charging assumption. Bus $1$ exhibits no split charging, whereas bus $26$ engages in split charging activities from 16:45 to 16:48, shown in Figure \ref{fig:bus_26_zoomed}. 

\begin{figure}[H]
  \centering
  \begin{subfigure}{0.49\textwidth}
  \centering
    \includegraphics[width=0.8\linewidth]{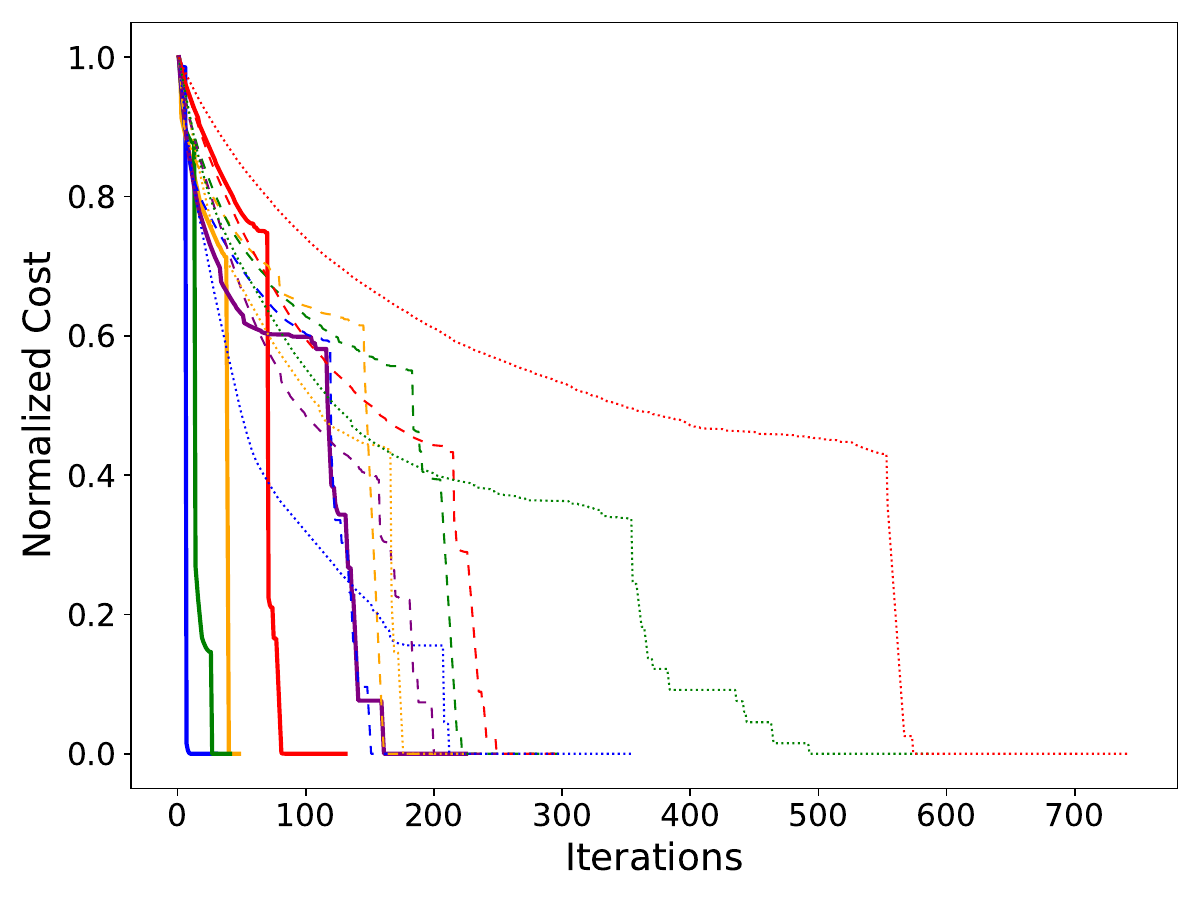}
    \caption{Sequential model}
  \end{subfigure}
  \hfill
  \begin{subfigure}{0.49\textwidth}
  \centering
    \includegraphics[width=0.8\linewidth]{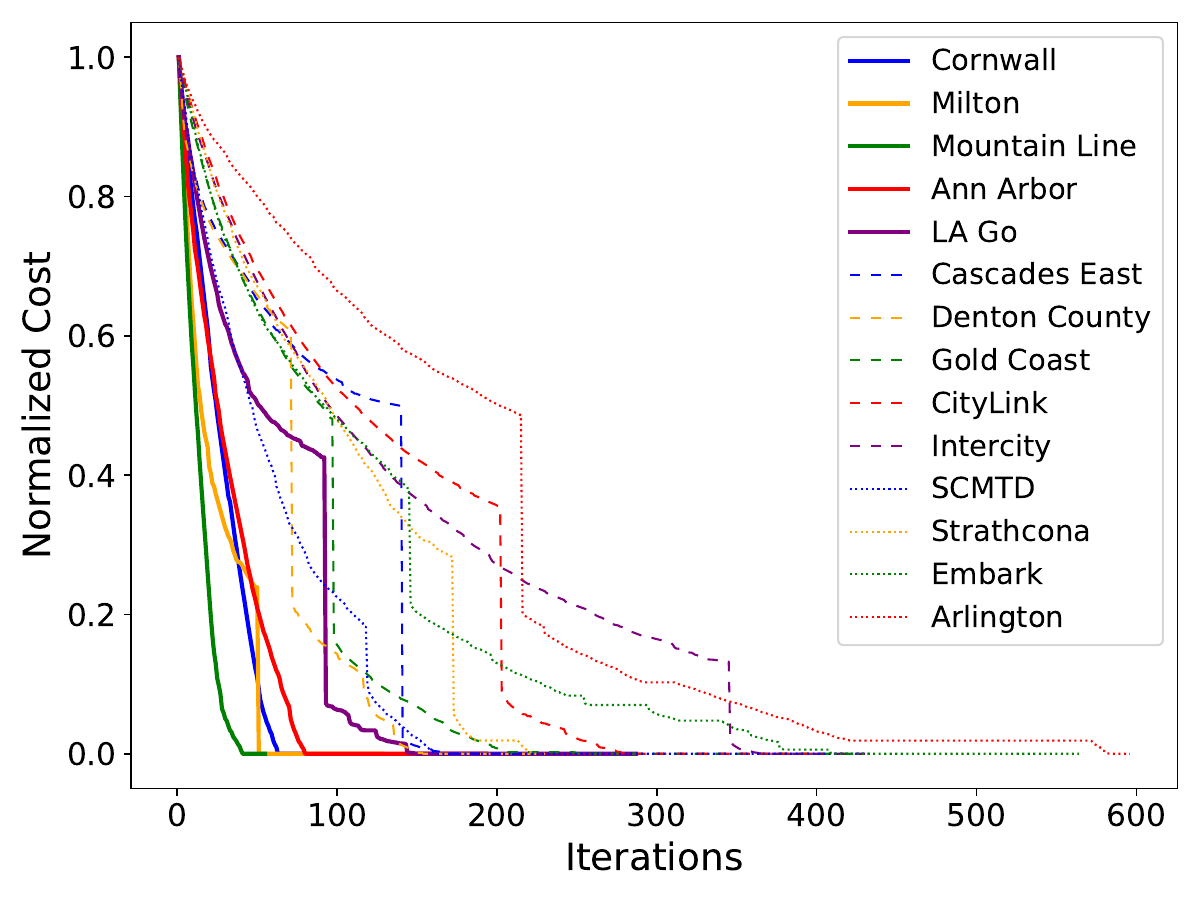}
    \caption{Joint model}
  \end{subfigure}
  \caption{Convergence of ILS models}
  \label{fig:convergence}
\end{figure}

\textbf{Convergence Analysis -- ILS Models:}
We examine the convergence patterns of ILS models using a normalized cost metric. The total cost of the CS model is normalized to $1$, and the normalized total cost of the sequential/joint model at termination is assumed to be $0$. In this analysis, CSP costs are excluded when studying the convergence of the sequential model, as the focus is exclusively on the CLP-EVSP. Convergence plots, as shown in Figure \ref{fig:convergence} for both sequential and joint models across all networks, reveal that smaller networks tend to converge faster than larger ones. Notably, most networks exhibit a sharp decline in normalized costs after a few iterations, primarily due to the closing of unused charging stations. The iteration count on the x-axis is assumed to increment each time a cost improvement occurs during calls to the \textsc{OptimizeRotations} function within Algorithm \ref{alg:fleet_optimization}.

\section{Conclusions and Future Work}
\label{sec:conc}
This study focused on a comprehensive approach to address the joint optimization of charging locations, bus-to-trip assignments, and charge scheduling. We first presented a sequential problem, which solves the CLP-EVSP, followed by charge schedule optimization using a new CSP MILP model. Additionally, we explored a fully integrated joint CLP-EVSP-CSP framework to better capture the interdependence between decision variables. Although MILP models can be formulated for both problems, they fall short of solving real-world instances. To overcome this, we employ a novel ILS framework for both models, with initial solutions generated using a concurrent scheduler algorithm. The joint model solves the CLP-EVSP-CSP, integrating the CSP during trip exchanges and shifts using innovative surrogate linear programming formulations. We applied the sequential and joint models to various real-world bus transit networks, with the joint model consistently outperforming the sequential model, achieving an average of 17.5\% and 14.1\% savings in scheduling and operational costs, respectively over all networks and maximum total cost savings of up to 7.7\%. 

This research sets the stage for combining learning-based approaches to identify parts of the state space worth exploring, thereby enabling the application of these techniques to larger transit networks. Potential areas for investigation could also include modifications to objectives and constraints based on battery life cycles, incorporation of realistic battery charging and discharging curves, addressing uncertainty in energy usage and travel time due to congestion, and enabling buses to take charging detours. Transit networks and schedules evolve over time; in this context, exploring multi-stage deterministic and stochastic models would be interesting. These models could account for changes in land use and travel patterns in urban settings, assessing their impact on long-term decisions regarding bus acquisition and charging infrastructure. 
Overall, our study makes a compelling case for joint modeling, suggesting that introducing additional features is worthwhile, provided that the trade-offs between solution quality and computation time are carefully managed.

\section*{Acknowledgments}
The authors thank the Indo-German Science and Technology Centre (IGSTC) for awarding the Paired Early Career Research Fellowship in Applied Research (PECFAR) which allowed us to undertake this collaborative work.

\section*{Author Contributions}
\textbf{Rito Brata Nath:} Methodology, Software, Writing -- Original Draft, Data Curation, Investigation, Visualization; \textbf{Tarun Rambha:} Conceptualization, Methodology, Software, Writing -- Original Draft, Investigation, Supervision, Funding acquisition; \textbf{Maximilian Schiffer:} Conceptualization, Methodology, Writing -- Review \& Editing, Supervision, Funding acquisition. 


\appendix
\renewcommand{\thesection}{Appendix \Alph{section}}
\section{Optimize Rotations}
\label{sec:appendix_optimize_rotations}

In this appendix, we detail the exchange and shift operators for service trips \textsc{Exst} and \textsc{Sst}, respectively. The sub-routine for exchanging depot trips, \textsc{Exd}, is similar to \textsc{Exst} and is hence not presented in detail. 

\subsection{Exchanging Service Trips}
Algorithm \ref{alg:exst} exchanges service trips between a pair of buses to reduce operational costs, primarily from deadheading. The algorithm computes cost savings for all feasible pairs based on trip compatibility (line 4) and charging levels (line 7). We finally perform the exchange yielding the highest cost saving (lines 10--12). 

\begin{figure}[h]
	\centering
	\includegraphics[scale=0.52]{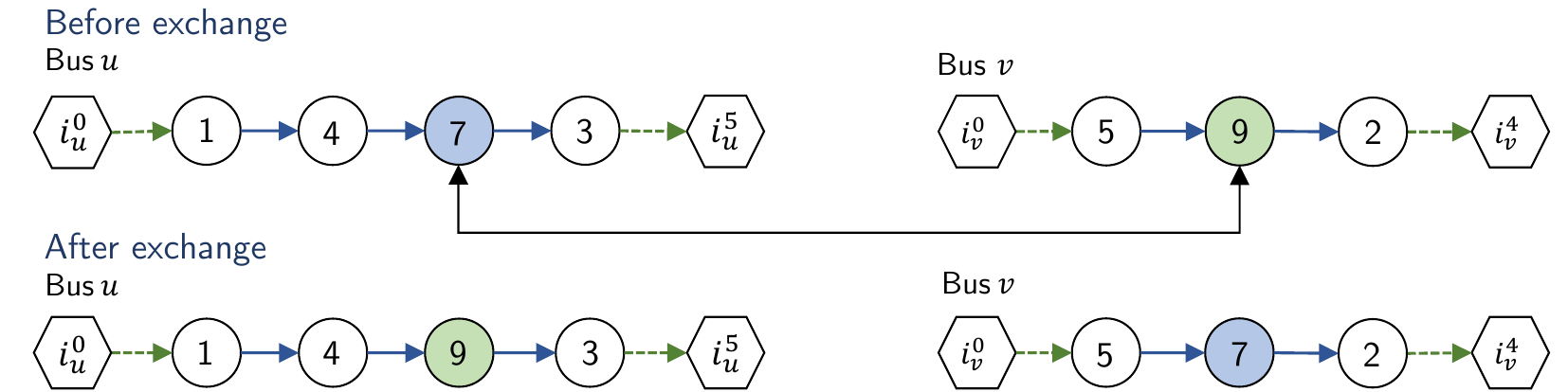}
	\caption{Example of exchanging service trips}
\label{fig:diagram_exst}
\end{figure}

For the joint CLP-EVSP-CSP model, this exchange sub-routine additionally calculates savings in electricity and power capacity costs using the charge scheduling LP (line 9). Figure \ref{fig:diagram_exst} illustrates a trip exchange step, where trips $7$ and $9$ are exchanged between two buses. Note that trip pairs ($4, 9$), ($9, 3$), ($5, 7$), and ($7, 2$) must be compatible, and the exchange should not make the rotations charge infeasible.

\begin{algorithm}[h]
\caption{\textsc{ExchangeServiceTrips (Exst)}}
\label{alg:exst}

\KwIn{$\busrotationList, \LocationSet$}
\KwOut{$\ExchangetriprotationList, \textsf{exchangeSavings}$}

$\ExchangetriprotationList \gets \busrotationList$ \;
$\textsf{exchangeSavings} \gets -\infty$\;
\tcp{\textsf{Iterate across all possible trip exchanges}}
\For{$\left\{ \{\busrotationList_u, \busrotationList_v\} : \busrotationList_u, \busrotationList_v \in \busrotationList, \busrotationList_u \neq \busrotationList_v \right\}$}{
    \For{$j \in \{1, 2, \ldots , \numTripsBus{u}\}, k \in \{ 1, 2, \ldots , \numTripsBus{v}\}: (\bustrip{u}{j-1}, \bustrip{v}{k}), (\bustrip{v}{k}, \bustrip{u}{j+1}), (\bustrip{v}{k-1}, \bustrip{u}{j}), (\bustrip{u}{j}, \bustrip{v}{k+1}) \in \arcSet$}{
        $\busrotationAnotherList \leftarrow \busrotationList$\;
        Update $\busrotationAnotherList$ by exchanging the $j^\text{th}$ trip of bus $u$ with $k^\text{th}$ trip of bus $v$\;
        \If{$\textsc{AreRotationsChargeFeasible}(\busrotationAnotherList, \LocationSet$)}{
            \If{the joint model is solved}{
                Solve the CSP LP (Section \ref{sec:surrogate_lp}) and add the charging costs to the $f$ values\;
            }
            \tcp{\textsf{Select the exchange that yields maximum savings}}
            \If{$\objectivefunction{\busrotationAnotherList}{\LocationSet} - \objectivefunction{\busrotationList}{\LocationSet} > \mathsf{exchangeSavings}$}{
                $\textsf{exchangeSavings} \gets \objectivefunction{\busrotationAnotherList}{\LocationSet} - \objectivefunction{\busrotationList}{\LocationSet}$\;
                $\ExchangetriprotationList \gets \busrotationAnotherList$\;
            }
        }
    }
}
\end{algorithm}

\subsection{Shifting Service Trips}

Algorithm \ref{alg:sst} performs an exhaustive search by shifting a service trip from one bus to another. We verify the possibility of trip insertion based on trip compatibility (line 4) and charge feasibility of the updated configuration (line 7). The cost savings are similarly to those for the exchange operator (lines 10--12), with the exception that if a bus performs no service trips, it is omitted, and the vehicle's fixed cost is recovered. As before, we solve the charge scheduling LP only for the joint CLP-EVSP-CSP model (line 9).

\begin{figure}[H]
	\centering
\includegraphics[scale=0.52]{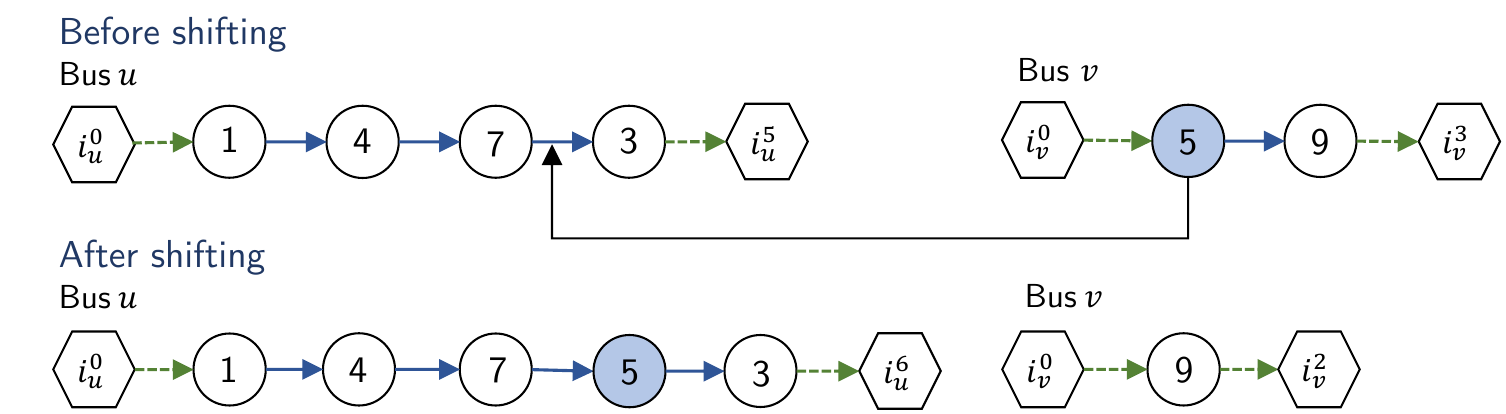}
	\caption{Example of shifting service trips}
\label{fig:diagram_sst}
\end{figure}

Figure \ref{fig:diagram_sst} illustrates an example of shifting of service trips. Suppose $(7, 5)$ and $(5, 3)$ are compatible trip pairs. Trip $5$, originally performed by bus $v$, is shifted after trip $7$ of bus $u$ after ensuring that the resulting rotation is charge feasible.

\begin{algorithm}[H]
\caption{\textsc{ShiftServiceTrips (Sst)}}
\label{alg:sst}

\KwIn{$\busrotationList, \LocationSet$}
\KwOut{$\ShifttriprotationList, \textsf{shiftSavings}$}

$\ShifttriprotationList \gets \busrotationList$ \;
$\textsf{shiftSavings} \gets -\infty$\;
\tcp{\textsf{Iterate across all possible shifts}}
\For{$\busrotationList_u, \busrotationList_v \in \busrotationList: \busrotationList_u \neq \busrotationList_v$}{
    \For{$j \in \{1, 2, \ldots , \numTripsBus{u}\}, k \in \{1, 2, \ldots , \numTripsBus{v}\}: (\bustrip{u}{j}, \bustrip{v}{k}), (\bustrip{v}{k}, \bustrip{u}{j+1}) \in \arcSet$}{
        $\busrotationAnotherList \gets \busrotationList$\;
        Update $\busrotationAnotherList$ by shifting the $k^\text{th}$ trip of bus $v$ after $j^\text{th}$ trip of bus $u$\;
        \If{$\textsc{AreRotationsChargeFeasible}(\busrotationAnotherList, \LocationSet)$}{
            \If{the joint model is solved}{
                Solve the CSP LP (Section \ref{sec:surrogate_lp}) and add the charging costs to the $f$ values\;
            }   
            \tcp{\textsf{Select the shift that yields maximum savings}}
            \If{$\objectivefunction{\busrotationAnotherList}{\LocationSet} - \objectivefunction{\busrotationList}{\LocationSet} > \mathsf{shiftSavings}$}{
                $\textsf{shiftSavings} \gets \objectivefunction{\busrotationAnotherList}{\LocationSet} - \objectivefunction{\busrotationList}{\LocationSet}$\;
                $\ShifttriprotationList \gets \busrotationAnotherList$\;
            }
        }
    }
}

\end{algorithm}

\section{ILS Routines for Shifting Multiple Trips}
\label{sec:appendix_multiple_shift}
Although we found that the rotations from the concurrent scheduler are usually tight, the closure of charging stations can lead to the splitting of rotations and the creation of rotations with fewer trips. In such cases, shifting one trip at a time in the $\textsc{Sst}$ subroutine may not help release these newly created rotations since it may be optimal to exchange and shift other trips in subsequent iterations. Since bus acquisition costs are one of the largest components of the overall objective, we use an operator that attempts to shift multiple trips in rotations that have fewer trips than a user-defined threshold. The locations at which these trips are to be inserted are chosen greedily. Even when they do not reduce the number of rotations, these operators offer some benefits to the deadheading costs while remaining computationally inexpensive.

\begin{algorithm}[h]
\caption{\textsc{OptimizeMultipleShifts}}
\label{alg:opt_multiple_shifts}

\KwIn{$\busrotationList, \LocationSet$}
\KwOut{$\busrotationList$}

$\textsf{improvement} \gets \infty$\;
$\bestbusrotationList \gets \busrotationList$\;

\tcp{\textsf{Apply shifts repeatedly until the benefits are exhausted}}
\While{$\mathsf{improvement} > 0$}{
    $\busrotationList \gets \bestbusrotationList$\;

    $\busrotationList^{few} \gets$ Rotations in $\busrotationList$ with fewer than $\fewThreshold$ service trips \;

    $\textsf{maxSavings} \gets -\infty$ \;

    \tcp{\textsf{Find the rotation that gives the maximum savings from shifting multiple trips}}
    \For{$\busrotationList_{b} \in \busrotationList^{few}$}{
    $\busrotationList^{msst}, \textsf{multiShiftSavings} \gets \textsc{ShiftMultipleTrips}(\busrotationList, \LocationSet, \busrotationList_{b})$

    \If{$\mathsf{multiShiftSavings} > \mathsf{maxSavings}$}{
    $\textsf{maxSavings} \gets \textsf{multiShiftSavings}$ \;
    $\busrotationList^{best} \gets \busrotationList^{msst}$ \;
    }
}
   
    $\textsf{improvement} \gets f(\busrotationList, \LocationSet) - f(\bestbusrotationList, \LocationSet)$\;
}
\end{algorithm}

Algorithm \ref{alg:opt_multiple_shifts} describes a sub-routine that iteratively checks if there is a rotation that improves the objective when multiple trips are shifted to other rotations. To keep these calculations tractable, we consider only those rotations with $\fewThreshold = 5$ or fewer service trips (see line 5). Each rotation in $\busrotationList^{few}$ is examined in the $\textsc{ShiftMultipleTrips}$ sub-routine (lines 7--8) for savings from shifting multiple trips, and the resulting rotations are saved in $\busrotationList^{msst}$. The most promising of these bus-to-trip assignments are saved in $\busrotationList^{best}$ in lines 9--11. We then calculate the change in the objective in line 12, and the process is repeated until such shifts are no longer advantageous. 

The \textsc{ShiftMultipleTrips} sub-routine (see Algorithm \ref{alg:shift_multiple_trips}) is similar to the \textsc{Sst} algorithm, except that it tries to shift multiple trips assigned to a bus to other rotations. The output includes an adjusted list of bus rotations ($\busrotationList^{msst}$) and the maximum savings possible from such greedy insertions ($\textsf{multiShiftSavings}$). The algorithm begins by creating a copy of the original bus rotations that is returned if no trip can be shifted (line 1). It then iterates over each trip in $\busrotationList_{b}$ and initializes variables to track the feasibility of a single shift and determines where to insert the displaced trip for maximum savings in the objective. Trip compatibility and charge feasibility are necessary during these calculations, and we embed the CSP linear program when the problem is jointly solved. Lines 4--13 carry out these steps and are similar to the \textsc{Sst} algorithm. If the shifted trip can be inserted elsewhere, the algorithm updates the rotation list, and proceeds to shift the next trip in $\busrotationList_b$. Note that some of the shifts can lead to negative savings. The overall objective of this procedure, however, can only improve as it is taken care of in line 3 of Algorithm \ref{alg:opt_multiple_shifts}.

\begin{algorithm}[H]
\caption{\textsc{ShiftMultipleTrips}}
\label{alg:shift_multiple_trips}

\KwIn{{$\busrotationList, \LocationSet, \busrotationList_{b}$}}
\KwOut{{$\busrotationList^{msst}, \textsf{multiShiftSavings}$}}

$\busrotationList^{org} \gets \busrotationList$ \;
\tcp{\textsf{Shift multiple trips in the rotation}}
\For{{$k=1$ \KwTo $\numTripsBus{b}$}}{
\textsf{multiShiftSavings} $\gets -\infty$\;
\tcp{\textsf{Iterate across all possible locations to insert the trip}}
    \For{{$\busrotationList_{u} \in \busrotationList : \busrotationList_{u} \neq \busrotationList_{b}$}}{
        \For{{$j \in \{1, 2, \ldots , \numTripsBus{u}\} : (\bustrip{u}{j}, \bustrip{b}{k}), (\bustrip{b}{k}, \bustrip{u}{j+1}) \in \arcSet$}}{
        $\busrotationAnotherList \gets \busrotationList$ \;
        Update $\busrotationAnotherList$ by shifting the $k^\text{th}$ trip of bus $b$ after $j^\text{th}$ trip of bus $u$\;
        \If{{$\textsc{AreRotationsChargeFeasible}(\busrotationAnotherList, \LocationSet)$}}{
            \If{the joint model is solved}{
                Solve the CSP LP (Section \ref{sec:surrogate_lp}) and add the charging costs to the $f$ values\;
            }     
            \tcp{\textsf{Select the shift that gives the maximum savings}}
            \If{{$\objectivefunction{\busrotationAnotherList}{\LocationSet} - \objectivefunction{\busrotationList}{\LocationSet} > \mathsf{multiShiftSavings}$}}{
                $\textsf{multiShiftSavings} \gets \objectivefunction{\busrotationAnotherList}{\LocationSet} - \objectivefunction{\busrotationList}{\LocationSet}$\;
                $\busrotationList^{msst} \gets \busrotationAnotherList$\;
            }
        }
        }
    }

    \tcp{\textsf{Update the rotations if shifting the trip is feasible}}
    \If{$\mathsf{multiShiftSavings = - \infty}$}{\Return $\busrotationList^{org}$}{
    $\busrotationList \gets \busrotationList^{msst}$\;
    }
}
\end{algorithm}

\bibliography{references}

\begin{thebibliography}{}

\bibitem[\protect\citeauthoryear{Abdelwahed, van~den Berg, Brandt, Collins, and
  Ketter}{Abdelwahed et~al.}{2020}]{abdelwahed2020evaluating}
Abdelwahed, A., P.~L. van~den Berg, T.~Brandt, J.~Collins, and W.~Ketter
  (2020).
\newblock Evaluating and optimizing opportunity fast-charging schedules in
  transit battery electric bus networks.
\newblock {\em Transportation Science\/}~{\em 54\/}(6), 1601--1615.

\bibitem[\protect\citeauthoryear{Adler and Mirchandani}{Adler and
  Mirchandani}{2017}]{adler2017vehicle}
Adler, J.~D. and P.~B. Mirchandani (2017).
\newblock The vehicle scheduling problem for fleets with alternative-fuel
  vehicles.
\newblock {\em Transportation Science\/}~{\em 51\/}(2), 441--456.

\bibitem[\protect\citeauthoryear{Ahuja, Ergun, Orlin, and Punnen}{Ahuja
  et~al.}{2002}]{ahuja2002survey}
Ahuja, R.~K., {\"O}.~Ergun, J.~B. Orlin, and A.~P. Punnen (2002).
\newblock A survey of very large-scale neighborhood search techniques.
\newblock {\em Discrete Applied Mathematics\/}~{\em 123\/}(1-3), 75--102.

\bibitem[\protect\citeauthoryear{Alvo, Angulo, and Klapp}{Alvo
  et~al.}{2021}]{alvo2021exact}
Alvo, M., G.~Angulo, and M.~A. Klapp (2021).
\newblock An exact solution approach for an electric bus dispatch problem.
\newblock {\em Transportation Research Part E: Logistics and Transportation
  Review\/}~{\em 156}, 102528.

\bibitem[\protect\citeauthoryear{Alwesabi, Liu, Kwon, and Wang}{Alwesabi
  et~al.}{2021}]{alwesabi2021novel}
Alwesabi, Y., Z.~Liu, S.~Kwon, and Y.~Wang (2021).
\newblock A novel integration of scheduling and dynamic wireless charging
  planning models of battery electric buses.
\newblock {\em Energy\/}~{\em 230}, 120806.

\bibitem[\protect\citeauthoryear{Alwesabi, Wang, Avalos, and Liu}{Alwesabi
  et~al.}{2020}]{alwesabi2020electric}
Alwesabi, Y., Y.~Wang, R.~Avalos, and Z.~Liu (2020).
\newblock Electric bus scheduling under single depot dynamic wireless charging
  infrastructure planning.
\newblock {\em Energy\/}~{\em 213}, 118855.

\bibitem[\protect\citeauthoryear{Bagherinezhad, Palomino, Li, and
  Parvania}{Bagherinezhad et~al.}{2020}]{bagherinezhad2020spatio}
Bagherinezhad, A., A.~D. Palomino, B.~Li, and M.~Parvania (2020).
\newblock Spatio-temporal electric bus charging optimization with transit
  network constraints.
\newblock {\em IEEE Transactions on Industry Applications\/}~{\em 56\/}(5),
  5741--5749.

\bibitem[\protect\citeauthoryear{Bertossi, Carraresi, and Gallo}{Bertossi
  et~al.}{1987}]{bertossi1987some}
Bertossi, A.~A., P.~Carraresi, and G.~Gallo (1987).
\newblock On some matching problems arising in vehicle scheduling models.
\newblock {\em Networks\/}~{\em 17\/}(3), 271--281.

\bibitem[\protect\citeauthoryear{Bodin, Rosenfield, and Kydes}{Bodin
  et~al.}{1978}]{bodin1978ucost}
Bodin, L., D.~Rosenfield, and A.~Kydes (1978).
\newblock {UCOST}: A micro approach to a transportation planning problem.
\newblock {\em Journal of Urban Analysis\/}~{\em 5\/}(1).

\bibitem[\protect\citeauthoryear{Bron and Kerbosch}{Bron and
  Kerbosch}{1973}]{bron1973algorithm}
Bron, C. and J.~Kerbosch (1973).
\newblock Algorithm 457: Finding all cliques of an undirected graph.
\newblock {\em Communications of the ACM\/}~{\em 16\/}(9), 575--577.

\bibitem[\protect\citeauthoryear{Bunte and Kliewer}{Bunte and
  Kliewer}{2009}]{bunte2009overview}
Bunte, S. and N.~Kliewer (2009).
\newblock An overview on vehicle scheduling models.
\newblock {\em Public Transport\/}~{\em 1\/}(4), 299--317.

\bibitem[\protect\citeauthoryear{BYD}{BYD}{2023}]{BYD}
BYD (2023).
\newblock {BYD} electric bus specifications.
\newblock \url{https://en.byd.com/bus/}.
\newblock Accessed: 2023-07-26.

\bibitem[\protect\citeauthoryear{Cavadas, de~Almeida~Correia, and
  Gouveia}{Cavadas et~al.}{2015}]{cavadas2015mip}
Cavadas, J., G.~H. de~Almeida~Correia, and J.~Gouveia (2015).
\newblock A {MIP} model for locating slow-charging stations for electric
  vehicles in urban areas accounting for driver tours.
\newblock {\em Transportation Research Part E: Logistics and Transportation
  Review\/}~{\em 75}, 188--201.

\bibitem[\protect\citeauthoryear{Dirks, Schiffer, and Walther}{Dirks
  et~al.}{2022}]{dirks2022integration}
Dirks, N., M.~Schiffer, and G.~Walther (2022).
\newblock On the integration of battery electric buses into urban bus networks.
\newblock {\em Transportation Research Part C: Emerging Technologies\/}~{\em
  139}, 103628.

\bibitem[\protect\citeauthoryear{Ellingsen, Thorne, Wind, Figenbaum, Romare,
  and Nordel{\"o}f}{Ellingsen et~al.}{2022}]{ellingsen2022life}
Ellingsen, L. A.-W., R.~J. Thorne, J.~Wind, E.~Figenbaum, M.~Romare, and
  A.~Nordel{\"o}f (2022).
\newblock Life cycle assessment of battery electric buses.
\newblock {\em Transportation Research Part D: Transport and
  Environment\/}~{\em 112}, 103498.

\bibitem[\protect\citeauthoryear{Forbes, Holt, and Watts}{Forbes
  et~al.}{1994}]{forbes1994exact}
Forbes, M., J.~Holt, and A.~Watts (1994).
\newblock An exact algorithm for multiple depot bus scheduling.
\newblock {\em European Journal of Operational Research\/}~{\em 72\/}(1),
  115--124.

\bibitem[\protect\citeauthoryear{Gairola and Nezamuddin}{Gairola and
  Nezamuddin}{2023}]{gairola2023optimization}
Gairola, P. and N.~Nezamuddin (2023).
\newblock Optimization framework for integrated battery electric bus planning
  and charging scheduling.
\newblock {\em Transportation Research Part D: Transport and
  Environment\/}~{\em 118}, 103697.

\bibitem[\protect\citeauthoryear{Gkiotsalitis, Iliopoulou, and
  Kepaptsoglou}{Gkiotsalitis et~al.}{2023}]{gkiotsalitis2023exact}
Gkiotsalitis, K., C.~Iliopoulou, and K.~Kepaptsoglou (2023).
\newblock An exact approach for the multi-depot electric bus scheduling problem
  with time windows.
\newblock {\em European Journal of Operational Research\/}~{\em 306\/}(1),
  189--206.

\bibitem[\protect\citeauthoryear{Haghani and Banihashemi}{Haghani and
  Banihashemi}{2002}]{haghani2002heuristic}
Haghani, A. and M.~Banihashemi (2002).
\newblock Heuristic approaches for solving large-scale bus transit vehicle
  scheduling problem with route time constraints.
\newblock {\em Transportation Research Part A: Policy and Practice\/}~{\em
  36\/}(4), 309--333.

\bibitem[\protect\citeauthoryear{Hansen, Mladenovi{\'c}, and
  Moreno~Perez}{Hansen et~al.}{2010}]{hansen2010variable}
Hansen, P., N.~Mladenovi{\'c}, and J.~A. Moreno~Perez (2010).
\newblock Variable neighbourhood search: methods and applications.
\newblock {\em Annals of Operations Research\/}~{\em 175}, 367--407.

\bibitem[\protect\citeauthoryear{He, Yin, and Zhou}{He
  et~al.}{2015}]{he2015deploying}
He, F., Y.~Yin, and J.~Zhou (2015).
\newblock Deploying public charging stations for electric vehicles on urban
  road networks.
\newblock {\em Transportation Research Part C: Emerging Technologies\/}~{\em
  60}, 227--240.

\bibitem[\protect\citeauthoryear{He, Liu, and Song}{He
  et~al.}{2020}]{he2020optimal}
He, Y., Z.~Liu, and Z.~Song (2020).
\newblock Optimal charging scheduling and management for a fast-charging
  battery electric bus system.
\newblock {\em Transportation Research Part E: Logistics and Transportation
  Review\/}~{\em 142}, 102056.

\bibitem[\protect\citeauthoryear{He, Song, and Liu}{He
  et~al.}{2019}]{he2019fast}
He, Y., Z.~Song, and Z.~Liu (2019).
\newblock Fast-charging station deployment for battery electric bus systems
  considering electricity demand charges.
\newblock {\em Sustainable Cities and Society\/}~{\em 48}, 101530.

\bibitem[\protect\citeauthoryear{Houbbadi, Trigui, Pelissier, Redondo-Iglesias,
  and Bouton}{Houbbadi et~al.}{2019}]{houbbadi2019optimal}
Houbbadi, A., R.~Trigui, S.~Pelissier, E.~Redondo-Iglesias, and T.~Bouton
  (2019).
\newblock Optimal scheduling to manage an electric bus fleet overnight
  charging.
\newblock {\em Energies\/}~{\em 12\/}(14), 2727.

\bibitem[\protect\citeauthoryear{Hu, Du, Liu, and Perez}{Hu
  et~al.}{2022}]{hu2022joint}
Hu, H., B.~Du, W.~Liu, and P.~Perez (2022).
\newblock A joint optimisation model for charger locating and electric bus
  charging scheduling considering opportunity fast charging and uncertainties.
\newblock {\em Transportation Research Part C: Emerging Technologies\/}~{\em
  141}, 103732.

\bibitem[\protect\citeauthoryear{Iliopoulou and Kepaptsoglou}{Iliopoulou and
  Kepaptsoglou}{2019}]{iliopoulou2019integrated}
Iliopoulou, C. and K.~Kepaptsoglou (2019).
\newblock Integrated transit route network design and infrastructure planning
  for on-line electric vehicles.
\newblock {\em Transportation Research Part D: Transport and
  Environment\/}~{\em 77}, 178--197.

\bibitem[\protect\citeauthoryear{Iliopoulou and Kepaptsoglou}{Iliopoulou and
  Kepaptsoglou}{2021}]{iliopoulou2021robust}
Iliopoulou, C. and K.~Kepaptsoglou (2021).
\newblock Robust electric transit route network design problem ({RE-TRNDP})
  with delay considerations: model and application.
\newblock {\em Transportation Research Part C: Emerging Technologies\/}~{\em
  129}, 103255.

\bibitem[\protect\citeauthoryear{Jahic, Eskander, and Schulz}{Jahic
  et~al.}{2019}]{jahic2019charging}
Jahic, A., M.~Eskander, and D.~Schulz (2019).
\newblock Charging schedule for load peak minimization on large-scale electric
  bus depots.
\newblock {\em Applied Sciences\/}~{\em 9\/}(9), 1748.

\bibitem[\protect\citeauthoryear{Janovec and Koh{\'a}ni}{Janovec and
  Koh{\'a}ni}{2019}]{janovec2019exact}
Janovec, M. and M.~Koh{\'a}ni (2019).
\newblock Exact approach to the electric bus fleet scheduling.
\newblock {\em Transportation Research Procedia\/}~{\em 40}, 1380--1387.

\bibitem[\protect\citeauthoryear{Jiang, Zhang, and Zhang}{Jiang
  et~al.}{2021}]{jiang2021multi}
Jiang, M., Y.~Zhang, and Y.~Zhang (2021).
\newblock Multi-depot electric bus scheduling considering operational
  constraint and partial charging: A case study in {S}henzhen, {C}hina.
\newblock {\em Sustainability\/}~{\em 14\/}(1), 255.

\bibitem[\protect\citeauthoryear{Ke, Lin, Chen, and Fang}{Ke
  et~al.}{2020}]{ke2020battery}
Ke, B.-R., Y.-H. Lin, H.-Z. Chen, and S.-C. Fang (2020).
\newblock Battery charging and discharging scheduling with demand response for
  an electric bus public transportation system.
\newblock {\em Sustainable Energy Technologies and Assessments\/}~{\em 40},
  100741.

\bibitem[\protect\citeauthoryear{Klein and Schiffer}{Klein and
  Schiffer}{2023}]{klein2023electric}
Klein, P.~S. and M.~Schiffer (2023).
\newblock Electric vehicle charge scheduling with flexible service operations.
\newblock {\em Transportation Science\/}~{\em 57\/}(6), 1605--1626.

\bibitem[\protect\citeauthoryear{Kliewer, Mellouli, and Suhl}{Kliewer
  et~al.}{2006}]{kliewer2006time}
Kliewer, N., T.~Mellouli, and L.~Suhl (2006).
\newblock A time--space network based exact optimization model for multi-depot
  bus scheduling.
\newblock {\em European Journal of Operational Research\/}~{\em 175\/}(3),
  1616--1627.

\bibitem[\protect\citeauthoryear{Kunith, Mendelevitch, and Goehlich}{Kunith
  et~al.}{2017}]{kunith2017electrification}
Kunith, A., R.~Mendelevitch, and D.~Goehlich (2017).
\newblock Electrification of a city bus network—{A}n optimization model for
  cost-effective placing of charging infrastructure and battery sizing of
  fast-charging electric bus systems.
\newblock {\em International Journal of Sustainable Transportation\/}~{\em
  11\/}(10), 707--720.

\bibitem[\protect\citeauthoryear{Lee, Shon, Papakonstantinou, and Son}{Lee
  et~al.}{2021}]{lee2021optimal}
Lee, J., H.~Shon, I.~Papakonstantinou, and S.~Son (2021).
\newblock Optimal fleet, battery, and charging infrastructure planning for
  reliable electric bus operations.
\newblock {\em Transportation Research Part D: Transport and
  Environment\/}~{\em 100}, 103066.

\bibitem[\protect\citeauthoryear{Leou and Hung}{Leou and
  Hung}{2017}]{leou2017optimal}
Leou, R.-C. and J.-J. Hung (2017).
\newblock Optimal charging schedule planning and economic analysis for electric
  bus charging stations.
\newblock {\em Energies\/}~{\em 10\/}(4), 483.

\bibitem[\protect\citeauthoryear{Li}{Li}{2014}]{li2014transit}
Li, J.-Q. (2014).
\newblock Transit bus scheduling with limited energy.
\newblock {\em Transportation Science\/}~{\em 48\/}(4), 521--539.

\bibitem[\protect\citeauthoryear{Li, Lo, and Xiao}{Li
  et~al.}{2019}]{li2019mixed}
Li, L., H.~K. Lo, and F.~Xiao (2019).
\newblock Mixed bus fleet scheduling under range and refueling constraints.
\newblock {\em Transportation Research Part C: Emerging Technologies\/}~{\em
  104}, 443--462.

\bibitem[\protect\citeauthoryear{Li, Wang, Li, Feng, Wang, and Cheng}{Li
  et~al.}{2020}]{li2020joint}
Li, X., T.~Wang, L.~Li, F.~Feng, W.~Wang, and C.~Cheng (2020).
\newblock Joint optimization of regular charging electric bus transit network
  schedule and stationary charger deployment considering partial charging
  policy and time-of-use electricity prices.
\newblock {\em Journal of Advanced Transportation\/}~{\em 2020}.

\bibitem[\protect\citeauthoryear{Liu and Ceder}{Liu and
  Ceder}{2020}]{liu2020battery}
Liu, T. and A.~A. Ceder (2020).
\newblock Battery-electric transit vehicle scheduling with optimal number of
  stationary chargers.
\newblock {\em Transportation Research Part C: Emerging Technologies\/}~{\em
  114}, 118--139.

\bibitem[\protect\citeauthoryear{Liu, Yao, Lu, and Yuan}{Liu
  et~al.}{2019}]{liu2019regional}
Liu, Y., E.~Yao, M.~Lu, and L.~Yuan (2019).
\newblock Regional electric bus driving plan optimization algorithm considering
  charging time window.
\newblock {\em Mathematical Problems in Engineering\/}~{\em 2019}, 1--9.

\bibitem[\protect\citeauthoryear{Metais, Jouini, Perez, Berrada, and
  Suomalainen}{Metais et~al.}{2022}]{metais2022too}
Metais, M.-O., O.~Jouini, Y.~Perez, J.~Berrada, and E.~Suomalainen (2022).
\newblock Too much or not enough? planning electric vehicle charging
  infrastructure: A review of modeling options.
\newblock {\em Renewable and Sustainable Energy Reviews\/}~{\em 153}, 111719.

\bibitem[\protect\citeauthoryear{Mladenovi{\'c} and Hansen}{Mladenovi{\'c} and
  Hansen}{1997}]{mladenovic1997variable}
Mladenovi{\'c}, N. and P.~Hansen (1997).
\newblock Variable neighborhood search.
\newblock {\em Computers \& Operations Research\/}~{\em 24\/}(11), 1097--1100.

\bibitem[\protect\citeauthoryear{Olsen and Kliewer}{Olsen and
  Kliewer}{2020}]{olsen2020scheduling}
Olsen, N. and N.~Kliewer (2020).
\newblock Scheduling electric buses in public transport: modeling of the
  charging process and analysis of assumptions.
\newblock {\em Logistics Research\/}~{\em 13\/}(1), 4.

\bibitem[\protect\citeauthoryear{Olsen and Kliewer}{Olsen and
  Kliewer}{2022}]{olsen2022location}
Olsen, N. and N.~Kliewer (2022).
\newblock Location planning of charging stations for electric buses in public
  transport considering vehicle scheduling: A variable neighborhood search
  based approach.
\newblock {\em Applied Sciences\/}~{\em 12\/}(8), 3855.

\bibitem[\protect\citeauthoryear{Pelletier, Jabali, Mendoza, and
  Laporte}{Pelletier et~al.}{2019}]{pelletier2019electric}
Pelletier, S., O.~Jabali, J.~E. Mendoza, and G.~Laporte (2019).
\newblock The electric bus fleet transition problem.
\newblock {\em Transportation Research Part C: Emerging Technologies\/}~{\em
  109}, 174--193.

\bibitem[\protect\citeauthoryear{Perumal, Dollevoet, Huisman, Lusby, Larsen,
  and Riis}{Perumal et~al.}{2021}]{perumal2021solution}
Perumal, S.~S., T.~Dollevoet, D.~Huisman, R.~M. Lusby, J.~Larsen, and M.~Riis
  (2021).
\newblock Solution approaches for integrated vehicle and crew scheduling with
  electric buses.
\newblock {\em Computers \& Operations Research\/}~{\em 132}, 105268.

\bibitem[\protect\citeauthoryear{Perumal, Lusby, and Larsen}{Perumal
  et~al.}{2022}]{perumal2022electric}
Perumal, S.~S., R.~M. Lusby, and J.~Larsen (2022).
\newblock Electric bus planning \& scheduling: A review of related problems and
  methodologies.
\newblock {\em European Journal of Operational Research\/}~{\em 301\/}(2),
  395--413.

\bibitem[\protect\citeauthoryear{PGEC}{PGEC}{2023}]{electricschedule}
PGEC (2023).
\newblock {P}acific {G}as and {E}lectric {C}ompany: {E}lectric schedule {BEV}.
\newblock
  \url{https://www.pge.com/tariffs/assets/pdf/tariffbook/ELEC_SCHEDS_BEV.pdf}.
\newblock Accessed: 2024-01-10.

\bibitem[\protect\citeauthoryear{Reuer, Kliewer, and Wolbeck}{Reuer
  et~al.}{2015}]{reuer2015electric}
Reuer, J., N.~Kliewer, and L.~Wolbeck (2015).
\newblock The electric vehicle scheduling problem: A study on time-space
  network based and heuristic solution.
\newblock In {\em Proceedings of the Conference on Advanced Systems in Public
  Transport (CASPT)}.

\bibitem[\protect\citeauthoryear{Ribeiro and Soumis}{Ribeiro and
  Soumis}{1994}]{ribeiro1994column}
Ribeiro, C.~C. and F.~Soumis (1994).
\newblock A column generation approach to the multiple-depot vehicle scheduling
  problem.
\newblock {\em Operations Research\/}~{\em 42\/}(1), 41--52.

\bibitem[\protect\citeauthoryear{Rogge, Van~der Hurk, Larsen, and Sauer}{Rogge
  et~al.}{2018}]{rogge2018electric}
Rogge, M., E.~Van~der Hurk, A.~Larsen, and D.~U. Sauer (2018).
\newblock Electric bus fleet size and mix problem with optimization of charging
  infrastructure.
\newblock {\em Applied Energy\/}~{\em 211}, 282--295.

\bibitem[\protect\citeauthoryear{Sadati, Moshtagh, Shafie-khah, Rastgou, and
  Catal{\~a}o}{Sadati et~al.}{2019}]{sadati2019operational}
Sadati, S. M.~B., J.~Moshtagh, M.~Shafie-khah, A.~Rastgou, and J.~P.
  Catal{\~a}o (2019).
\newblock Operational scheduling of a smart distribution system considering
  electric vehicles parking lot: A bi-level approach.
\newblock {\em International Journal of Electrical Power \& Energy
  Systems\/}~{\em 105}, 159--178.

\bibitem[\protect\citeauthoryear{Sadeghi-Barzani, Rajabi-Ghahnavieh, and
  Kazemi-Karegar}{Sadeghi-Barzani et~al.}{2014}]{sadeghi2014optimal}
Sadeghi-Barzani, P., A.~Rajabi-Ghahnavieh, and H.~Kazemi-Karegar (2014).
\newblock Optimal fast charging station placing and sizing.
\newblock {\em Applied Energy\/}~{\em 125}, 289--299.

\bibitem[\protect\citeauthoryear{Sadeghian, Oshnoei, Mohammadi-Ivatloo,
  Vahidinasab, and Anvari-Moghaddam}{Sadeghian
  et~al.}{2022}]{sadeghian2022comprehensive}
Sadeghian, O., A.~Oshnoei, B.~Mohammadi-Ivatloo, V.~Vahidinasab, and
  A.~Anvari-Moghaddam (2022).
\newblock A comprehensive review on electric vehicles smart charging:
  solutions, strategies, technologies, and challenges.
\newblock {\em Journal of Energy Storage\/}~{\em 54}, 105241.

\bibitem[\protect\citeauthoryear{Schettini, Dell’Amico, Fumero, Jabali, and
  Malucelli}{Schettini et~al.}{2023}]{schettini2023locating}
Schettini, T., M.~Dell’Amico, F.~Fumero, O.~Jabali, and F.~Malucelli (2023).
\newblock Locating and sizing electric vehicle chargers considering multiple
  technologies.
\newblock {\em Energies\/}~{\em 16\/}(10), 4186.

\bibitem[\protect\citeauthoryear{Schiffer, Klein, Laporte, and
  Walther}{Schiffer et~al.}{2021}]{schiffer2021integrated}
Schiffer, M., P.~S. Klein, G.~Laporte, and G.~Walther (2021).
\newblock Integrated planning for electric commercial vehicle fleets: A case
  study for retail mid-haul logistics networks.
\newblock {\em European Journal of Operational Research\/}~{\em 291\/}(3),
  944--960.

\bibitem[\protect\citeauthoryear{Schiffer, Schneider, Walther, and
  Laporte}{Schiffer et~al.}{2019}]{schiffer2019vehicle}
Schiffer, M., M.~Schneider, G.~Walther, and G.~Laporte (2019).
\newblock Vehicle routing and location routing with intermediate stops: A
  review.
\newblock {\em Transportation Science\/}~{\em 53\/}(2), 319--343.

\bibitem[\protect\citeauthoryear{Shahraki, Cai, Turkay, and Xu}{Shahraki
  et~al.}{2015}]{shahraki2015optimal}
Shahraki, N., H.~Cai, M.~Turkay, and M.~Xu (2015).
\newblock Optimal locations of electric public charging stations using real
  world vehicle travel patterns.
\newblock {\em Transportation Research Part D: Transport and
  Environment\/}~{\em 41}, 165--176.

\bibitem[\protect\citeauthoryear{Sierzchula, Bakker, Maat, and
  Van~Wee}{Sierzchula et~al.}{2014}]{sierzchula2014influence}
Sierzchula, W., S.~Bakker, K.~Maat, and B.~Van~Wee (2014).
\newblock The influence of financial incentives and other socio-economic
  factors on electric vehicle adoption.
\newblock {\em Energy Policy\/}~{\em 68}, 183--194.

\bibitem[\protect\citeauthoryear{Stumpe, R{\"o}{\ss}ler, Schryen, and
  Kliewer}{Stumpe et~al.}{2021}]{stumpe2021study}
Stumpe, M., D.~R{\"o}{\ss}ler, G.~Schryen, and N.~Kliewer (2021).
\newblock Study on sensitivity of electric bus systems under simultaneous
  optimization of charging infrastructure and vehicle schedules.
\newblock {\em EURO Journal on Transportation and Logistics\/}~{\em 10},
  100049.

\bibitem[\protect\citeauthoryear{Tang, Shi, and Liu}{Tang
  et~al.}{2023}]{tang2023optimization}
Tang, C., H.~Shi, and T.~Liu (2023).
\newblock Optimization of single-line electric bus scheduling with skip-stop
  operation.
\newblock {\em Transportation Research Part D: Transport and
  Environment\/}~{\em 117}, 103652.

\bibitem[\protect\citeauthoryear{Tang, Lin, and He}{Tang
  et~al.}{2019}]{tang2019robust}
Tang, X., X.~Lin, and F.~He (2019).
\newblock Robust scheduling strategies of electric buses under stochastic
  traffic conditions.
\newblock {\em Transportation Research Part C: Emerging Technologies\/}~{\em
  105}, 163--182.

\bibitem[\protect\citeauthoryear{Teng, Chen, and Fan}{Teng
  et~al.}{2020}]{teng2020integrated}
Teng, J., T.~Chen, and W.~Fan (2020).
\newblock Integrated approach to vehicle scheduling and bus timetabling for an
  electric bus line.
\newblock {\em Journal of Transportation Engineering, Part A: Systems\/}~{\em
  146\/}(2), 04019073.

\bibitem[\protect\citeauthoryear{Wan, Sperling, and Wang}{Wan
  et~al.}{2015}]{wan2015china}
Wan, Z., D.~Sperling, and Y.~Wang (2015).
\newblock China’s electric car frustrations.
\newblock {\em Transportation Research Part D: Transport and
  Environment\/}~{\em 34}, 116--121.

\bibitem[\protect\citeauthoryear{Wang and Shen}{Wang and
  Shen}{2007}]{wang2007heuristic}
Wang, H. and J.~Shen (2007).
\newblock Heuristic approaches for solving transit vehicle scheduling problem
  with route and fueling time constraints.
\newblock {\em Applied Mathematics and Computation\/}~{\em 190\/}(2),
  1237--1249.

\bibitem[\protect\citeauthoryear{Wang and Wang}{Wang and
  Wang}{2010}]{wang2010locating}
Wang, Y.-W. and C.-R. Wang (2010).
\newblock Locating passenger vehicle refueling stations.
\newblock {\em Transportation Research Part E: Logistics and Transportation
  Review\/}~{\em 46\/}(5), 791--801.

\bibitem[\protect\citeauthoryear{Wen, Linde, Ropke, Mirchandani, and
  Larsen}{Wen et~al.}{2016}]{wen2016adaptive}
Wen, M., E.~Linde, S.~Ropke, P.~Mirchandani, and A.~Larsen (2016).
\newblock An adaptive large neighborhood search heuristic for the electric
  vehicle scheduling problem.
\newblock {\em Computers \& Operations Research\/}~{\em 76}, 73--83.

\bibitem[\protect\citeauthoryear{Xu, Miao, Zhang, and Shi}{Xu
  et~al.}{2013}]{xu2013optimal}
Xu, H., S.~Miao, C.~Zhang, and D.~Shi (2013).
\newblock Optimal placement of charging infrastructures for large-scale
  integration of pure electric vehicles into grid.
\newblock {\em International Journal of Electrical Power \& Energy
  Systems\/}~{\em 53}, 159--165.

\bibitem[\protect\citeauthoryear{Xylia, Leduc, Patrizio, Kraxner, and
  Silveira}{Xylia et~al.}{2017}]{xylia2017locating}
Xylia, M., S.~Leduc, P.~Patrizio, F.~Kraxner, and S.~Silveira (2017).
\newblock Locating charging infrastructure for electric buses in {S}tockholm.
\newblock {\em Transportation Research Part C: Emerging Technologies\/}~{\em
  78}, 183--200.

\bibitem[\protect\citeauthoryear{Yang, Lou, Yao, and Xie}{Yang
  et~al.}{2017}]{yang2017charging}
Yang, C., W.~Lou, J.~Yao, and S.~Xie (2017).
\newblock On charging scheduling optimization for a wirelessly charged electric
  bus system.
\newblock {\em IEEE Transactions on Intelligent Transportation Systems\/}~{\em
  19\/}(6), 1814--1826.

\bibitem[\protect\citeauthoryear{Yao, Liu, Lu, and Yang}{Yao
  et~al.}{2020}]{yao2020optimization}
Yao, E., T.~Liu, T.~Lu, and Y.~Yang (2020).
\newblock Optimization of electric vehicle scheduling with multiple vehicle
  types in public transport.
\newblock {\em Sustainable Cities and Society\/}~{\em 52}, 101862.

\bibitem[\protect\citeauthoryear{Zhang, Moura, Hu, and Song}{Zhang
  et~al.}{2016}]{zhang2016pev}
Zhang, H., S.~J. Moura, Z.~Hu, and Y.~Song (2016).
\newblock {PEV} fast-charging station siting and sizing on coupled
  transportation and power networks.
\newblock {\em IEEE Transactions on Smart Grid\/}~{\em 9\/}(4), 2595--2605.

\bibitem[\protect\citeauthoryear{Zhang, Liu, Wang, and Yu}{Zhang
  et~al.}{2022}]{zhang2022long}
Zhang, L., Z.~Liu, W.~Wang, and B.~Yu (2022).
\newblock Long-term charging infrastructure deployment and bus fleet transition
  considering seasonal differences.
\newblock {\em Transportation Research Part D: Transport and
  Environment\/}~{\em 111}, 103429.

\bibitem[\protect\citeauthoryear{Zhou, Xie, Zhao, and Lu}{Zhou
  et~al.}{2020}]{zhou2020collaborative}
Zhou, G.-J., D.-F. Xie, X.-M. Zhao, and C.~Lu (2020).
\newblock Collaborative optimization of vehicle and charging scheduling for a
  bus fleet mixed with electric and traditional buses.
\newblock {\em IEEE Access\/}~{\em 8}, 8056--8072.

\bibitem[\protect\citeauthoryear{Zhou, Liu, Wei, and Golub}{Zhou
  et~al.}{2020}]{zhou2020bi}
Zhou, Y., X.~C. Liu, R.~Wei, and A.~Golub (2020).
\newblock Bi-objective optimization for battery electric bus deployment
  considering cost and environmental equity.
\newblock {\em IEEE Transactions on Intelligent Transportation Systems\/}~{\em
  22\/}(4), 2487--2497.

\bibitem[\protect\citeauthoryear{Zhou, Meng, and Ong}{Zhou
  et~al.}{2022}]{zhou2022electric}
Zhou, Y., Q.~Meng, and G.~P. Ong (2022).
\newblock Electric bus charging scheduling for a single public transport route
  considering nonlinear charging profile and battery degradation effect.
\newblock {\em Transportation Research Part B: Methodological\/}~{\em 159},
  49--75.

\end{thebibliography}

\appendix

\end{document}